\documentclass[12pt]{amsart}
\usepackage{amssymb,epsfig,amsmath,latexsym,amsthm}
\usepackage{graphicx}
\usepackage{color}
\theoremstyle{plain}
\newtheorem{theorem}{Theorem}[section]
\newtheorem{lemma}[theorem]{Lemma}

\newtheorem{proposition}[theorem]{Proposition}
\newtheorem{corollary}[theorem]{Corollary}
\newtheorem{question}[theorem]{Question}

\newtheorem{conjecture}[theorem]{Conjecture}
\newtheorem{problem}[theorem]{Problem}
\theoremstyle{definition}
\newtheorem{definition}[theorem]{Definition}
\newtheorem{definition/construction}[theorem]{Definition/Construction}
\newtheorem{construction}[theorem]{Construction}
\newtheorem{remark}[theorem]{Remark}

\newtheorem{acknowledgements}[theorem]{Acknowledgements}

\newtheorem{remarks}[theorem]{Remarks}
\newtheorem{data memo}[theorem]{Data Memo}

\newtheorem{notation}[theorem]{Notation}

\newtheorem{operation}[theorem]{Operation}
\DeclareMathOperator{\Image}{Im}
\DeclareMathOperator{\Diff}{Diff}
\DeclareMathOperator{\Dim}{Dim}

\DeclareMathOperator{\Emb}{Emb}

\DeclareMathOperator{\id}{id}

\DeclareMathOperator{\Diam}{Diam}
\DeclareMathOperator{\diam}{diam}

\DeclareMathOperator{\genus}{genus}

\DeclareMathOperator{\inte}{int}

\DeclareMathOperator{\fix}{fix}
\DeclareMathOperator{\std}{std}

\DeclareMathOperator{\tw}{tw}
\DeclareMathOperator{\pt}{pt}
\DeclareMathOperator{\Max}{Max}

\newcommand{\solidtorus}{D^2\times S^1}

\newcommand{\ginv}{g^{-1}}

\newcommand{\pinv}{p^{-1}}
\newcommand{\piinv}{\pi^{-1}}
\newcommand{\qinv}{q^{-1}}

\newcommand{\phiinv}{\phi^{-1}}

\newcommand{\Psiinv}{\Psi^{-1}}

\newcommand{\BN}{\mathbb N}

\newcommand{\BQ}{\mathbb Q}
\newcommand{\BR}{\mathbb R}
\newcommand{\BS}{\mathbb S}
\newcommand{\BZ}{\mathbb Z}

\newcommand{\mRs}{\mR^{\std}}
\newcommand{\mGs}{\mG^{\std}}
\newcommand{\Rs}{R^{\std}}
\newcommand{\Gs}{G^{\std}}

\newcommand{\mA}{\mathcal{A}}
\newcommand{\mB}{\mathcal{B}}
\newcommand{\mC}{\mathcal{C}}
\newcommand{\mD}{\mathcal{D}}
\newcommand{\mE}{\mathcal{E}}
\newcommand{\mF}{\mathcal{F}}
\newcommand{\mG}{\mathcal{G}}
\newcommand{\mH}{\mathcal{H}}

\newcommand{\mL}{\mathcal{L}}

\newcommand{\mN}{\mathcal{N}}
\newcommand{\mP}{\mathcal{P}}

\newcommand{\mR}{\mathcal{R}}
\newcommand{\mS}{\mathcal{S}}

\newcommand{\mW}{\mathcal{W}}

\newcommand{\sonesthree}{S^1\times S^3}
\newcommand{\rsthree}{\BR\times S^3}
\newcommand{\stwostwo}{S^2\times S^2}

\begin{document}

\title{3-spheres in the 4-sphere and pseudo-isotopies of $S^1\times S^3$}

\author{David Gabai}
\address{Department of Mathematics\\Princeton
University\\Princeton, NJ 08544}
\email{gabai@math.princeton.edu}

\thanks{Version 0.695 May 29, 2024.\newline Partially supported by NSF grants
DMS-1006553, DMS-1607374, DMS-2003892 and DMS-2304841.
\newline\noindent\emph{Primary class:} 57M99
\newline\noindent\emph{secondary class:} 57R52, 57R50, 57N50
\newline\noindent\emph{keywords:} 4-manifolds, Schoenflies, pseudo-isotopy}

\begin{abstract}  We offer an approach to the smooth 4-dimensional Schoenflies conjecture via pseudo-isotopy theory.  Along the way we show that a diffeomorphism $\phi$ of a compact oriented 4-manifold is stably isotopic to $\id$ if and only if there is a pseudo-isotopy $f$ from $\id$ to $\phi$ such that $\Sigma(f)=0$ where $\Sigma$ is the first Hatcher - Wagoner pseudo-isotopy obstruction. \end{abstract}

\maketitle

\setcounter{section}{-1}

\section{Introduction}\label{S0}

The 4-dimensional smooth Schoenflies conjecture (SS4) asserts that every embedded smooth 3-sphere in the 4-sphere bounds a smooth 4-ball.  In 1958 Barry Mazur \cite{Ma1} proved that such spheres bound topological 4-balls, as a special case of a more general result.  A consequence, using \cite{Ce2}, of his elementary but strikingly original proof is the following that was known to Mazur and to topologists in the 1960's.

\begin{theorem}\label{mazur} SS4 is false if and only if there exists a diffeomorphism $\phi:S^1\times S^3\to S^1\times S^3$ such that $\phi$ is homotopic to $\id$ but $\phi(x_0\times S^3)$ is not isotopic to $x_0\times S^3$, even after lifting to any finite sheeted covering of $S^1\times S^3$.\end{theorem}

On the other hand it was recently proved in \cite{BG} and in \cite{Wa} that there exists diffeomorphisms $\phi\in \Diff_0(S^1\times S^3)$ with $\phi(x_0\times S^3)$  not isotopic to $x_0\times S^3$.  Here $\Diff_0(\sonesthree)$ denotes the topological group of diffeomorphisms of $S^1\times S^3$ whose elements are homotopic to $\id$, a  subgroup of index-8 in $\Diff(\sonesthree)$.  

This paper investigates $\Diff_0(S^1\times S^3)$ via pseudo-isotopy theory in an attempt to address the Schoenflies conjecture. Theorem \ref{main sequivalent} and Proposition \ref{whitney homotopy} equate SS4 to a certain interpolation problem involving a regular homotopy whose finger and Whitney discs coincide near their boundaries.  We also use this to give a condition for showing that a Poincare 4-ball (i.e. a compact contractible 4-manifold with boundary $S^3$) is a Schoenflies 4-ball (i.e. a closed complementary region of an embedded 3-sphere in $S^4$).  Our methods yield the following applications.

\begin{theorem} \label{diffs stable} If $\phi\in \Diff_0(S^1\times S^3)$,  then $\phi$ is stably isotopic to $\id$.  \end{theorem}

This means that for $k$ sufficiently large, we can isotope $\phi$ to be the $\id$ on a $\sqcup_k B^4$, and then extend to $\hat\phi\in\Diff_0( S^1\times S^3\#_k S^2\times S^2$) where the sums are taken inside the $B^4$'s and $\hat\phi$ is $\id$ within the summands so that $\hat\phi$ is isotopic to $\id$.  Recall that the first Hatcher - Wagoner \cite{HW} pseudo-isotopy obstruction $\Sigma$, is roughly the obstruction to having a nested eye 1-parameter family with no handle slides.  More generally we have:

\begin{theorem} \label{stable isotopy} Let $\phi\in \Diff_0(M)$, where $M$ is a compact oriented 4-manifold.  Then $\id$ is stably isotopic to $\phi$ if and only if $\id$ is pseudo-isotopic to $\phi$ by a pseudo-isotopy $f$ with $\Sigma(f)=0$. \end{theorem} 

\begin{remarks} i) The proof shows that the minimal \emph{stabilization number} is equal to the \emph{nested eye number}, i.e. the minimal number of stabilizations needed is equal to the minimal number of nested eyes required. See Corollary \ref{stabilization nested}.  

ii) For $M$ simply connected, a diffeomorphism homotopic to $\id$ is pseudo-isotopic to id, \cite{Kr}, \cite{Qu} and the pseudo-isotopy itself is stably isotopic to $\id$   \cite{Qu}, \cite{GGHKP}.\end{remarks}

We prove the following characterization of Schoenflies balls.

\begin{theorem} Every Schoenflies ball has a carving/surgery presentation.\end{theorem}

This means that it is obtained by a finite process starting with the 4-ball, attaching finitely many  2-handles, then carving finitely many 2-handles, then attaching finitely many 2-handles, etc., with every step happening in the 4-sphere.  An attached or carved 2-handle may nest a previously carved or attached 2-handle and so on.  Actually we show that the presentation can be chosen to be of a special type called an \emph{optimized $F|W$-carving/surgery presentation}.  See Definitions \ref{fwcs} and \ref{optimized} and Theorem \ref{optimized theorem}.  A key feature of a $F|W$-carving/surgery presentation is that when viewed as a surgery presentation of a 3-manifold $M\subset S^4$, our $M$ is \emph{obviously} the 3-sphere.  On the other hand, there are many compact 4-manifolds in $S^4$ with carving/surgery presentations.    Section 9 gives the definition of an optimized $F|W$-carving/surgery presentation and can be read independently of the rest of the paper.  
\vskip 10pt

\noindent\textbf{Underlying Coventions}: Unless said otherwise, this paper works in the smooth category, manifolds are orientable and diffeormorphisms between oriented manifolds are orientation preserving.

\begin{acknowledgements}  We thank Valentin Poenaru for many conversations about mathematics over the last 30 years and for introducing me to Barry Mazur's work.  We thank Bob Edwards for long conversations over the last 15 years.  We thank Toby Colding for discussions leading to the rediscovery of Mazur's unpublished theorem.   During 2016 part of this work, namely the reduction of the Schoenflies problem to a certain interpolation problem Theorem \ref{germs theorem}, was presented at the  Bonn Max Planck Institute ``Conference on 4-manifolds and knot concordance", the ``Benjamin Pierce Centennial Conference" at Harvard and at a meeting at Trinity College.  In the interim have appeared the papers \cite{Ga1}, \cite{ST}, \cite{BG} and \cite{Gay} which have enabled a more concise exposition of that theorem.  We thank Ryan Budney, David Gay and Hannah Schwartz for helpful conversations.  This research was  carried out in part during visits to Trinity College, Dublin, the Institute for Advanced Study and the IHES.  We thank these institutions for their hospitality. \end{acknowledgements}

\section {The Schoenflies conjecture and $\Diff_0(\sonesthree)$} \label{test}

Before proving Mazur's Theorem \ref{mazur} we record the following basic fact. 

\begin{lemma} \label{parallel} Let $E,W\in S^4$ be distinct points. Two oriented 3-spheres $\Sigma_0, \Sigma_1\subset S^4\setminus \{E,W\}$ are isotopic in $S^4\setminus \{E,W\}$ if and only if they are isotopic in $S^4$ and they represent the same class in $H_3(S^4\setminus \{E,W\})$.  \qed\end{lemma}

\noindent\emph{Proof of Theorem \ref{mazur}}: View $S^4$ as $[-\infty,\infty]\times S^3$ where $-\infty\times S^3$ and $\infty\times S^3$ are identified to points, which are respectively denoted $W$ and $E$.  The $\BZ$-action 
$t\times S^3\to t+1\times S^3$ induces the covering $\pi:\BR\times S^3\to S^1\times S^3$.  Let $S_0'=\pi(S_0)$, where $S_0=0\times S^3$.

Let $\phi\in \Diff_0(S^1\times S^3) $ and $\Sigma'=\phi(S_0')$.  If SS4 is true, then by Lemma \ref{parallel} $\Sigma'$ lifts to $\Sigma\subset \BR\times S^3\subset S^4$ which is isotopic within $\BR\times S^3$ to $S_0$.  Thus $\Sigma'$ would be isotopic to $S_0'$ after lifting to a sufficiently high finite cover. 

Conversely if $\Sigma\subset S^4$ is an oriented 3-sphere, then after isotopy, we can assume that $\Sigma\subset (-1/2,1/2)\times S^3$ and is homologous to $S_0$ there.   Let $\Sigma'=\pi(\Sigma)$ and $S_0'=\pi(S_0)$.  By Theorem 3.13 \cite{BG}, there exists a diffeomorphism $\phi\in \Diff_0(S^1\times S^3) $ such that $\phi (S_0')=\Sigma'$. If $\Sigma$ is not standard, then $\phi$ satisfies the conclusion of Theorem \ref{mazur}. \qed
\vskip 8pt

\noindent\emph{Proof of Mazur's Schoenflies Theorem}:  With notation as in the previous proof, let $\mF$ denote the $S^3$-fibration of $S^1\times S^3$ by pushing forward the standard fibration by $\phi$.  This lifts to a $S^3$-fibration of $\BR\times S^3$ with $\Sigma$ as a leaf.  Thus, as in Mazur's original Schoenflies proof, the closure of each component of $S^4\setminus  \Sigma$ is a smooth $S^3\times [0,\infty)$ whose end limits on the missing point.  \qed

\begin{lemma}\label{product}  If $\Sigma_0, \Sigma_1$ are disjoint, isotopic and non separating 3-spheres in $S^1\times S^3$, then they bound a smooth product.\end{lemma}

\begin{proof}  This is immediate if $\Sigma_0 = x_0\times S^3$.  By \cite{BG} we can assume this is the case.\end{proof}

\begin{definition}  Call a compact oriented contractible 4-manifold $\Delta$ with $\partial \Delta=S^3$ a \emph{Poincare} 4-ball.  Also call it a \emph{Schoenflies} 4-ball if it embeds in $S^3$.  A \emph{Schoenflies} sphere is an oriented embedded 3-sphere in $S^4$.\end{definition}

The following is a restatement of a theorem of Bob Gompf.

\begin{theorem} \label{gompf} (Gompf \cite{Go}) Two Schoenflies balls $\Delta_0, \Delta_1\subset S^4$ are diffeomorphic if and only if they are ambiently isotopic.\end{theorem}

\begin{proof}  It suffices to consider the case $\Delta_0\cap \Delta_1=\emptyset$.  Gompf shows that there exists a diffeomorphism $\phi:(S^4,\Delta_0)\to (S^4,\Delta_1)$.  By precomposing $\phi$ if necessary by a diffeomorphism supported away from $\Delta_0\cup\Delta_1$ we can assume that $\phi$ is isotopic to $\id$.  \end{proof}

\begin{notation} If $f\in\Diff_0(\sonesthree)$, then $\hat f$ will denote a lift to a finite sheeted covering space that will be determined by context and $\tilde f$ will denote a lift to $\rsthree$.   Unless said otherwise the particular lift is immaterial.  \end{notation}

\begin{definition}  Two diffeomorphisms $f,g\in \Diff_0(\sonesthree)$ are \emph{S-equivalent} if there are lifts $\hat f $ and $\hat g$ of $f$ and $g$ to a finite sheeting covering space such that $\hat f(Q)$ is isotopic to $\hat g(Q)$  for some $Q$ of the form $t\times S^3$.

Let $f,g \in \Diff_0(\sonesthree)$.  We say that $f$ \emph{interpolates to} $g$ if there exists a $\tilde h\in \Diff_0(\rsthree)$ such that $\tilde h$ coincides with $\tilde f$ (resp. $\tilde g$) on the $-\infty$ (resp. $+\infty$) end of $\rsthree$. \end{definition}

\begin{proposition}  \label{s equivalent} Let $f,g \in \Diff_0(\sonesthree)$.  The following are equivalent.

i)   $f$ and $g $ are S-equivalent.

ii)  If $P\subset \sonesthree$ is a non separating 3-sphere, then $\tilde f(\tilde P)$ is isotopic to $\tilde g(\tilde P)$.  Here $\tilde P$ is any lift of $P$ to $\rsthree$.

iii)  If $P\subset \sonesthree$ is a non separating 3-sphere, then there exists a finite sheeted covering space $\hat V$ such that $\hat f(\hat P)$ is isotopic to $\hat g(\hat P)$.  Here $\hat P$ is any lift of $P$ to $\hat V$.

iv) $f$ interpolates to $g$ and $g$ interpolates to $f$.

v) There exists a finite sheeted covering space such that $\hat f$ is isotopic to $\hat g$ modulo $\Diff_0(B^4 \fix\partial)$.\end{proposition}

\begin{proof} i) implies v):  Let $Q=t\times S^3$.  Pass to a finite sheeted covering $\hat V$ of $\sonesthree$ so that $\hat f(\hat Q)$ is isotopic to $\hat g(\hat Q)$.  By \cite{BG} $\hat f$ and $\hat g$ are isotopic modulo $\Diff (B^4\fix\partial)$.  \qed

\vskip 8pt

\noindent v) implies iii):  Pass to a finite sheeted covering so that v) holds.  After isotopy $\hat f = \hat g$ modulo $ \Diff (B^4)$.  Since we can assume that the 4-ball is disjoint from $\hat P$ the result follows.\qed

\vskip8pt

\noindent iii) if and only if  ii): The only if  direction is immediate and the other follows from the fact that an ambient isotopy of $\tilde f(\tilde P) $ to $\tilde g(\tilde P)$ in $\rsthree$ can be taken to be compactly supported.\qed

\vskip 8pt

\noindent iii) implies i):  \qed

\vskip 8pt

\noindent ii) implies iv):  We show $f$ interpolates to $g$, Consider $\tilde f, \tilde g:\rsthree\to \rsthree$.  By Cerf \cite{Ce2} and uniqueness of regular neighborhoods, after a compactly supported isotopy of $\tilde f$ to $\tilde f'$, we can assume that $\tilde f'|(N(\tilde Q))=\tilde g|N(\tilde Q)$.  Let $L$ denote the closed complementary region lying to the $-\infty$-side of $Q$ and $R$ the component to the $+\infty$ side.  Define $h:\rsthree\to \rsthree$ by $h|L=\tilde f'|L$ and $h|R=\tilde g$.  A similar argument shows that $g $ interpolates to $f$.  \qed

\vskip 8pt

\noindent iv) implies ii):  Consider an interpolation $\tilde h$ of  $\tilde f$ to $\tilde g$.  Since $\tilde h(t\times S^3)$ is isotopic to $\tilde h(t'\times S^3)$ our assertion follows.  \end{proof}

\begin{corollary}  Interpolation is an equivalence relation.\qed\end{corollary}

Any $f\in \Diff_0(\sonesthree)$ is isotopic to one supported in $S^1\times B$, for any 3-ball $B\subset S^3$ and hence $\Diff_0(\sonesthree)$ is abelian and composition is isotopic to the contenation of two such maps supported on disjoint $S^1\times B^3$'s.  E.g. see \cite{BG}.  From this the next result follows.   

\begin{proposition}  \label{equivalence} The groups of  Schoenflies spheres,  Schoenflies balls and S-equivalence classes of $\Diff_0(\sonesthree)$ are naturally isomorphic.   The isomorphism between  Schoenflies balls and Schoenflies 3-spheres is induced by passing to the boundary, using the outward first orientation convention.  The bijection between S-equivalence classes and spheres is given by $f\to i\circ \tilde f(x_0\times S^3)$.  Furthermore, 

i) (law of composition) boundary connect sum for Schoenflies balls and composition for S-equivalence classes. 
 
ii) (inverse) passing to complementary Schoenflies ball in $S^4$ for Schoenflies balls, reversing orientation of Schoenflies spheres and inverse in $\Diff_0(\sonesthree)$ for S-equivalence classes.\end{proposition}

\begin{notation}  If $M$ is an oriented manifold,  then $\bar M $ denotes $M$ oppositely oriented.  If $\Delta$ is a Schoenflies ball, then $-\Delta$ denotes the complementary Schoenflies ball. \end{notation}

\begin{proposition}  \label{bar delta equivalence} If  $\Delta^4, \Sigma^3$, f correspond under Proposition \ref{equivalence}, then $\bar\Delta^4$ is diffeomorphic to $r_{S^3}(\Delta^4)$, $\bar \Sigma^3$ corresponds to $r_{S^3}(\Sigma^3)$ and  $\bar f$ corresponds to $r_{S^3} f r_{S^3}$.  Here  $ r_{S^3}$  denotes reflection of the $S^3$ factor either in $\BR\times S^3$ or $S^1\times S^3$.  Also $r_{S^3}(\Delta^4)$ is oriented as a subspace of $S^4$ and $r_{S^3}(\Sigma^3)$ has the corresponding boundary orientation.   \qed\end{proposition}

\begin{conjecture} \label{negative equals reverse} If $\Delta$ is a Schoenflies ball, then $-\Delta = \bar \Delta$.  Equivalently, if $\Delta$ is a Schoenflies ball, then $\Delta\times I=B^5$.   \end{conjecture}

The following characterization of the Schoenflies problem is well known.  

\begin{definition}  If $x\in \inte(M)$, then let $M_x$ denote $M\setminus x$.  We say that a diffeomorphism $\phi:M_x\to N_y$ \emph{induces} the diffeomorphism $\phi':M\to N$ if there exists a compact 4-ball $B\subset M$ such that $x\in \inte(B)$ and $\phi|M\setminus B=\phi'|M\setminus B$.\end{definition}

\begin{theorem} \label{tame end} The Schoenflies conjecture is true if and only if for every pair of compact 4-manifolds $M, N$ a diffeomorphism $\phi:M_x\to N_y$ induces a diffeomorphism $\phi':M\to N$.\qed  \end{theorem}

\begin{remark}  Deleting an interior point is not in general sufficient to make homeomorphic but non diffeomorphic manifolds diffeomorphic.  Indeed, Akbulut's original cork \cite{Ak} is an example of a compact contractible 4-manifold $A$ with an involution $f$ on $\partial A$ such that $f$ extends to a homeomorphism but not a diffeomorphism on $A$.  This property continues to hold even if $A$ is punctured for Akbulut shows that there exists a curve $\beta\subset \partial A$ that slices in $A$ such that $f(\beta)$ does not slice.  \end{remark}


\section{Pseudo-Isotopy vs Stable Isotopy}\label{pi vs si}

\begin{definition}  Let $M$ be a compact oriented 4-manifold.  A \emph{pseudo-isotopy}  from $\phi_0$ to $\phi_1$ is a diffeomorphism $f:M\times I \to M\times I$ such that $f|M\times 0=\phi_0\times 0$, $f|\partial M\times I=\phi_0\times \id$  and $f|M\times 1=\phi_1\times 1$. In all cases in this paper, $\phi_0|\partial M=\id$.   Let $M_k$ denote $M\#_k S^2\times S^2$, where $k\in \BN$.  Here all the sums are taken in disjoint 4-balls $B_1, \cdots, B_k$.   We say $\Phi: M_k\times [0,1]\to M_k$ is a \emph{stable isotopy from $\phi_0$ to $\phi_1$} if for $t=0,1$ and all $i$, then $\Phi_t| B_i\#S^2\times S^2=\id$ and  the induced maps $\hat \Phi_t:M\to M$ which are $\id$ on each $B_i$, are isotopic to $\phi_0$ and $\phi_1$ respectively.  When $\phi_0=\id$ we say that $\phi_1$ is the \emph{mapping class induced by the stable isotopy $\Phi$}. Let $\Diff_0^\Sigma(M)$ denote the subgroup of $\Diff_0(M)$ generated by elements stably isotopic to $ \id$. \end{definition}

\begin{lemma} If $ \Phi$ is a stable isotopy from $\phi_0$ to $\phi_1$, then $\Phi^{-1}$ defined by $\Phi^{-1}_t=(\Phi_t)^{-1}$  is a stable isotopy from $ \phi_0^{-1}$ to $\phi_1^{-1}$ and $\bar\Phi$ defined by $\bar\Phi_t=\Phi_{1-t}$ is a stable isotopy from $\phi_1$ to $\phi_0$.\qed\end{lemma}

\begin{lemma} Pseudo-isotopy and stable isotopy are equivalence relations on $\Diff_0(M) $ where $M$ is a compact oriented 4-manifold.\qed\end{lemma}

If $M$ is a closed oriented simply connected 4-manifold and $\phi\in \Diff_0(M)$, then by  Kreck \cite{Kr} p. 645 or \cite{Qu} $\phi$ is pseudo-isotopic to $\id$.  Quinn \cite{Qu} as corrected in \cite{GGHKP} used this to show  that  $\phi$ is stably isotopic to $\id$.   In the topological category Quinn \cite{Qu}, \cite{GGHKP} showed that $\phi$ is isotopic to $\id$, so no stabilization is needed.  Slightly earlier Perron \cite{Pe} proved this for manifolds  without 1-handles and using an argument of Siebenmann proved it for closed simply connected manifolds in general.  Ruberman \cite{Ru} constructed diffeomorphisms of closed simply connected 4-manifolds that are topologically but not smoothly isotopic to $\id$.

Using the seminal work of Cerf \cite{Ce3};
Hatcher, Wagoner and Igusa found obstructions for a pseudo-isotopy of a compact manifold $M$ to be isotopic to $\id$ when $\dim(M)\ge 7$. See \cite{HW}, \cite{Ha2}, \cite{Ig1}.    As detailed in \cite{HW} various elements of the Hatcher-Wagoner theory function when $\dim(M)\ge 4$. The first obstruction, $\Sigma(f)\in$ Wh$_2(\pi_1(M))$, is defined when $\dim(M)\ge 4$ and is the exact obstruction to having a \emph{nested eye} Cerf diagram with only critical points of index 2 and 3 and a gradient like vector field (glvf) that does not involve handle slides and has independent birth and death points.  Chapter 5 \S 6 of \cite{HW} shows that under these circumstances $\Sigma(f)=0$.  The converse follows from  Proposition 3, P. 214 \cite{HW} and its proof.    Call such a nested eye 1-parameter family $q_t$ together with its glvf $v_t$ a \emph{Hatcher - Wagoner family}.

Since it will introduce our basic set up and notation we briefly outline the proof \cite{Qu}, \cite{GGHKP} that if $M$ is a compact simply connected 4-manifold and $\phi\in \Diff_0(M)$ where   $\phi$ is pseudo-isotopic to $\id$, then $\phi$ is stably isotopic to $\id$.  Let $f$ denote a pseudo-isotopy from $\id$ to $\phi$.   Consider a 1-parameter family $q_t:M\times I \to [0,1], t\in [0,1]$ so that $q_0$ is the standard projection to $[0,1]$, $q_1$ is the projection to $[0,1]$ induced from the pseudo-isotopy and $q_t$ is a path from $q_0$ to $q_1$ in the space of smooth maps $M\times I \to [0,1]$ that agrees with $q_0$ on $N(\partial M\times I)$.  Since $\pi_1(M)=1$ it follows that $\Sigma(f)=0$ and hence we can assume that this is a Hatcher - Wagoner family.  Quinn observed that after some modification of the 1-parameter family, analogous to reordering critical points in the proof of the h-cobordism theorem, the essential information of the family can be succinctly captured in what we call the \emph{middle middle level picture}, see p. 353-354 \cite{Qu}. While Quinn did this for the innermost eye, one readily proves the general case which we now state.  Start with $M_k$, where $k$ is the number of eye components, each $\stwostwo$ summand is standardly parametrized and $(x_0,y_0)\in \stwostwo$.  Call the $x_0\times S^2$  2-spheres the \emph{standard red}  spheres $\mR^{\std} :=\{R^{\std}_1, \cdots R^{\std}_k\}$ and the $S^2\times y_0$  2-spheres the \emph{standard green} spheres $\mG^{\std}:=\{G^{\std}_1, \cdots, G^{\std}_k\}$.  Now do a finite sequence of pairwise disjoint finger moves to the red spheres to obtain $\mR=\{R_1, \cdots, R_k\}$ with   new intersections with $\mGs$, one pair for each finger move. Suppose that there is a set of pairwise disjoint Whitney discs such that applying the corresponding Whitney moves to $\mR$  yields a new system of red  spheres that  $\delta_{ij}$ pairwise geometrically intersect the components of $\mGs$. The middle middle level picture is the situation after the finger moves. Thus in $M_k$  we have two sets of spheres $\{R_1, \cdots, R_k\},\{G^{\std}_1, \cdots, G^{\std}_k\}$ and two sets $\mathcal{F} =\{f_p\}, \mathcal{W}= \{w_q\}$ of 
\emph{Whitney discs}  that  cancel the excess $R_i/G^{\std}_j$ intersections. The $\mathcal{F}$ Whitney discs are called \emph{finger discs}. Doing Whitney moves using these discs undoes the finger moves.  Starting with the middle middle level picture we create a 1-parameter family of maps $q_t: M\times [0,1] \times t\to [0,1]$ and glvf's $v_t$,  $t\in [0,1]$, where $q_{1/2}$ has $k$ critical points of index-2 and $k$ critical points of index-3, $G^{\std}_j\subset q_{1/2}^{-1}(1/2)$ is the ascending sphere of the $j$'th index-2 critical point and $R_i  \subset q_{1/2}^{-1}(1/2)$ is the descending sphere of the $i$'th index-3 critical point.  Using $\{w_q\}$ the handle structure on $M\times I\times 1/2$ is modified, as is the corresponding $(q_t, v_t)$, to one whose ascending and descending spheres intersect geometrically $\delta_{ij}$.  The index-2 and index-3 critical points are then cancelled at death critical points after which $q_t $ is nonsingular.  These modifications enable  an extension of $(q_{1/2}, v_{1/2})$ to $(q_t, v_t),  t\in  [1/2,1]$.  Similarly $\{f_p\}$ enables an extension $(q_t, v_t), t\in [0,1/2]$.  This  can be done so that the resulting the resulting $(q_t, v_t)$ is a Hatcher - Wagoner family and the resulting pseudo isotopy is isotopic to $f$.  

\begin{remarks} i) See Chapter 1 \cite{HW} for basic facts about 1-parameter families including descriptions of their low dimensional strata as well as terminology used in this section.

ii)  Using the light bulb theorem \cite{Ga1}, \cite{ST} the red spheres obtained by doing the Whitney moves can be isotoped back to $ \mRs$ and thus the pseudo-isotopy is determined, up to isotopy, by a loop in the embedding space of red spheres.  Our original motivation for proving the light bulb theorem was constructing such a loop when $M=S^1\times S^3$.\end{remarks}

When $\pi_1(M)=1$, \cite{GGHKP} shows that after finitely many stabilizations of the pseudo-isotopy $f$ and modification of the glvf, the ascending and descending spheres from the critical points of the innermost eye component  have no excess intersections.   That component can be eliminated by Cerf's \emph{unicity of death lemma} Chapter 3 \cite {Ce3}, \cite{Ch} or P. 170 \cite{HW}.  By \emph{stabilization of the pseudo-isotopy $f$} we mean first isotope $f$  to be $\id$ on a $B^4\times I$, then replace the $B^4\times I$ with a $(B^4\#S^2\times S^2)\times I$ and extend $f$ to be the identity on the $(B^4\#S^2\times S^2)\times I$.  Finally, by induction on components, a 1-parameter family $q_t$ can be constructed having no critical points and hence a sufficiently stabilized $f$ is isotopic to $\id$ and therefore so is a stabilized $\phi$.    

The following is the main result of this section.  Its proof will be crucial for applications.

\begin{theorem}  \label{stable main} Let $\phi\in \Diff_0(M)$, where $M$ is a compact oriented 4-manifold.  Then $\id$ is stably isotopic to $\phi$ if and only if $\id$ is pseudo-isotopic to $\phi$ by a pseudo-isotopy $f$ with $\Sigma(f)=0$. In particular, when $\id$ is stably isotopic to $\phi$, the stable isotopy $\Phi:M_k\times [0,1]\to M_k$ can be chosen so that $\Phi_0=\id$.
\end{theorem}

We start with the following well known result.

\begin{lemma} \label{chenciner} Let $M$ be a compact 4-manifold and $(q_t, v_t), (q'_t, v'_t)$ be two Hatcher - Wagoner 1-parameter families.  If for some $s\in[0,1]$, not a birth or death point, $q_t=q'_t$ and $v_t=v'_t$ for $t\le s$ and neither $(q_t, v_t)$ nor $(q'_t, v'_t)$ have excess 3/2 intersections  for $t\ge s$, then the associated pseudo-isotopies $f, f'$ are isotopic.  \end{lemma}

\begin{proof}  It suffices to  assume that $s=1/2$,  is after all the births and before the deaths.   The 1-parameter family $r_t :=q_{1-t}$ for $t\le 1/2$ and $r_t:=q'_t$ for $t\ge 1/2$ with corresponding glvf's is of the nested eye type without handle slides or excess 3/2 intersections, hence the singular locus can be eliminated by uniqueness of death.  Thus the induced pseudo-isotopy $g$ is isotopic to $\id$.  Since $f' = g\circ f$, the result follows.\end{proof}

\noindent\emph{Proof of Theorem \ref{stable isotopy}}: First assume that $\id$ is pseudo-isotopic to $\phi$ by the pseudo-isotopy $ f$ with $\Sigma(f)=0$, hence is realized by a Hatcher - Wagoner 1-parameter family $(q_t, v_t)$.   

Using elements of the proof of Theorem 9 \cite{Gay} we can arrange the following.  For $t\in [0,1/8]\cup [7/8,1]$, $q_t$ is non singular.  For $t\in(1/8, 1/4)$ (resp. $(3/4, 7/8)$) $k$ births (resp. deaths) occur.  All the excess 3/2 intersections occur when $t\in (1/4,3/4)$.  For $t\in [1/4, 3/4], q_t=q_{1/4}$.  Also,  for such $t$, $v_t = v_{1/4}$ when restricted to $\qinv_t([0,1/4]\cup [3/4,1])$.  The $k$ critical points of index-2 (resp. index-3) lie in $\qinv_t(1/8, 1/4)$ (resp. $\qinv_t(3/4,7/8)$) all with distinct critical values.  This requires a bit of preliminary work, e.g. as noted in \cite{Gay} the descending spheres of the 2-handles in $M\times 0\times t$ are simple closed curves that may follow non trivial paths in $\Emb(S^1, M)$, however it can be arranged that these paths are constant for $t\in [1/4,3/4]$ and all the movement shoved into $t\in(3/4, 3/4+\epsilon)$.  We can assume that $v_{1/4}$ is the model gradient like vector field on $M\times I\times 1/4$ arising from  the births and has the following features.  Both $\qinv_{1/4}(1/4)$ and $\qinv_{1/4}(3/4)$ are diffeomorphic to $M_k  $ and  respectively denoted $M_k\times 1/4\times 1/4$ and $M_k\times 3/4\times 1/4$ with the first factors equated via the flow of $v_{1/4}$.  Let $\Gs_i\subset M_k\times 1/4\times 1/4$  (resp. $\Rs_j\subset M_k\times 3/4\times 1/4)$ denote the ascending (resp. descending) sphere of the i'th 2-handle (resp. j'th 3-handle).  We can assume that the $\Gs_i$'s (resp $\Rs_j$'s) are the standard green and red spheres in $M_k$.  

\emph{Notational alert}:  In this proof the third coordinate will always refer to the $t$ of $(q_t, v_t)$.  The first two coordinates refer to a point in $M\times I$, however for $t\in [1/4, 3/4]$ the subset $\qinv_t[1/4, 3/4]$ of $M\times I\times t$ is diffeomorphic to and denoted by $M_k\times [1/4, 3/4]\times t$.  Thus a point $(x,s,t)$ may refer to a point in $M\times I\times t$ or a point in $M_k\times [1/4,3/4]\times t$.  The distinction should be clear from context.

Under $v_{1/4}$ the $\Gs_i$'s flow to the standard green spheres in $M_k\times 3/4\times 1/4$ that $\delta_{ij} $ intersect the $\Rs_j$'s, the standard red spheres.  Under $v_{1/4}$ (resp. $-v_{1/4}$) the $\Gs_i$'s (resp. $\Rs_j$'s) flow to discs $\{E_1, \cdots, E_k\}\subset M\times 1\times 1/4$  (resp.$ \{D_1,\cdots, D_k\}\subset M\times 0\times 1/4)$ spanning the ascending (resp. descending) spheres of the 3-handles (resp. 2-handles).  When projected to $M$ the $D_i$'s  intersect $\delta_{ij}$ the $E_j$'s and if $N_i$ denotes a neighborhood of $D_i\cup E_i$, then $v_{1/4}|(M\setminus \cup_{i=1}^k \inte(N_i))\times I=v_0$.  Use $v_{1/4}$ to define a product structure on $\qinv_{1/4}([1/4,3/4])$.  The induced map $M_k\times 1/4\times 1/4\to  M_k\times 3/4\times 1/4$ when projected to the first factor is $\id$ and denoted $\Phi_{1/4}$.  It is constructed from the $\id$ on $M$ (which is induced from $v_0$), replacing each $N_i$ by an $\stwostwo\setminus \inte(B^4)$ and then extending the $\id$, i.e. the map $\Phi_{1/4}$ on $M_k$ is a stabilization of $\id$ on $M$.

    We now show how to \emph{remember} the $\id$ pseudo-isotopy as we pass from $(q_0, v_0)$ to $(q_1, v_1)$ and use it to construct a stable isotopy $\Phi$ from $\id$ to a $\phi'$ which we will show is isotopic to $\phi$ in the paragraph after next.     Since $q_t=q_{1/4}$ for $t\in [1/4,3/4]$,\ $\cup_{t\in [1/4,3/4]}\qinv_t([1/4,3/4])$ is diffeomorphic to $M_k\times [1/4,3/4]\times [1/4,3/4]:=\mC$, the \emph{core} of the 1-parameter family.  The core has two natural parameterizations, the first using the glvf $v_{1/4}$ on each $\qinv_t([1/4,3/4])$ which we call the \emph{standard parametrization} and a second using $v_t$ on $\qinv_t([1/4,3/4])$ which we call the \emph{$\Psi$ parametrization}.  Having defined $\Gs_i\times 1/4\times 1/4$ and $\Rs_j\times 3/4\times 1/4$ the standard parametrization allows us to identify the spheres $\Gs_i\times s\times t$ and $\Rs_j\times s\times t$.  Since $v_{1/4}$ and $v_t$ agree on $\qinv_t([0,1/4])$ we canonically identify $\cup_{t\in [1/4,3/4]} \qinv_t(1/4)$ with $M_k\times 1/4\times [1/4,3/4]$.  We define the parametrization $\Psi:\mC\to \mC$ by $\Psi|M_k\times 1/4\times [1/4,3/4]=\id$, $\Psi|M_k\times [1/4,3/4]\times 1/4=\id$, $\Psi$ fixes each $M_k\times s\times t$ setwise and $\Psi$ takes $v_{1/4}$ flow lines to $v_t$ flow lines.  Since $v_{1/4}$ and $v_t$ agree on $\qinv_t[3/4,1]$, each $\qinv_{s_0}(3/4)$ is canonically identified with $\qinv_{s_1}(3/4)$ for $s_0, s_1\in [1/4,3/4]$ and these identifications agree with the ones given by the standard parametrization.  In particular, $\Rs_j\times 3/4\times 3/4$ is the descending sphere of the $j$'th 3-handle.   We will now assume that $v_t$ has been rescaled so that for $t\in [1/4, 3/4],  v_t$ directly induces $\Psi$, i.e. the flow takes $M_k\times s\times t$ to $M_k\times s\times t$.
    
Since  $\Psiinv(\Rs_j\times 3/4\times 3/4)$ has geometric $\delta_{ij}$ intersection with  $\Gs_i\times 3/4\times 3/4$,  we can isotope $\Psiinv|\Rs_j\times 3/4\times 3/4, j=1, \cdots, k$ to $\id$ via an isotopy staying transverse to $\cup \Gs_i\times 3/4\times 3/4$ by the light bulb theorem \cite{Ga1}, \cite{ST}.  It follows that we can homotope the $v_t$'s $t\in (3/4-\alpha/2, 3/4+\alpha/2)$, keeping $v_t=v_{1/4}$ on $\qinv_t([0,1/4]\cup[3/4,1])$, such that with respect to the new $v_t$'s, $\Psiinv|\Rs_j\times 3/4\times 3/4=\id$.  Also, no new intersections between the ascending and descending spheres of the critical points are created.  Continuing to call our new glvf family $v_t$, the resulting pseudo-isotopy $f$ is unchanged.  Note that when $\pi_1(M)$ has 2-torsion \cite{ST}  applies since $\Psiinv(\Rs_j\times 3/4\times 3/4)$ is isotopic to $\Rs_j\times 3/4\times 3/4$ in $M_k\times 3/4\times 3/4$ and so the Freedman-Quinn obstruction $=0$.   After a second application of the light bulb theorem we can additionally assume that $\Psi|\cup\Gs_i\times 3/4\times 3/4 = \id$.  Here we use the fact that homotopy implies isotopy keeping the dual sphere fixed pointwise provided that  the intersection is preserved under the original map.  Using uniqueness of regular neighborhoods we can further assume that $$(*)\quad\qquad \Psi|\cup_{i=1}^k N(\Gs_i\cup\Rs_i)\times 3/4\times 3/4=\id.$$ 
Note that for  $s\in [1/4,3/4]$, $N(\Gs_i\cup\Rs_i)\times s\times 3/4=\qinv_{3/4}(s)\cap (N_i\times I\times 3/4)$.  
Define $\Phi_s=\Psi|M_k\times s\times 3/4, s\in [1/4,3/4]$, viewed as a map from $M_k$ to $M_k$.  It is a stable isotopy from $\id$ to some $\phi'\in \Diff(M)$.  Call a vector field $v_t$ as above satisfying $(*)$ \emph{stable-inducing}.

We now show that $\phi'$ is isotopic to $\phi$.  Define $\hat\Psi:M\times I\times 3/4\to M\times I\times 3/4$ by $\hat\Psi|\qinv_{3/4}([0,1/4])=\id, \hat \Psi|\qinv_{3/4}([1/4,3/4])=\Psi, \hat\Psi|(\qinv_{3/4}([3/4,1])\cap (N_i\times I))=\id$ all $i$ and $\hat\Psi|\qinv_{3/4}([3/4,1])\cap(M\times I\setminus (\inte(\cup N_i)\times I))$ the extension of $\Psi$ which takes $v_{1/4}$ flow lines to $v_{3/4}$ flow lines.  Recall that $v_{1/4}=v_{3/4}$ in that region.  Next define a 1-parameter family $(q'_t, v'_t), t\in [0,1]$ by $(q'_t, v'_t)=(q_t, v_t)$ for $t\in [0, 3/4]$ and for $t\in [3/4,1]$ first define $\hat\Psi_t:M\times I\times (1-t)\to M\times I\times t$ to be the map which agrees with $\hat\Psi$ on the $M\times I$ factor.  Next define $q'_t=q_{1-t}\circ \hat\Psi_t^{-1}$ and $v'_t=(\hat\Psi_t)_*(v_{1-t})$.  By construction the pseudo-isotopy $f'$ arising from $(q'_t, v'_t)$ is from $\id$ to $\phi'$.  By Lemma \ref{chenciner}, $\phi'$ is isotopic to $\phi$.  
\vskip 10pt
Conversely, suppose that we are given a stable isotopy $\Phi_s, s\in [1/4, 3/4]$ from $\id$ to $\phi$.  We can assume that $\Phi_s=\Phi_{1/4}$ (resp. $\Phi_{3/4}$) for $s$ $\epsilon$-close to $1/4$ (resp. $3/4$). 
 Construct a 1-parameter family $(p_t, \omega_t),\ t\in [0,1]$ as follows.  First, for $t\in [0,3/4]$   define $p_t=q_t$, where $q_t$ is as in the second paragraph of this proof.  
Next define a vector field $\omega$ on $M\times I\times [0,3/4]$ which restricts to a glvf $\omega_t$ on each $M\times I\times t$ as follows.  For $t\in [0,1/4]$ let $\omega_t=v_t$, where $v_t$ is as above.  Also for $t\in [1/4,3/4]$ define $\omega_t=v_{1/4}$ when restricted to $\pinv_t([0,1/4]\cup[3/4,1])$.  For $t\in [3/4-\epsilon, 3/4]$ define $\omega_t$ along $M_k\times [1/4,3/4]\times t$ by $\hat\Phi_*(v_{1/4})$, where $\hat\Phi:M_k\times [1/4,3/4]\times 1/4\to M_k\times [1/4,3/4]\times t$ by  $\hat\Phi(x,s,1/4)=(\Phi_s(x),s,t)$. Use a partition of unity argument to extend our partially defined $\omega$ to a vector field on $M\times I\times [0,3/4]$ which restricts to a glvf on each $M\times I\times t$.  Since $p_{1/4}=p_{3/4} $, $v_t$ agrees with $\omega_t$ as above and there is a diffeomorphism from $M\times [1/4,3/4]\times 1/4$ to $M\times [1/4,3/4]\times 3/4$ which takes $\omega_{1/4}=v_{1/4}$ flow lines to $\omega_{3/4}$ flow lines, the previous paragraph shows how extend the 1-parameter family on $M\times I \times [0,3/4]$ to one on $M\times I\times [0,1]$ yielding a pseudo-isotopy from $\id$ to $\phi$.  In words, we have an identification of $M\times I\times 1/4$ with $M\times I \times 3/4$ which takes the data of one to the other.  Our $(p_t, \omega_t), t\in [0,1/4]$ corresponds to  $k$ standard births.  Reversing this and using the identification informs how to extend our 1-parameter family to have $k$ standard deaths and to see the pseudo-isotopy from $\id$ to $\phi$.  For the latter see the last sentence of the fourth paragraph which shows how to relate the identity pseudo-isotopy with $(q_{1/4}, v_{1/4})$.\qed

\begin{definition} The operation of obtaining a stable isotopy (resp. 1-parameter family) from a 1-parameter family (resp. stable isotopy) as above is called \emph{extracting} (resp. \emph{transplanting}). \end{definition}

\begin{remark}  Our stable isotopy is found by peering inside $(M\times I)\times I$.  In contrast, assuming $\pi_1(M)=1$, \cite{Qu} and \cite{GGHKP} stablize the pseudo-isotopy $f$ itself and then  inductively modify the 1-parameter family to one without critical points to turn the pseudo-isotopy into an isotopy.  \end{remark}

\begin{definition}  If $\phi$ is stably isotopic to $\id$, then the \emph{stabilization number} is the minimal number needed.  If $\phi$ is pseudo-isotopic to $\id$ via a pseudo-isotopy with $\Sigma=0$, then define the \emph{nested eye number} to be the minimal number of nested eyes needed among all Hatcher - Wagoner families yielding a pseudo-isotopy from $\phi$ to $\id$.\end{definition}

The next result follows directly from Theorem \ref{stable isotopy} and its proof.

\begin{corollary} \label{stabilization nested} If $\phi\in \Diff_0(M)$, then the stabilization number is equal to the nested eye number if either is defined.  In particular, both numbers are either defined or undefined.\qed\end{corollary}

\begin{definition}Let $\textrm{PI}(M)$ denote the group of isotopy classes of pseudo-isotopies of $M$ starting at $\id$.  Let $\textrm{PI}^\Sigma(M)$ denote the subgroup of classes $f$ with $\Sigma(f)=0$.  Let $\Diff^\Sigma_0(M)$ denote the subgroup of $\Diff_0(M)$ generated by elements pseudo-isotopic to $\id$ by a pseudo-isotopy $f$ with $\Sigma(f)=0$.  \end{definition}

\begin{remark} \label{wh2} By \cite{HW} p.12, Wh$_2(G)=0$ if $G$ is either free abelian or free,  hence if 
$\pi_1(M)$ is either free or free abelian, then by \cite{HW} $\Sigma(f)=0$ for all pseudo-isotopies of $M$.  \end{remark}

\begin{corollary} \label{stable free} If $\phi\in\Diff_0(M)$ and $\pi_1(M)$ is either free or free abelian, then $\phi$ is pseudo-isotopic to $\id$ if and only if $\phi$ is stably isotopic to $\id$.\end{corollary}

\begin{theorem}  \label{diff stable} If $\phi\in \Diff_0(S^1\times S^3)$ or $\phi\in\Diff_0(M)$ where $M$ is a closed simply connected 4-manifold, then $\phi$ is stably isotopic to $\id$.\end{theorem}

\begin{proof} If $M$ is simply connected and $\phi\in \Diff_0(M)$, then by Kreck \cite{Kr} $\phi$ is pseudo-isotopic to $\id$.   In 1968 Lashoff - Shaneson \cite{LS} and Sato \cite{Sa} proved that if $\phi\in\Diff_0(S^1\times S^q)$, $q=3,4$ then $\phi$ is pseudo-isotopic to $\id$.  Now apply Corollary \ref{stable free}.\end{proof} 

\begin{remark} The punchline of Sato's proof uses Theorem III from the Appendix of \cite{Ke} to prove that a certain homotopy q+1-sphere $\Sigma^{q+1}$ is smoothly standard, but Theorem III is not applicable here for dimensional reasons.  However, he shows that $\Sigma^{q+1}$ is the union of two $B^{q+1}$'s glued along their boundary, hence is standard by Cerf \cite{Ce2}) when $q=3$ or Kervaire - Milnor \cite{KM} and Smale \cite{Sm2} when $q=4$.\end{remark}

\begin{remark} \label{KK} Very recently and after the first version of this paper,  Krannich and Kupers \cite{KK2} proved that there exist diffeomorphisms of certain closed orientable 4-manifolds that are homotopic to but not pseudo-isotopic to $\id$.  It then follows by Theorem \ref{stable main} that such diffeomorphisms are never stably isotopic to $\id$. In addition, elaborating on work of Shaneson \cite{Sh1},\cite{Sh2}, they proved that if $M$ is closed, orientable and $\pi_1(M)$ is free and $\phi\in \Diff_0(M)$, then $\phi$ is pseudo-isotopic to $\id$.  Thus their result together with Corollary \ref{stable free} gives the following extension of Theorem \ref{diff stable}.\end{remark}

\begin{theorem}  \textrm{(Krannich - Kupers)} If $M$ is a smooth closed orientable 4-manifold with free fundamental group and $\phi\in \Diff_0(M)$, then $\phi$ is stably isotopic to $\id$.\end{theorem}

The following result was obtained by David Gay \cite{Gay} for $M=S^4$ and generalized to  $M$  simply connected in \cite{KK1}.  Here we prove the general case with an independent proof.

\begin{theorem} \label{loop to diff} Let $M$ be a compact oriented 4-manifold.  Let $\mR_k$ denote a disjoint union of k ordered 2-spheres.  

i) There is a homomorphism $\gamma_k:\pi_1(\Emb(\mR_k,M_k; \mRs))\to \pi_0(\Diff_0^\Sigma(M))$.

ii) $\pi_0(\Diff_0^\Sigma(M))=\cup_{k=1}^\infty\Image(\gamma_k)$.

iii)  If $ [\alpha]\in \pi_1(\Emb(\mR_k,M_k;\mRs))$ and $[\beta]=[h*\alpha*k]$ where for all $t, i,j, |h_t(\Rs_i)\cap\Gs_j|=|k_t(\Rs_i)\cap\Gs_j|=\delta_{ij}$, then $\gamma_k[\alpha]=\gamma_k[\beta]$.\end{theorem}

\begin{proof}  We abuse notation by identifying $\mR_k$ with $\mRs$.  Let $\alpha_t:\mRs\to M_k$ with $\alpha_0=\alpha_1=\id$.  We can assume that $\alpha_t=\id$ for $t$ near 0 and 1.  Use parametrized isotopy extension to obtain an extension $\alpha_t:M_k\to M_k$.   By \cite{Ga1}, \cite{ST} and uniqueness of regular neighborhoods we can further assume that $\alpha_1|N(\mRs\cup\mGs)=\id$ via an isotopy which does not modify any $\alpha_t|\mRs$.  Thus $\alpha_t$ extends to a stable isotopy of $\id$ and hence induces a $\phi\in\Diff_0(M)$.  The isotopy class of $\phi$ is independent of choices since composing one with the inverse of the other produces a stable isotopy isotopic to one which fixes  $\mRs$ and hence induces a diffeomorphism isotopic to $\id$ by Lemma \ref{chenciner}.  It similarly follows that the isotopy class is independent of the representative of $[\alpha]$.  Since the concatenation of loops gives rise to the composition of stable isotopies the proof of i) is complete.  Given $\phi\in \pi_0(\Diff_0^\Sigma(M))$, construct a pseudo-isotopy $f$ from $\id$ to $\phi$ with $\Sigma(f)=0$ and then extract a stable isotopy $\Phi$ from an associated 1-parameter family  to produce an $\alpha$ with $\gamma_k([\alpha])=\phi$.  Here $\alpha_s=\Phi_s|\mRs$.  Finally for iii), create 1-parameter families by transplanting  stable isotopies arising from $\alpha$ and $\beta$, then use Lemma \ref{chenciner} to conclude that the associated  pseudo-isotopies are isotopic and hence so are the induced elements of $\Diff_0(M)$.\end{proof}

\begin{remarks}\label{wf to fw} i) There is the analogous result with $\mGs$ in place of $\mRs$.

ii) The original motivation to prove the light bulb theorem was to prove i), ii) for $M=\sonesthree$. \end{remarks}

\section {Finger - Whitney Systems}

Let $\Phi:M_k\times I\to M_k$ be a stable isotopy from $\id$ to $\phi\in \Diff_0(M)$ and let $\mRs,\mGs$ be defined as in \S2.  We can  assume that  $\Phi_s(\mRs)$ is transverse to $\mGs$ except for finitely many $s\in I$ corresponding to finger and Whitney moves and all the finger moves occur before the Whitney moves.  For more details see \cite{FQ} pp. 19-20 and \cite{Qu} p. 353.  

\begin{definition}\label{fw} A \emph{finger-Whitney} or \emph{$F|W$ system} $(\mG, \mR, \mF, \mW)$ on $M_k$ consists of two transverse sets of algebraically dual, pairwise disjoint embedded 2-spheres $\mG:=\{G_j\},  \mR:=\{R_i\}$,  (i.e. $\langle R_i,G_j\rangle=\delta_{ij})$  with trivial normal bundles in $M_k$ together with two complete collections of Whitney discs $\mF:=\{f_p\}, \mW:=\{w_q\}$.  This means that performing Whitney moves to $\mR$ using either set of discs produces a set of spheres geometrically dual to $\mG$.  We require that $M$ is diffeomorphic to the manifold obtained by surgering $M_k$ along the $\mG$ family, i.e. replacing neighborhoods of the $G_i$'s by $S^1\times B^3$'s. Call $\mF$ the set of \emph{finger discs} and $\mW$ the \emph{Whitney discs}.   Call a $F|W$ system arising from a stable isotopy of $M$ as in the introduction an \emph{induced $F|W$ system}\end{definition}

\begin{remark} Finger-Whitney systems have their origin in \S 4 \cite{Qu}.  \end{remark}

\begin{lemma} \label{fw unique}  Let $M$ be a compact orientable 4-manifold.  A $F|W$ system on $M_k$ determines a conjugacy  class $\phi(\mG,\mR,\mF,\mW)\subset \pi_0(\Diff_0(M))$. Conversely, given $\phi\in\pi_0(\Diff^\Sigma_0(M))$, there exists an $F|W$ system on some $M_k$ inducing $\phi$.  \end{lemma}

\begin{proof}  Let $\mR'$ (resp. $\mR^{''}$) be obtained from $\mR$ by doing the Whitney moves $\mF$ (resp. $\mW$). Since $M$ is diffeomorphic to $M_k$ surgered along $\mG$, there is a $\zeta_k\in \Diff(M_k)$ such that $\zeta_k( \mRs) =\mR'$ and $\zeta_k(\mGs) =\mG$. To see this note that after surgery the components of $\mGs$ and $\mG$ become circles respectively spanned by 0-framed discs arising from $\mRs$ and $\mR'$ and that a $\zeta\in\Diff_0(M)$ takes one set to the other.   It therefore suffices to consider the case that $(\mG, \mR')=(\mGs, \mRs)$.  

The $F|W$ system gives rise to a loop $\alpha_t\in \Emb(\mR_k, M_k;\mRs)$ obtained by first isotoping to $\mR$ by the finger moves, then to $\mR''$ by the Whitney moves and then back to $\mRs$ by an isotopy fixing $\mG$ setwise.  The last isotopy follows from \cite{Ga1} or \cite{ST}.  Again, since $\mR''$ is isotopic to $\mRs \subset M_k$, it follows that FQ$(\mRs, \mR'')=0$.  By Theorem \ref{loop to diff} $\alpha_t$ induces an element of $\pi_0(\Diff_0(M))$.  Note that   different identifications of $(\mG, \mR')$ with $(\mGs, \mRs)$,  using the method above, will change the resulting class in $\pi_0(\Diff_0(M))$ by conjugation in $\pi_0(\Diff_0(M))$.

Conversely, if $\Phi$ is a stable isotopy of $\id$ to $\phi$, then after isotopy we can assume that $\Phi_s(\mRs)$ is a generic isotopy such that all the finger (resp. Whitney) moves with $\mGs$ occur before (resp. after) $s=1/2$.  lt follows that $(\mGs,\Phi_{1/2}(\mRs),  \mF, \mW)$ is a $F|W$ system on $M_k$, inducing $\phi$ where $\mW$ (resp. $\mF)$ is the set of Whitney (resp. finger) discs.  \end{proof}

\begin{corollary} \label{fw to stable isotopy} A $F|W$ system on $M$ induces a stable isotopy whose associated element of $\pi_0(\Diff_0(M))$ is unique up to conjugacy.\qed\end{corollary}

\begin{corollary} \label{phi to fw} If $\phi\in \pi_0(\Diff_0(S^1\times S^3))$, then there exists a $F|W$ system inducing $\phi$.  A given $F|W$ system on $\sonesthree$ induces a  well defined element of $\pi_0(\Diff_0(\sonesthree))$. \end{corollary}

\begin{proof} Since $\Diff_0(\sonesthree) $ is abelian, the second statement follows from Lemma \ref{fw unique}.  The first follows from the fact that $\pi_0(\Diff_0(\sonesthree))=\pi_0(\Diff_0^\Sigma(\sonesthree))$.\end{proof}

\begin{remark} \label{stable remarks} 

Our original argument for the existence of a $F|W$ system for $\phi\in\Diff_0( S^1\times S^3)$ is as follows.  A given $\phi$ is pseudo-isotopic to $\id$ by \cite{LS}, \cite{Sa}.  It has a Hatcher - Wagoner 1-parameter family by \cite{HW}.  Now apply Quinn's method as outlined in \S \ref{pi vs si}.   \end{remark}

We now specialize to $\phi\in \Diff_0(S^1\times S^3)$, although much of what follows applies more generally. Denote $S^1\times S^3$ by $V$, $S^1\times S^3\#_k S^2\times S^2$ by $V_k$ and $V_\infty:=\BR\times S^3\#_\infty \stwostwo$, where the sums are locally finite and go out both ends.  In practice $V_\infty = \tilde V_k$ for some $k$.  

\begin{notation}\label{winding}  Denote the $S^2\times S^2$ factors of $V_k$ by $S^2\times S^2_i$.  By $S^2\times S^2_i$ we mean a $S^2\times S^2\setminus \inte B^4$, where $B^4$ is disjoint from $\Rs_i\cup \Gs_i $.  Denote by $\Sigma_i$ the 3-sphere that separates off  $S^2\times S^2_i$ from the rest.  When $ (\mG, \mR, \mF, \mW)$ is an $F|W$ system on $V_k$, then we will assume that $\mG=\mGs$ and $\mR$ is isotopic to $\mRs$. Let $Z=V_k\setminus \cup_{i=1}^k \inte(\stwostwo_i)$.  Let $\pi:Z\to S^1$ be the restriction of the projection map from $S^1\times S^3\to S^1$.  Viewing $S^1$ as $[0,k]$ with 0 identified with $k$, we can assume that $\pi(\Sigma_i)=N_\epsilon(i)$. Indices are often chosen mod k.  Here $\tilde \pi$, $\tilde Z$, $\widetilde\stwostwo_j$, $\tilde R_i$, $\tilde G_j$ denote fixed lifts to $\tilde V_k$, where the indices $\in \BZ$.   Let $\rho:\rsthree\to\rsthree$ be the covering translation that shifts everything $k$ units, e.g. $\rho(\tilde R_j)=\tilde R_{j+k}$.\end{notation}

\begin{definition}  Let $w \in \mF\cup \mW$.  We say the \emph{winding} $\omega(w)=r\in \BQ$ if $r=(q-p)/k$ where w lifts to a Whitney disc between $\tilde R_p$ and $\tilde G_q$. Say that $w$ and $w'$ are \emph{winding equivalent} if both $w,w'$ connect $R_i$ to $G_j$ and $\omega(w)=\omega(w')$.  \end{definition}

\begin{remark} Winding depends on how $\pi$ was chosen, however it is well defined when w is between $R_i$ and $G_i$.  Also, when $w,w'$ are Whitney discs between $R_i$ and $G_j$, then the truth of $\omega(w)\neq\omega(w')$ is independent of $\pi$. It follows that the winding equivalence relation is well defined.\end{remark}

\begin{figure}
 \centering
\begin{tabular}{ c c }
 $\includegraphics[width=5.0in]{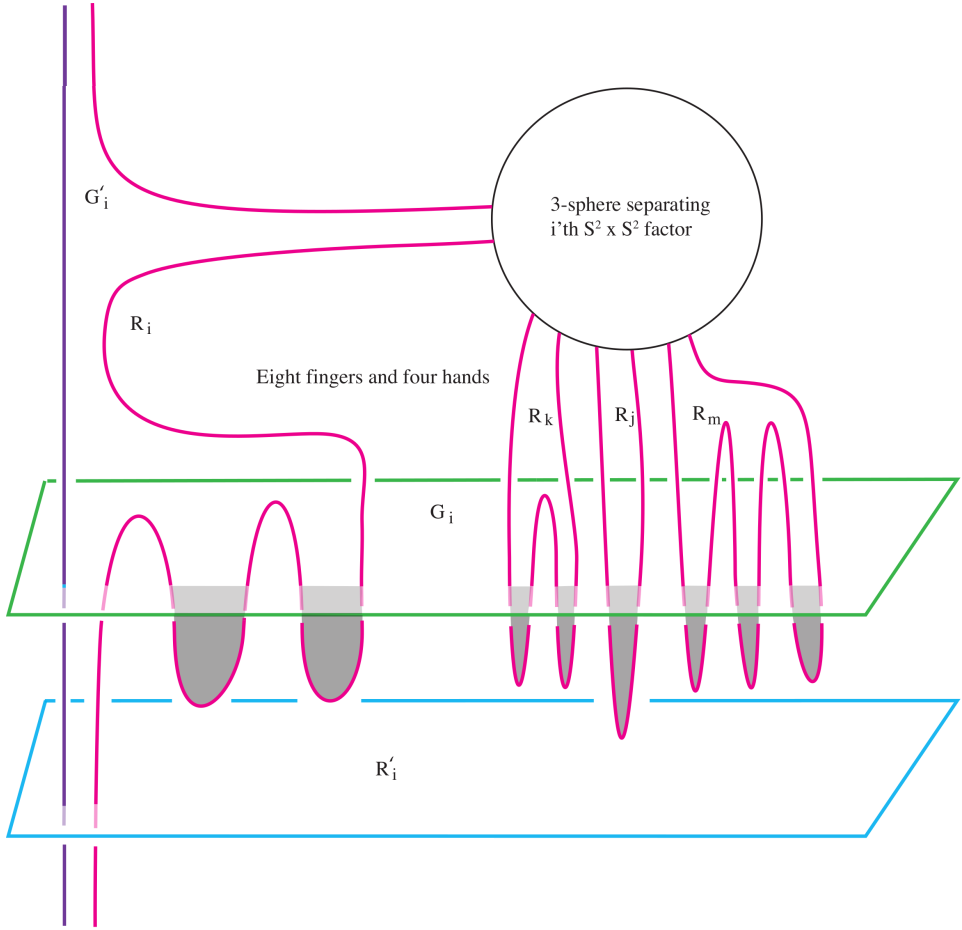}$  
\end{tabular}
 \caption[(a) X; (b) Y]{\label{hand}
 \begin{tabular}[t]{ @{} r @{\ } l @{}}
Arm Hand Finger Form
\end{tabular}}
\end{figure}

\begin{definition} Let $\mF=\{f_1, \cdots, f_p\}$ and $\{[f_{i_1, j_1, \omega_1}], \cdots, [f_{i_m, j_m, \omega_m}]\}$ the winding equivalence classes, where $f_{i,j,\omega}$ denotes the class where a finger goes from $R_i$ to $G_j$ with winding $\omega$.  We say that $\mR$ is in \emph{arm hand finger (AHF) form} with respect to $\mF$ if it is constructed as follows.  For the class $[f_{i,j,\omega}]$ with $|[f_{i,j,\omega}]|$ elements, remove a disc from $\Rs_i$ and replace it by a disc $a_{i,j,\omega}$ called an \emph{arm} which consists of a thin annulus $\subset \stwostwo_i$ that goes essentially straight from $R_i$ to $\Sigma_i$, then a thin annulus $\subset Z$ that goes essentially straight to $\Sigma_j$ and winds $\omega$ about $Z$, then a disc in $\stwostwo_j$ called a \emph{hand} with $|[f_{i,j,\omega}]|$ fingers. Finally, $\mF$ are the discs associated to the fingers. When $ i=j$ and $\omega=0$, then $a_{i,i,0}\subset \stwostwo_i$.  See Figure \ref{hand}.

In similar manner define AHF form of $\mR$ with respect to $\mW$.  And in analogous manner define AHF form of $\mGs$ with respect to either $\mF$ or $\mW$.  This form is the result of an ambient isotopy taking $\mR$ to $\mRs$ and $\mGs$ to $\mG$ where the latter is obtained from $\mGs$ by  attaching arms hands and fingers which poke into $\mRs$, where the discs corresponding to the fingers are the isotoped $\mF$ or $\mW$ discs as applicable. \end{definition}

\setlength{\tabcolsep}{00pt}

\begin{lemma} \label{arm} $\mR$ can be isotoped into AHF form with respect to  $\mF$ as well as $\mW$ via an ambient isotopy that fixes $\mGs$ setwise.  Also $\mGs$ can be isotoped into AHF form with respect to  $\mF$ as well as $\mW$ via an ambient isotopy. \end{lemma}

\begin{proof} We show how to isotope $\mR$ into AHF form with respect to $\mF$.  There exists a set  of pairwise disjoint \emph{finger} arcs from $\mRs$ to $\mGs$ such that if $\mR_1$ is the result of doing finger moves to $\mRs$ along these arcs and $\mF_1$ is the resulting set of finger discs, then $(\mR_1, \mF_1)$ is isotopic to $(\mR, \mF)$ via an isotopy fixing $\mGs$ setwise.  To see this let $\mR'$ be the result of applying the $\mF$ Whitney moves to $\mR$ and noting  $(\mR, \mF)$ arises from $\mR'$ by finger moves along finger arcs to $\mGs$.  Now apply \cite {Ga1} to isotope $\mR'$ to $\mRs$ fixing $\mGs$ setwise and consider the  finger arcs resulting from applying isotopy extension to the finger arcs from $\mR'$.  

Since $\pi_1(V_k\setminus \mGs\cup\mRs)=\BZ$, the result follows because we can isotope the finger arcs so that each goes from an $\Rs_i$ to $\Sigma_i$ to $ \Sigma_j$ to $\Gs_j$ by essentially straight arcs with the correct winding and any two finger arcs corresponding to equivalent fingers are parallel.\end{proof}

\begin{definition}  If $\mR$ is in arm hand finger form, then its \emph{auxiliary discs} are the Whitney discs near the hands as in Figure \ref{hand}.  Note that a hand with $f$ fingers has $f-1$ auxiliary discs.\end{definition}

\begin{definition} \label{fw moves} ($F|W$ Moves) We define the following operations on $F|W$ systems.

i)  \emph{reverse}: $ (\mG,\mR,\mF, \mW)\to (\mG,\mR,\mW, \mF)$.

ii)  \emph{upside down}: $(\mG,\mR, \mF,\mW)\to (\mR,\mG,\mW, \mF)$

iii) \emph{concatenation}:  $(\mG_1, \mR_1, \mF_1, \mW_1)\cup (\mG_2, \mR_2,\mF_2, \mW_2)$.  This an $F|W$ system on $V_{k_1+k_2}$ where the first $F|W$ system is on $V_{k_1}$ and the second on $V_{k_2}$.  Here view $S^1\times S^3=S^1\times B^3\cup_\partial S^1\times B^3$ with each system  of the concatenation supported in its own $S^1\times B^3$.  \end{definition}

\begin{lemma}\label {fw operations} ($F|W$ \textrm{Operations})

i) $\phi(\mG,\mR,\mW, \mF)=(\phi(\mG, \mR, \mF, \mW))^{-1}:=\phi^{-1}(\mG, \mR, \mF, \mW)$

ii) $\phi(\mR,\mG, \mW,\mF)=\phi(\mG, \mR, \mF, \mW)$

iii) $\phi(\mG_1, \mR_1, \mF_1, \mW_1)\cup (\mG_2, \mR_2, \mF_2, \mW_2)) =\phi(\mG_1, \mR_1, \mF_1, \mW_1)\circ\phi(\mG_2, \mR_2, \mF_2, \mW_2)$.\end{lemma}

\begin{proof}  i) If $\phi(\mG, \mR, \mF, \mW)$ arises from the 1-parameter family $(q_t, v_t)$, then $\phi(\mG, \mR, \mW, \mF)$ arises from $(q_{1-t}, v_{1-t})$.  This produces a pseudo-isotopy from $\id$ to $\phiinv$.

ii)  The loop $\beta_s:=\Phi_s|\mGs=\pi_{V_k}\circ\Psi|\mGs\times s\times 3/4$ induces $\phi$.  Let $\beta_{s,t} =\pi_{V_k}\circ \Psi|\mGs\times s\times t$.  Then $\beta_s=\beta_{s,3/4}$ is homotopic to $\beta_{1/4,1-u}*\beta_{v,1/4}*\beta_{3/4,t}$, $u,v,t\in [1/4, 3/4]$, which is homotopic to $\beta_{3/4,t}$ since the first two loops are constant.  Under $\omega_t$ when projected to $V_k$,\ $(\mGs, \alpha_{1/4,t}(\mRs))$ flows to $(\beta_{3/4,t}(\mGs),\mRs)$.  Since $\alpha_{1/4,s}$ approximates $\alpha_{1-s}$, as $s$ increases $\alpha_{1/4, t}(\mRs)$ undergoes inverse Whitney and then inverse finger moves corresponding to the $F|W$ system associated to $\alpha_s$.  This process flows under $\omega_t$ to $\mGs$ undergoing inverse Whitney moves and then inverse finger moves corresponding to the $F|W$ system obtained by flowing the original system seen in $V_k\times 1/4\times 1/2$ to one in $V_k\times 3/4\times 1/2$.

iii) The pseudo-isotopy arising from $\phi(\mG_1, \mR_1, \mF_1, \mW_1)\circ\phi(\mG_2, \mR_2, \mF_2, \mW_2)$ is isotopic to the concatenation of the pseudo-isotopies from $(\mG_1, \mR_1, \mF_1, \mW_1)$ and $(\mG_2, \mR_2, F_2, W_2)$.   By \emph{concatenation} we mean, supported on disjoint vertical $S^1\times B^3\times I$'s $\subset S^1\times S^3\times I$.\end{proof}

\begin{lemma} (factorization lemma) \label{factorization} If $(\mG, \mR, \mF_1, \mF_2), (\mG, \mR, \mF_2, \mF_3)$ are $F|W$ systems on $V_k$, then $\phi(\mG, \mR, \mF_2, \mF_3)\circ\phi(\mG, \mR, \mF_1, \mF_2)=\phi(\mG, \mR, \mF_1, \mF_3)$.\end{lemma}

\begin{proof} The map $\phi(\mG, \mR, \mF_1, \mF_2)$ is induced by the following loop $\alpha_{12}$ in $\Emb(\mR_1,V_k)$ where $ \mR_1$ is the result of applying the $\mF_1$ Whitney moves to $\mR$.  First undo the $\mF_1$ moves to get back to $\mR$, then do the $\mF_2$ Whitney moves to get $ \mR_2$, then isotope back to $\mR_1$, say using the path $\beta$ transverse to $\mG$.  The map  $\phi(\mG, \mR, \mF_2, \mF_3)$ is induced by the loop $\alpha_{23}$ defined by starting at $\mR_1$, next applying the reverse of $\beta$ to get back to $\mR_2$, then applying the reverse  $\mF_2$ moves to get back to $\mR$, then applying the $\mF_3$ Whitney moves, then isotoping to $R_1$ via an isotopy transverse to $\mG$.  But $\alpha_{12}*\alpha_{23}$ is homotopic to the loop $\alpha_{13}$ that induces $\phi(\mG, \mR, \mF_1, \mF_3)$.\end{proof}

The next result follows from Proposition \ref{bar delta equivalence} and the methods of this section.

\begin{lemma}  Let $r$ denote a reflection of $\sonesthree$ across the $S^3$ factor.  If $\mR$, $\mG$ and $\mF$ are invariant under r, then $\phi(\mG, \mR, \mF, \bar\mW)=\bar\phi(\mG, \mR, \mF, \mW)$ where $\bar\mW=r(\mW)$.\qed\end{lemma}
 
\section {Interpolation}

Recall from \S1 that $\phi, \phi'\in \pi_0(\Diff_0(\sonesthree))$ are S-equivalent if they interpolate when lifted to $\rsthree$.  In a similar manner we define S-equivalence for pseudo-isotopies, stable isotopies, 1-parameter families and $F|W$ systems and show that S-equivalence for any two such structures is equivalent to  S-equivalence for the associated element of $\pi_0(\Diff_0(\sonesthree))$.  We show that every structure of a given type is S-equivalent to the trivial one if and only if the Schoenflies conjecture is true.  We give conditions for interpolation of $F|W$ systems and give sufficient conditions for $F|W$ systems to be S-equivalent to the trivial one.  Finally we state a \emph{slice missing slice disc} problem related to $F|W$ interpolation.  Though they can be often stated in more generality, results in this section are given for $S^1\times S^3$ which we again denote by $V$.  As before $V_k$ will denote $S^1\times S^3\#_k S^2\times S^2, k\in \BZ_{\ge 0}$ and $V_\infty$ will denote $\BR\times S^3\#_\infty\stwostwo$.

\begin{definition} \label{sequivalence} The pseudo-isotopy $f$ is \emph{S-equivalent} to $g$ if there exists a pseudo-isotopy $\tilde h :\rsthree\times I\to \rsthree\times I$ such that $\tilde h$ coincides with $\tilde f$ near the $-\infty$ end and with $\tilde g$ near the $+\infty$ end and we say that \emph{$\tilde f$ interpolates to $\tilde g$}.

Let  $p_t, q_t, t\in [0,1]$ be 1-parameter families of $V\times I \to [0,1]$ such that $p_0=q_0$ is the standard projection and both $p_1$ and $q_1$  are non singular.   Then $p_t$ is \emph{S-equivalent} to $q_t$ if there exists a 1-parameter family $\tilde r_t:\rsthree\times I\to [0,1], t\in[0,1]$, such that $\tilde r_t$ coincides with $p_t $ (resp. $q_t$) near  the $-\infty$ (resp. $+\infty$) end, $\tilde r_0$ is the standard projection and $\tilde r_1$ is non singular.  

The stable isotopies $\Phi, \Phi'$ are S-equivalent if their lifts interpolate on $V_\infty$, i.e. if $\Phi$ is defined on $V_k$ and $\Phi'$ on $V_{k'}$, then after identifying the $-\infty$ (resp. $+\infty$) end of $V_\infty$ with the $-\infty$ (resp. $+\infty$) end of $\tilde V_k$ (resp. $\tilde V_{k'}$), there exists a stable isotopy on $V_\infty$ that coincides with $\tilde \Phi$ (resp. $\tilde \Phi'$) near $-\infty$ (resp. $+\infty$).

The $F|W$ systems $(\mG, \mR, \mF, \mW)$, $(\mG', \mR', \mF', \mW')$ are S-equivalent if their lifts interpolate on $V_\infty$, i.e. there exists a  $F|W$ system on $V_\infty$ that agrees with $(\tilde\mG, \tilde\mR, \tilde\mF, \tilde\mW)$ near $-\infty$ and $(\tilde\mG',\tilde \mR', \tilde\mF', \tilde\mW')$ near $+\infty$.  Here we require that the manifolds obtained by surgering $\mG$ and $\mG'$ are diffeomorphic to $\BR\times S^3$.  If these systems are S-equivalent and $\mG'=\mG, \mR'=\mR$ and $\mF'=\mF$, then we say that $\mW$ interpolates to $\mW'$.  \end{definition}

\begin{proposition}\label{pi sequivalence} If for $i=0,1, f_i $ is a pseudo-isotopy from $\id$ to $\phi_i$, then $f_0$ is S-equivalent to $f_1$ if and only if $\phi_0$ is S-equivalent to $\phi_1$.  S-equivalence on pseudo-isotopies on $V\times I$ from $\id$ to elements of $\Diff_0(V)$ is an equivalence relation. \end{proposition}

\begin{proof}  Since S-equivalence is an equivalence relation on elements of $\Diff_0(V)$ and transitivity for S-equivalence of pseudo-isotopies is routine, the first assertion implies the second.  The forward direction of the first assertion is immediate by definition.  

We now show that if $f$ is a pseudo-isotopy from $\id$ to $\phi$ and $\phi$ is S-equivalent to $\id$ then so is $f$.  Let $J=S^1\times y\times I$, some $y\in S^3$.  Keeping $S^1\times S^3\times 0$ fixed isotope $f$ such $f|N(J) =\id$.  
Let $g=f|V\times I\setminus\inte(N(J))$ and consider $\tilde g:\BR\times B^3\times I\to \BR\times B^3\times I$.  In what follows isotopies of $g$ will be constant near the $-\infty$ end as well as on $\BR\times \partial B^3\times I\cup \BR\times B^3\times 0$.  Let $B_t:= t\times B^3\subset \BR\times B^3$ and $\Delta^4_0:=\tilde g(B_0\times I)$.  It follows from the proof of Theorem 9.11 \cite{BG} that since $\phi$ is S-equivalent to $\id, \Delta^4_0\cap (\BR \times B^3\times 1)$ is ambiently properly isotopic to $B_0\times 1$ and hence we can assume by \cite{Ce2} that  $\tilde g|N(B_0\times 1)=\id$.  Let $s>0$ be sufficiently large so that $B_s\times I\cap \Delta^4_0=\emptyset$.  Let $\Delta^5$ be the closure of the region between $B_s\times I$ and $\Delta^4_0$.  Note that $\Delta^5$ is connected and contractible.  Also, $\partial \Delta^5$ is diffeomorphic to $S^4$, again by \cite{Ce2}, since $\partial \Delta^5$ is the union of two 4-balls glued along their boundaries.  It follows \cite{Sm2} that $\Delta^5$ is  diffeomorphic to a 5-ball.  Essentially by uniqueness of regular neighborhoods we can ambiently isotope $\Delta^4$, and hence $\tilde g$, so that  $\Delta^4=B_0\times I:=B^4_0$.  Since $\tilde g|\partial B^4_0=\id$ it follows that  $\tilde g|B^4_0\in \Diff_0(B^4_0\fix \partial) $ and hence is pseudo-isotopic to $\id$ by \cite{Br}.  Therefore, there is an $h:\BR\times B^3\times I\to \BR\times B^3\times I$ so that $h|(-\infty,0]\times B^3\times I=\tilde g|(-\infty,0]\times B^3\times I$, $h|[1,\infty)\times B^3\times I =\id$  and $h|(\BR\times B^3\times 0)\cup \BR\times(\partial B^3\times I)=\id$.  It follows that $f$ is S-equivalent to $\id$.\end{proof}

\begin{proposition}\label{one parameter sequivalence}  The 1-parameter families $q'_t$ and $q_t$ are S-equivalent if and only if their corresponding pseudo-isotopies $f', f$ are S-equivalent. S-equivalence of 1-parameter families is an equivalence relation.\end{proposition}

\begin{proof} If $q_t$ and $q'_t$ are S-equivalent, then by restriction so are $f$ and $f'$.  The converse follows from the contractibility of the space of smooth functions  to $[0,1]$.\end{proof}

\begin{proposition} \label{stable sequivalence} The stable isotopies $\Phi, \Phi'$ are S-equivalent if and only their corresponding transplantations  $f, f'$ are S-equivalent. \end{proposition}

Before the proof we give the following.

\begin{corollary} \label{fw sequivalence} The $F|W$ systems $(\mG, \mR, \mF, \mW)$ and  $(\mG', \mR', \mF', \mW')$ are S-equivalent if and only if $\phi:=\phi(\mG, \mR, \mF, \mW)$ and  $\phi':=\phi(\mG', \mR', \mF', \mW')$ are S-equivalent. \end{corollary}

\begin{proof}  If the $F|W$ systems are $S$-equivalent, then the proof of Lemma \ref{fw unique} shows that they induce isotopic diffeomorphisms.  Conversely, if $\phi$ and $\phi'$ are S-equivalent, then the $F|W$ systems induce S-equivalent stable isotopies by the Proposition, hence the $F|W$ systems are S-equivalent.\end{proof}

\noindent\emph{Proof of Proposition \ref{stable sequivalence}}
The forward direction is immediate.  To minimize notation we will consider the case that $\Phi'$ is the trivial stable isotopy and hence $f'=\id$, since general case is similar.

Suppose that $(q_t, v_t)$ is a transplanted 1-parameter family arising from $\Phi$.  It  has a nested eye Cerf diagram only involving critical points of index-2 and 3 without  2/2 or 3/3 intersections, i.e. the corresponding 1-parameter family of handle structures on $V\times I$ has no handle slides.  Let $(\tilde q_t, \tilde v_t)$ be the lift to $\rsthree\times I$.  Consider the \emph{$\BR$-Cerf diagram} on $\BR\times I\times I$.  Here a point $(u,s,t)$ is labeled $i$, if $\tilde q_t$  has a critical point  $\in u\times S^3\times s$ of index-$i$.  The $\BR$-Cerf diagram of $(\tilde q_t, \tilde v_t)$ is a locally finite union of single eye components, only involving critical points of index-2 and 3, without  2/2 or 3/3 gradient intersections.

There is a generic 1-parameter family $\tilde p_t:\BR\times S^3\times I\to [0,1]$, $t\in [0,1]$ with glvf $\tilde w_t$, such that  $\tilde p_t|(\infty, -1]\times S^3\times I=\tilde q_t$  and $\tilde p_t|[1,\infty)\times S^3\times I$ is the standard projection to $[0,1]$ with $\tilde w_t$ the vertical vector field.  Now and in the future $\tilde p_t(\BR\times S^3\times i)=i$ for $i$ close to 0 or 1; $\tilde p_t$ is the standard projection to $[0,1]$ for $t$ close to 0 and $\tilde p_t$ is non singular for $t$ close to 1.  That such a $(\tilde p_t, \tilde w_t)$ exists follows from the contractibility of smooth maps to $[0,1]$ and the hypothesis that $f$ is S-equivalent to $\id$.  Also all modifications of $\tilde p_t $ will be compactly supported, hence will always respectively agree with $\tilde q_t$ and $\id$ near the negative and positive ends of $\BR\times S^3\times I$.   To complete the proof it suffices to show that $(\tilde p_t, \tilde w_t)$ can be modified  so that the components of the $\BR$-Cerf diagram  are eyes with edges labeled 2 and 3 with independent births and deaths and without gradient intersections of type 2/2 or 3/3.   The methods of \S 2 and \S 3 then show that after a compactly supported modification of $\tilde p_t$ that the resulting proper stable isotopy  is an interpolation of $\Phi$ to the trivial stable isotopy.  

We  show that the proof of Proposition 3, Chapter VI \cite{HW} extends to our  setting.   Since  $\tilde p_t$  is generic and for $t\in [0,1]$ agrees with either $\tilde q_t$ or the standard projection off of a uniform compact set, all the components of the $\BR$-Cerf diagram are compact, only finitely many are involved with 2/2 or 3/3 gradient intersections, or have critical points not of index-2 or 3 or have more than one birth or death, or have non independent births or deaths.    Let $C$ denote the union of these components.  Since the proof establishing the first paragraph of p. 214 \cite{HW} only requires that $C$ be compact we can assume that the births and deaths of  $\tilde p_t$ are independent and the non birth/death critical points are of index-2 or 3. In a similar manner we can assume that conclusions of Step 1 and Steps 2 of the proof of Proposition 3 hold for $\tilde p_t$.  Note that since $\dim(V)=4$, the 3/3 gradient intersections are traded for 2/2 intersections.  We continue to denote by $C$ the components of the $\BR$-Cerf diagram involved with 2/2 or 3/3 gradient intersections.

We now adapt Step 3 of the proof including its notation and terminology.  Since the 2/2 intersections are supported in the compact set $C$, the argument of \S 6, Chapter V of \cite{HW} shows that the $2/2$ gradient intersections correspond to a word $x\in K_2[Z[\pi_1(\rsthree\times I)]] $ which is a subgroup of the Steinberg group St$(\BZ[1])= $St$(\BZ)$.  Since $W_2$ of the trivial group equals 1 \cite{HW} there is a word $w\in W(\pm 1)\subset St(\BZ)$ such that $x\cdot w = 1$ in  St$(r_1,\BZ)$ for $r_1$ sufficiently large.  (See p. 9-10 \cite{HW} for definitions of these algebraic objects. Note that multiplicative group notation is used.) If $r_0$ is the number of edges of $C$ labeled 2 we can assume that $r_0\le r_1$.  By modifying $\tilde p_t$ to add $r_1-r_0 $  trivial eye components we can assume that $r_1=r_0$.  By Lemma 2.7 Chapter IV \cite{HW}, we can modify $\tilde p_t$ to add a new  graphic representing the word w and then using (0.1) V \cite{HW}  to join the graphics to create a new one, still denoted $C$, whose word corresponding to the 2/2 gradient intersections is now $x\cdot w$.  By the proof of Theorem 1.1, Chapter II \cite{HW} $p_t$ can be homotoped to eliminate these 2/2 gradient intersections via a homotopy that is supported on a compact region of $\BR\times S^3\times I$.  No 3/3 gradient intersections are created in the above process.

Note that both after the births and before the deaths there is geometric $\delta_{ij} $ intersection between the ascending spheres of the 2-handles and the descending spheres of the 3-handles.   The 3/2 gradient intersections preserve the $\delta_{ij}$ pairing algebraically, so if there are no 2/2 or 3/3 gradient intersections,  edges that are born together  must die together.  These facts together with the fact that  all the births and deaths are independent enable us to modify $(\tilde p_t, \tilde w_t)$, with compact support, so that the usual Cerf diagram is of the nested eye type, though infinitely many of the components of the $\BR$-Cerf diagram may map to a single component of the usual diagram.\qed
\vskip 8pt

 Putting this all together we have the following 

\begin{theorem}  \label{main sequivalent} The following are equivalent

i) The Schoenflies conjecture is true,

ii) every $\phi\in \Diff_0(S^1\times S^3)$ is S-equivalent to $\id$, 

iii) every pseudo-isotopy on $S^1\times S^3\times I$ from $\id$ is S-equivalent to $\id$,

iv) every $F|W$ system on $S^1\times S^3$ is S-equivalent to the trivial system,

v) every stable isotopy from the $\id $ on $S^1\times S^3$ is S-equivalent to $\id$.\qed\end{theorem}

\begin{definition} Define an abelian group structure on $F|W$ systems on $V$ whose elements are S-equivalence classes and whose addition is induced by the bijection with S-equivalence classes of  $\Diff_0(V)$.  \end{definition}

We now describe interpolation operations on $F|W$ systems.

\begin{lemma}  (Disjoint Replacement)\label{disjoint replacement}  Let $(\mG, \mR, \mF, \mW)$ and $(\mG, \mR, \mF, \mW')$ be $F|W$ systems on $V$.  Suppose that $\partial\mW=\partial \mW'$  and $\inte(\mW)\cap\inte(\mW')=\emptyset$.  Then the two systems are S-equivalent. \end{lemma}

\begin{proof}  Replace all the elements of $\tilde \mW$ near the $+\infty$ end of $V_\infty$ by corresponding elements of $\tilde\mW'$.\end{proof} 

\begin{remarks}  i) Often a weaker version of this condition suffices, e.g. $\partial\mW=\partial \mW'$ and $\inte(\mW)\cap\inte(\mW')\neq\emptyset$ but after swapping  $\tilde \mW$ discs by $\tilde \mW'$ discs near the $+\infty$ end of $V_\infty$ the resulting system of Whitney discs in $V_\infty$ is embedded.

ii)  Examples arise when $\mW'$ is obtained from $\mW$ by introducing local knotting and linking to its components.  
\end{remarks}

   Geometrically dual spheres disjoint from the finger and Whitney discs provide another useful tool.
  
 \begin{lemma} (Dual Sphere Lemma) \label{dual spheres} Let $(\mG,\mR,\mF,\mW)$ be an $F|W$ system supported on $V_k$.  Suppose there exist pairwise disjoint embedded spheres $H_{i_1}, \cdots, H_{i_m}$ such that $i_1, \cdots, i_m$ are distinct elements of $\{1, 2, \cdots, k\}$ and 

i) $H_{i_j}\cap (\mR\cup \mG)=H_{i_j}\cap (R_{i_j}\cup G_{i_j})=1$ and

ii) $H_{i_j}\cap (\mF\cup\mW)=\emptyset$ all $i_j$, 

\noindent then $\phi(\mG,\mR,\mF,\mW)$ is S-equivalent to an explicit $F|W$ system supported on some $V_{k-m}$. In particular if $m=k$, then this is the trivial system and $\phi$ is S-equivalent to $\id$.  \end{lemma}
 
 \begin{proof}  By reordering the $S^2\times S^2_i$ factors we can assume that $i_j=j$.  Let $\widetilde \stwostwo_i, i\in \BZ$ denote the preimages of the $S^2\times S^2_i$ factors in $\tilde V_k$ with the covering translation $\rho$ shifting everything $k$ units.  We can assume that $\mR$ is in arm hand finger form.  

If possible, reorder so that $H_1\cap G_1\neq\emptyset$, otherwise proceed as in four paragraphs below.   If $|G_1\cap \mR|=1$, then replace $H_1$ by a translate of $G_1$ and begin again.  Observe that $H_1$ has trivial normal bundle since $[H_1]=[R_1]\in H_2(V_k)$ and is embedded.  Modify $\tilde\mF|\tilde\mW$ as follows.  First, isotope $\tilde \mR$ to $\tilde\mR_1$ by doing Whitney moves using exactly all $f\in \tilde \mF$ such that $\partial f\cap \tilde G_j\neq\emptyset$ where  $j=1$ modulo $k$ and $j\ge 1$.  Denote by $\hat\mF\subset\tilde\mF$ these $f$'s and $\tilde\mF_1:=\tilde\mF\setminus \hat\mF$.  Define $\hat\mW=\{w\in \tilde \mW|\partial w\cap G_j\neq\emptyset$ where  $j=1$ modulo $k$ and $j\ge 1\}$ and $\tilde \mW':=\tilde\mW\setminus \hat W$.  

The isotopy from $\tilde\mR$ to $\tilde\mR_1$ is supported very close to $\hat\mF$, in particular so that $\tilde\mR_1$ remains disjoint from the $\tilde H_i$'s.  If $w\in \tilde \mW' $ and $w\cap \hat\mF\neq \emptyset$, then under isotopy extension $w$ gets moved to $w_1$.  Let $\mW_1'=\tilde\mW'$ with these $w$'s replaced by their $w_1$'s.  

This creates two types of problems for $(\tilde \mG, \tilde \mR_1, \tilde \mF_1, \tilde \mW_1') $ being an $F|W$ system.  If $f\in \hat \mF$ and $w\in \tilde\mW'$ with $|\inte f\cap w|=p$, then $w_1$ will have $2p$ new intersections with some $\tilde G_j$ while if $|w\cap\partial f|=p$, then $w_1$ will have $p$ additional new intersections with some $\tilde G_j$.  In both cases, $j\ge 1$ and $j=1$ modulo $k$.     Now modify the $w_1$'s by tubing off with copies of the component of $\tilde H_1$ that intersects $\tilde G_j$.  Let $\tilde\mW_1$ denote the modified $\tilde\mW_1'$.  

Do all the modifications $\rho$-equivariantly.  This means that if $f\in \hat \mF$,  then the Whitney move associated to $\rho(f)$ is $\rho$ of the Whitney move associated to $f$.  Similarly, if $w\in \tilde\mW'$ and $w\cap \hat\mF\neq\emptyset$, then $(\rho(w))_1\cap N(\hat\mF)=\rho(w_1)\cap N(\hat\mF)$.  Here the isotopy of $\tilde \mR$ to $\tilde \mR_1$ is supported in $N(\hat\mF)$.   Also the tubings are done $\rho$-equivariantly as well.   It follows that if $(\tilde \mG, \tilde\mR_1, \tilde \mF_1, \tilde \mW_1)$ is the new system, then the restriction to the $+\infty$ end projects to a $F|W$ system $(\mG, \mR_1, \mF_1, \mW_1)$ on $V_k$ such that  $|\mR_1\cap G_1|=1$ and  $H_1,  \cdots, H_m$ satisfy the same properties as before with this $F|W$ system.  Now replace $H_1$ by a translate $H'_1$ of $G_1$.  

By induction on $|\cup H_i\cap \mG|$ we now assume that for all $i, H_i\cap \mG=\emptyset$. Abuse notation by denoting the new $F|W$ system arising from the construction of the previous paragraphs by $(\mG,\mR,\mF,\mW)$.  Now modify $(\tilde \mG, \tilde\mR, \tilde \mF, \tilde \mW)$ as follows. For $i\ge 1$ and $i = j$ modulo $k$, where $1\le j\le m$, obtain $\tilde R_i' $ from $\tilde R_i$ by doing Whitney moves using all $f\in \tilde\mF$ such that $\partial f\cap\tilde R_i\neq\emptyset$.  Again, the Whitney moves are done very close to these $f$'s and $\rho$-equivariantly. Let $\tilde \mR_1$ denote the modified $\tilde \mR$.  Let $\hat\mF\subset \tilde \mF$ denote the union of discs used in these moves and $\tilde \mF_1=\tilde\mF\setminus \hat\mF$.  Let  $\hat W=\{w\in \tilde\mW|w\cap \tilde R_i\neq\emptyset$ for $i\ge 1$ and $i = j $ modulo $k$, where $1\le j\le m \}$.   Let $\tilde \mW'=\tilde\mW\setminus \hat W$.  

The problem with $(\tilde \mG, \tilde\mR_1, \tilde \mF_1, \tilde \mW')$ is that if $w\in \tilde \mW'$, then $w\cap \tilde \mR_1$ may have intersections.  Indeed if $f\in \hat F$, then each point of $\inte(f)\cap w$ will give rise to two such intersections and each point of $\partial f\cap w$ will give rise to one.  Modify $\tilde \mW'$ to $\tilde \mW_1$ by $\rho$-equivariantly tubing off each $\tilde\mW'\cap \tilde \mR_1$ intersection with a copy of some $\tilde H_i$.  Project the $+\infty$ end of $(\tilde\mG, \tilde \mR_1, \tilde\mF_1, \tilde\mW_1)$ to $V_k$ to obtain an $F|W$ system satisfying the conclusions of the lemma.

If $m=k$, then the resulting $F|W$ system will have $\mF=\emptyset $ and hence the corresponding $\phi$ is S-equivalent to $\id$.  Also, the interpolation is done explicitly and algorithmically.  Indeed, the final $F|W$ system on $V_k$ can be explicitly constructed directly in $V_k$.  \end{proof}

\begin{definition} A component \emph{component} $\mC=(\mG_C,\mR_C,\mF_C,\mW_C)$ of $(\mG,\mR, \mF, \mW)$ consists of those elements lying in a connected component of $\mG\cup\mR\cup \mF\cup \mW$.  It is supported in $V_{k-m}$ where $m=|\mG\setminus \mG_C|$ and is obtained from $V_k$ by surgering the components of $\mG\setminus \mG_C$.   \emph{Deletion} of $\mC$ is the $F|W$ system obtained by removing $\mC$ from  $(\mG,\mR, \mF, \mW)$ and surgering $\mG_C$.  Note that it is supported on $V_m$.  We say that $\mC$ is \emph{S-trivial} if the induced map $\phi(\mC)$ is S-equivalent to $\id$.  The $F|W$ system is \emph{connected} if it has one component.    \end{definition}

\begin{remarks} i) If each component of $\tilde\mC$ is compact, then $\mC$ is S-trivial.  

ii) If $h\in \Diff_0(B^4 \fix \partial)$, then $h $ is pseudo-isotopic to the identity so by \cite{HW} and  \cite{Qu} it arises from an $F|W$ system on $B^4\#_k\stwostwo$ which can be viewed as an S-trivial $F|W$ system on $\sonesthree$.  

iii) Concatenation can be viewed as the operation of addition of unions of components.

\end{remarks}

\begin{lemma}  If the  $F|W$ system $(\mG,\mR,\mF,\mW) $ is the disjoint union of components $\cup_{i=1}^m(\mG_i,\mR_i,\mF_i, \mW_i)$, then $\phi(\mG,\mR,\mF,\mW)=\sum_{i=1}^m \phi(\mG_i,\mR_i,\mF_i, \mW_i)$.\end{lemma}

\begin{proof} By disentangling in $V_\infty$ we will show that $(\mG,\mR,\mF,\mW)$ is S-equivalent to the concatenation of its components.  We give the proof for $m=2$, the general case follows by induction.  Suppose that $|\mG_i|=k_i$.  View $\mR$ in arm hand finger form and let $L$ denote an unknotted $S^1\times S^2$ which separates $\mG_1\cup \mR_1\cup\mF_1\cup \mRs_1$ from $\mG_2\cup \mR_2\cup\mF_2\cup\mRs_2$.  This means that the closed complementary regions are copies of $S^1\times B^3\#_{k_i}\stwostwo$, $i=1,2$ respectively denoted $V'_{k_1}, V'_{k_2}$.  Note that $\mW_1$ is disjoint from $\mG_2\cup \mR_2\cup \mF_2$ which deformation expands to a space $X_2\subset V'_{k_2}$ which is $N(\mG_2\cup\mRs_2)$ union finitely many 1-handles.  Using isotopy extension we can assume that $\mW_1\cap X_2=\emptyset$.  Let $\mW_1' $ denote the result of isotoping $\mW_1$ into $V'_{k_1}$ via an isotopy fixing $\mG_1\cup \mR_1\cup \mRs_1$ pointwise.   The track of this isotopy may cross $X_2$.   In a similar manner construct $\mW_2'$.  

We now show that $(\mG,\mR,\mF,\mW)$ is S-equivalent to the concatenation of $(\mG_1,\mR_1,\mF_1, \mW_1') $ and $(\mG_2,\mR_2,\mF_2, \mW_2')$.  Consider $\tilde V_k$ where $k=k_1+k_2$.  Let  $\widetilde \stwostwo_i$, $i\in \BZ$ denote the preimages of the $S^2\times S^2$ factors in $V_k$ with the covering translation $\rho$ shifting everything $k$ units.  Order the factors so that when $i=j$ modulo $k$, where $1\le j\le k_1$ (resp. $k_1+1\le j\le k$), then these factors lift from $V'_{k_1}$ (resp. $V'_{k_2}$).  Now let $\tilde \mW^1_i$ be a Whitney system for $\tilde \mG_i\cup\tilde R_i $ that coincides with $\tilde\mW_i$ near $-\infty$ and $\tilde\mW_i'$ near $+\infty$ and $\tilde\mW^1_i\cap \tilde X_j=\emptyset$ where $i\neq j$.   Such systems exist because we can $\rho$-equivariently isotope  $\tilde \mW_i$ to $\tilde \mW_i'$ only near the $+\infty$ end, so  don't have unwanted intersections between moved and unmoved components of $\tilde \mW_i^1$.

To complete the proof we modify $\tilde W^1_1$ to eliminate the finite set $\tilde W^1_1\cap \tilde W^1_2$.  First isotope $\tilde \mW^1_1$ off of  $\tilde \mW^1_2$ to obtain $\tilde\mW^2_1$ at the cost of creating twice as many $\tilde\mW^2_1\cap \tilde\mR_2$ intersections.  Now $\tilde\mR_2$ has pairwise disjoint geometrically dual spheres that are disjoint from $\tilde \mG_2\cup \tilde \mW^1_2\cup \tilde\mW^2_1\cup\tilde\mG_1\cup\tilde\mR_1 $ by doing Whitney moves to $\tilde \mG_2$ using $\tilde W^1_2$ and isotoping slightly.  Finally create $\tilde\mW^3_1$ by tubing $\tilde\mW^2_1$ to copies of the dual spheres, one for each point of $\tilde\mW^2_1\cap \tilde\mR_2$.  \end{proof}

\begin{corollary}  (Deletion Lemma) \label{deletion} If $\mF'|\mW'$ is the subsystem of $\mF|\mW$ with all S-trivial components deleted, then $\phi(\mF|\mW)$ is S-equivalent to $\phi(\mF'|\mW')$.\qed\end{corollary}

\begin{definition} Let $\mG_1\subset \mG, \mR_1\subset \mR$ be such that $\mG_1\cup\mR_1$ is a union of components of $\mG\cup \mR$.  We say that $(\mG,\mR',\mF', \mW')$ is obtained by \emph{contracting $(\mG,\mR,\mF, \mW)$ along $\mG_1, \mR_1$}, if $\mR'$ is obtained by replacing $\mR_1$ by $\mRs_1$; $\mF' =\{f\in \mF|\partial f\cap \mR_1=\emptyset\}$; $\mW'' =\{w\in \mW|\partial w\cap \mR_1=\emptyset\}$ and $\mW' = \mW''$ modified as follows.  If $w\in \mW''$   and $\inte (w)\cap \mRs_1\neq\emptyset$, then modify by tubing each intersection of $w \cap \Rs_j$ with a copy of $G_j$. \end{definition}

\begin{remark} Since we can assume that $\mR$ is in AHF form with respect to $\mF$, $\mRs_1\cap (\mR\setminus \mR_1)=\emptyset$.\end{remark}

\begin{lemma} (Contraction Lemma) \label{contraction} If $(\mG,\mR',\mF', \mW')$ is obtained from $(\mG,\mR,\mF, \mW)$ by contracting along those components of $\mG\cup\mR$ whose induced maps into $\pi_1(V_k)$ are trivial, then $\phi(\mG,\mR,\mF, \mW)$ is S-equivalent to $\phi(\mG,\mR',\mF', \mW')$.  In particular, if all components of $\mG\cup\mR$ are $\pi_1$-inessential, then $\phi(\mG,\mR,\mF, \mW) $ is S-equivalent to $\id$. \end{lemma}

\begin{proof} Let $U$ be a neighborhood of the $+\infty$ end of $\tilde V_k$ that is disjoint from some neighborhood of the $-\infty$ end.  Replace all the lifts of the red spheres contained in all the $\pi_1$-inessential components of $\mG\cup\mR$ which intersect $U$ by the corresponding standard red spheres.  Obtain the new $\tilde\mF$ by deleting those elements whose boundaries intersected the red spheres that were replaced.  Obtain the new $\tilde\mW$ by also deleting those elements whose boundaries intersected the replaced red spheres.  Modify those remaining discs $\in \tilde\mW$ that intersect the new red spheres by tubing off intersections using copies of their dual green spheres.  \end{proof}

\begin{definition}  We say that $(\mG,\mR,\mF,\mW)$ has \emph{fingers monotonically pointing up} (resp. \emph{down}) if the $\stwostwo$ factors can be ordered in $V_k$ so that in $\tilde V_k, \tilde R_i\cap \tilde G_j\neq\emptyset$ implies $j\ge i$ (resp. $j\le i$).\end{definition}
\begin{proposition}  If $(\mG,\mR,\mF,\mW)$ has either fingers monotonically pointing up or down, then $\phi(\mG,\mR,\mF,\mW)$ is S-equivalent to $\id$.\end{proposition}

\begin{proof} Start with $\mR$ in arm hand finger form and the $\stwostwo$ factors of $\tilde V_k $ ordered as usual with the covering translation $\rho$ shifting  $\widetilde \stwostwo_i$ to $\widetilde\stwostwo_{i+k}$. It suffices to consider the case that all the fingers point down.  Replace $\tilde R_i, i\ge 0$ by $\tilde R^{\std}_i$ and eliminate those elements of $\tilde F\cup \tilde W$ which intersect $\tilde R_i, i\ge 0$.  The cost is that some of the remaining elements of $\tilde W$ may intersect $\tilde R^{\std}_i, i\ge 0$, however there are only finitely many such intersections.  Tube off each of  intersections with $\tilde R^{\std}_i$ with a copy of $\tilde G_i$.  Since $(\tilde \mG,\tilde \mR,\tilde\mF,\tilde\mW)$ interpolates to $(\tilde\mG, \tilde\mR^{\std}, \emptyset,\emptyset)$ the result follows.\end{proof}

\begin{remark}  This argument more generally shows the following.  Given $(\mG,\mR,\mF,\mW)$ as in Notation {3.5}, define a graph $G\subset S^1\times D^2$ as follows, where $S^1=[0,k]/\sim$.  The vertices are $\{i\times 0|i\in k\}$.   To each finger $f$ from $R_i$ to  $G_j$, construct an embedded directed edge from $i\times 0$ to $j\times 0$ which winds $\omega(f)$ about the $S^1$.   If no component of $G$ contains both a cycle representing a positive element of $\pi_1$ and a cycle representing a negative element, then $\phi(\mG,\mR,\mF,\mW)$ is S-equivalent to $\id$.\end{remark}

The next lemma shows that if $\mF$ and $\mW$ have common sets of dual spheres then they are interpolable.

\begin{definition}  A set of pairwise disjoint embedded 2-spheres $\mN_f=\{N_1, \cdots, N_n\}$  with trivial normal bundles is said to be \emph{dual} to $\mF=\{F_1, \cdots, F_n\}$ if $\mN_f\cap \mG\cup \mR=\emptyset$ and $|N_i\cap F_j|=\delta_{ij}$.  In a similar manner we define the notion of dual spheres to $\mW$.\end{definition}

\begin{lemma} \label{whitney duals} (Whitney duals exist) Given the $F|W $ system $(\mG,\mR,\mF,\mW)$, there are dual spheres $\mN_f$ and $\mN_w$ to $\mF$ and $\mW$.  \end{lemma}

\begin{proof} We argue as in \cite{FQ}.  Given $F_i$, let $D_i$ be a small  2-sphere near $F_i\cap R_i$ that intersects $ F_i$ once and $R_i$ twice.  For each $i$, tube $D_i$ with two copies of the dual sphere $R_i'$ to $R_i$, where $R_i'$ is as in Figure \ref{hand}, to eliminate these points and thereby construct $N_i$ and hence $\mN_f$.  In a similar manner construct $\mN_w$.\end{proof}

\begin{remarks} We can also construct $\mN_f$ or $\mN_w$ by starting with  spheres that each intersect $\mG$ twice. In general $\mN_f$ is  not obviously equal $\mN_g$ and $\mN_f\cap\mN_w\neq\emptyset$, which essentially is the cause of our difficulties. \end{remarks}  

\begin{lemma} \label{common duals} Let $(\mG,\mR,\mF,\mW)$ be a $F|W$ system such that the boundary germs of $\mF$ coincide with that of $\mW$. 

 i) If $\mN_f=\mN_w$, then $\phi(\mG,\mR,\mF,\mW)$ is S-equivalent to $\id$.

ii) If $\mW'$ is another set of Whitney discs with $\mN_w=\mN_{w'}$, then $\phi(\mG,\mR,\mF,\mW)$ is S-equivalent to $\phi(\mG,\mR,\mF,\mW')$.\end{lemma}

\begin{proof}  To prove i) we show that $\tilde\mW$ interpolates to $\tilde\mF$.  In $\tilde V_k$ attempt to construct a interpolating system $\tilde \mW'$ by using $\tilde\mW$ near the $-\infty$ end and $\tilde \mF$ near the $+\infty$ end.  While each component of $\tilde \mW'$ is embedded a component coming from $\tilde \mW$ may intersect one from $\tilde \mF$.  Since $\tilde \mN_f$ is a common system of dual spheres  to both $\tilde \mF$ and $\tilde \mW$ we can use copies of components of $\tilde \mN_f$ to tube away these intersections.  In a similar manner  prove ii) by showing that $\mW$ interpolates to $\mW'$.\end{proof}

\begin{question} \label{second test case} (Second Test Case) Let $(\mG,\mR,\mF,\mW)$ be such that $\{\mG\}=G_1$, $R_1$ has exactly two fingers respectively of winding $\pm 1$ and each Whitney disc coincides with a finger disc in a neighborhood of its boundary.  Is $\phi(\mG,\mR,\mF,\mW)$ S-equivalent to $\id$?\end{question}

\begin{problem} (Slice missing slice disc problem).  The knot $K\subset S^1\times S^2=\partial S^1\times B^3$  shown in Figure \ref{Cloverleaf1} bounds two obvious ribbon discs $D_1$ and $D_2$ such that the simple closed curve $\alpha\subset S^1\times S^2\setminus K$ (resp. $\beta$) slices in $S^1\times B^3$ with a slice disc disjoint from $D_1$ (resp. $D_2$).  Is it true that for any smooth disc $D$ bounded by $K$, one of $\alpha$ or $\beta$ slices in the complement of $D$?\end{problem}

\setlength{\tabcolsep}{00pt}
\begin{figure}
 \centering
\begin{tabular}{ c c }
 $\includegraphics[width=2.0in]{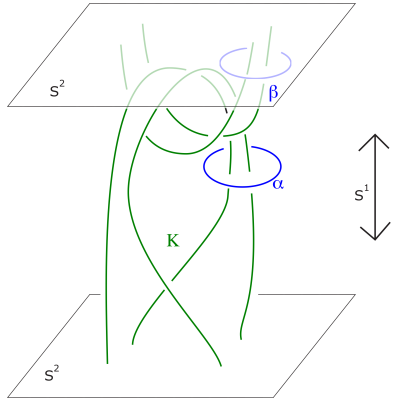}$  
\end{tabular}
 \caption[(a) X; (b) Y]{\label{Cloverleaf1}
 \begin{tabular}[t]{ @{} r @{\ } l @{}}
Slice Missing Slice Disc
\end{tabular}}
\end{figure}

\section   {Twisted Whitney discs} \label{twisted}

The goal of this section and the next is to show that a $F|W$ system interpolates to one whose  finger and Whitney discs coincide along neighborhoods of their boundaries, i.e. their \emph{boundary germs} coincide.   Such a $F|W$ system is called \emph{boundary germ coinciding}.  To do this we first arrange for their boundaries to coincide and then get neighborhoods of their boundaries to coincide. The  proof of the latter uses the main result of this section, Lemma \ref{twist}  which subject to a certain technical condition that is satisfied by first passing to a finite cover, asserts that given any system $\mF$ of Whitney discs for $\mG$ and $\mR$ in $\sonesthree\times \#_k\stwostwo$, there exists another system $\mF'$ such that $\mF'$ has prescribed twisting relative to that of $\mF$ and $\phi(\mG,\mR,\mF,\mF')$ is S-equivalent to the identity.  

We now define the twisting of one Whitney disc relative to another, by putting a neighborhood of the boundary of one into a normal form relative to that of the other.  

\begin{definition} \label{whitney normal form} Let $w_0$ and $w_1$ be Whitney discs for the oriented, possibly disconnected surfaces $G$ and $R$  in the oriented 4-manifold $M$ such that $\partial w_0=\partial w_1$.  Let $x$ and $y$ denote the points of $w_i\cap G\cap R$ where $x$ is the point of $-1$ intersection.  Define $\beta=w_i\cap G$ and $\alpha =w_i\cap R$ with $\alpha$ (resp. $\beta$) oriented from $x$ to $y$ (resp. $y$ to $x$).   See Figure \ref{Figure2,1}  a) which shows a 3-dimensional slice of $M$ that contains $w_0$.  We assume that the orientation of $G$ is given by $(\epsilon_1,\epsilon_2)$ and at $y$, $R$ is oriented by $(\epsilon_3, \epsilon_4)$.  After an isotopy of $w_1$ we can assume that it coincides with $w_0$ near both $x$ and $y$.      It follows that after a further isotopy, a neighborhood of $\partial w_1$ rotates $p\in \BZ$ times along $\alpha$ and $q\in \BZ$ times along $\beta$.  Here $q$ (resp. $p$) is the number of full right hand twists about $\beta$ (resp. $\alpha$), the former using the convention of Figure \ref{Figure2,1} b).  Figure \ref{Figure2,1} c) shows the projection to the 3-dimensional slice of a neighborhood of the boundary of a $(3,1)$-twisted disc to the $(x,y,z)$ plane.  We call $w_1$ a $(p,q)$-\emph{twisted Whitney disc} rel $w_0$ and we define $\tw(w_1,w_0)=(p,q)$.  \end{definition}

\setlength{\tabcolsep}{60pt}
\begin{figure}
 \centering
\begin{tabular}{ c c }
 $\includegraphics[width=5.5in]{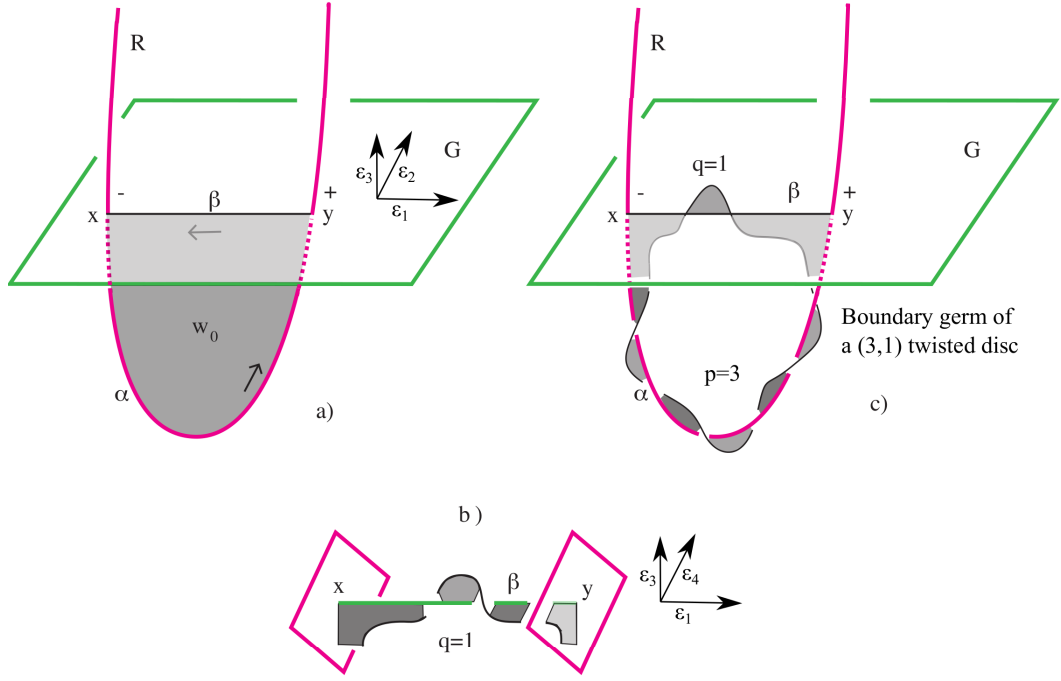}$  
\end{tabular}
 \caption[(a) X; (b) Y]{\label{Figure2,1}
 \begin{tabular}[t]{ @{} r @{\ } l @{}}
The boundary germ of a twisted Whitney disc
\end{tabular}}
\end{figure}

\begin{lemma}  \label{twist sum} If $w_0, w_1$, and $w_2$ are Whitney discs with $\partial w_2=\partial w_1=\partial w_0$, then $\tw(w_2,w_0)=\tw(w_2,w_1)+\tw(w_1, w_0)$.  Also $\tw(w_0,w_1)=-\tw(w_1,w_0)$.  \qed\end{lemma}

\begin{lemma}  \label{parity} If $w_1\subset M$ is a $(p,q)$-twisted Whitney disc rel $w_0$ and $M$ has trivial second Stiefel - Whitney class, then $p+q$ is even.\end{lemma}

\begin{proof} The normal bundle of the Whitney disc $w_1$ for surfaces $R$ and $G$ has a framing which when restricted to $\partial D$ has, where applicable, one vector tangent to $R$ and transverse to $G$ and the other tangent to $G$ and transverse to $R$.  On the other hand, starting with $w_0$ perform the boundary twisting operation, e.g. see \cite{E} P. 216, to obtain an embedded pre-Whitney disc $E$ whose boundary has a neighborhood that coincides with that of $w_1$ but whose framing differs from that of a genuine Whitney disc by $p+q$ mod 2.  By gluing $E$ to $w_1$ along their 
boundaries and smoothing near the gluing we obtain a smoothly immersed 2-sphere in $M$ with normal bundle of Euler class equal to $p+q$ mod 2.  Since $w_2(M)=0$ it follows that $p+q =0$ mod 2.  \end{proof}

The following is the main result of this section.  As before $V_k$ denotes $\sonesthree\times \#_k\stwostwo$ and $\mG$ and $\mR$ denote sets of algebraically dual embedded 2-spheres as in Definition \ref{fw}.  

\begin{lemma} \label{twist}  Let $\mF=\{f_1, \cdots, f_n\}$ be a complete set of Whitney discs for $\mG=\{G_1, \cdots, G_k\}$ and $\mR=\{R_1, \cdots, R_k\}$ in $V_k$ such that the winding of any hand from any $R_i$ to $G_i$ is equal to zero.     Let $((p_1,q_1), \cdots, (p_n,q_n))\in \BZ\oplus\BZ$ be such that for all $i, p_i+q_i$ is even.  Then there exists a system of Whitney discs $\mF'=\{f_1', \cdots, f_n'\}$ such that $\partial \mF=\partial \mF', \phi(\mG,\mR, \mF, \mF')$ is S-equivalent to the identity and for $i=1,2, \cdots, n, \tw(f_i',f_i)=(p_i,q_i)$. \end{lemma}

\begin{remark}   Let $z_r\in G_r\cap R_r$ denote the point disjoint from all the $f_i$'s.  The winding hypothesis implies that if there is a Whitney disc $f_i$ between $G_r$ and $R_r$, then the element of $\pi_1(V_k)$ corresponding to a loop starting at $z_r$ that follows $R_r$ to a point of $G_r\cap R_r\cap f_i$, then follows $G_r$ back to $z_r$ is homotopically trivial.\end{remark}

\begin{proof}  \emph{Step 1}: The lemma holds when $(p_i,q_i)=(0,0)$ for all $i> 1$ and either $(p_1,q_1)=(0,\pm 2)$ or $(\pm 2, 0)$.

\vskip 8pt

\noindent\emph{Proof of Step 1.}   Suppose that $f_1$ cancels points of $G_r\cap R_s$. We consider the $(p_1,q_1)=(2,0)$ case as the other cases are similar.  Use the boundary twisting operation \cite{E}) to obtain an embedded $(2,0)$ twisted pre-Whitney disc $E_0$ which fails to be a genuine Whitney disc because its framing is off by two and it intersects $R_s$ twice.  Correct the framing by replacing a small disc with one of self intersection $\pm 1$ (see \cite{FQ} p. 14) and then push the intersection off the $\alpha$ boundary to obtain the embedded disc $E_1$ which has the correct framing but $|E_1\cap R_s|=4$.

  Let $R_s'$ be a geometrically dual sphere to $R_s$ disjoint from $(\mG\cup (R\setminus R_s)\cup \mF\cup E_1)$.  Such a sphere can be obtained from $G_s$ by doing Whitney moves to $G_s$ using the components of $\mF$ that intersect $G_s$ and then isotoping slightly.  Here we assume that  $E_1$ is constructed to lie very close to $f_1$.  Next eliminate the four $\inte(E_1)\cap R_s$ intersections by tubing $E_1$ to four parallel copies of $R_s'$ along arcs in $R_s$  disjoint from $\mF$.  Thus we obtain a Whitney disc $f_1'$ with $\tw(f_1',f_1)=(2,0)$ and $f_1'\cap (\mF\setminus f_1)=\emptyset$.  Let $\mF_1=\{f_1', f_2, \cdots, f_n\}$.  Now $(\tilde G, \tilde R, \tilde F, \tilde F_1)\subset \tilde V_k$ interpolates to $(\tilde G, \tilde R, \tilde F, \tilde F)$ by replacing lifts of $f_1'$ with lifts of $f_1$ near the $+\infty$ end of $\tilde V_k$.  This uses the fact that a lift of $f_1'$ is disjoint from $\tilde F$ except for the single component that has the same boundary.  It follows that $\phi(\mG, \mR, \mF, \mF_1)$ is $S$-equivalent to $\phi(\mG, \mR, \mF, \mF)$ which is the class of the identity.\qed

\vskip 8pt
\noindent\emph{Step 2}:  If the Lemma holds for $(p_1,q_1) = (p,q)$ and $(p_i,q_i)=(0,0)$ for $i>1$, then the Lemma holds for $(p_1,q_1)=(p \pm 2, q)$ and $(p, q\pm 2)$ with $(p_i,q_i)=(0,0)$ for $i>1$.

\vskip 8pt\
\noindent\emph{Proof of Step 2.}  Let $\mF' $ be a set of  Whitney discs for which the hypothesis of Step 2 holds.  Now apply Step 1 to  $\mF'$ to obtain  $\mF_1$.  By Lemma \ref{twist sum} $\mF_1$ satisfies the twisting  conclusion of Step 2 relative to $\mF$.  Also by  Lemma \ref{factorization}, $\phi(\mG,\mR,\mF,\mF_1)=\phi(\mG,\mR,\mF', \mF_1)\circ\phi\mG,\mR,\mF, \mF')$ and so is $S$-equivalent to the identity.\qed

\vskip 8pt

\noindent\emph{Step 3:} The lemma holds when $(p_1,q_1)=(1,1), (p_i,q_i)=(0,0)$ for $i>1$  and $f_1 \cap G_r\cap R_s\neq\emptyset$ where $r\neq s$.

\vskip 8pt

\noindent\emph{Proof of Step 3.}  The proof is a modification of the proof of Step 1. 
Construct a framed, embedded pre-Whitney disc $E_0$ with $\tw(E_0,f_1)=(1,1)$ disjoint from $G_r'\cup R_s'$ by starting with $f_1$ and then doing boundary twisting operations by twisting  once about $f_1\cap G_r$ and once about $f_1\cap R_s$.
If necessary, correct the framing by first replacing a small embedded disc by one with self intersection $\pm 1$ and then making the resulting pre-Whitney disc embedded  by pushing the self intersection off the $f_1\cap R_s$ component.  Except for its intersections with $G_r$ and $R_s$ the resulting disc $E_1$ is a genuine Whitney disc. 

Construct geometric dual spheres $R_s'$ (resp. $G_r'$) from $G_s$ (resp. $R_r$) disjoint from $(\mG\cup (\mR\setminus R_s)\cup \mF\cup E_1)$ (resp. $(\mG\setminus G_r)\cup\mR \cup \mF\cup E_1)$. Construct $f_1'$ by tubing off $E_1\cap(R_s\cup G_r)$ using copies of $R_s'$ and $G_r'$ and tubes that avoid the discs of $\mF$.  The argument of Step 1 shows that if $\mF'=\{f_1', f_2, \cdots, f_n\}$, then $f(\mG,\mR, \mF, \mF')$ is $S$-equivalent to the identity.\qed

\vskip 8pt
\noindent\emph{Step 4}: The lemma holds when $(p_1,q_1)=(1,1), (p_i,q_i)=(0,0)$ for $i>1$  and $f_1 \cap G_1\cap R_1\neq\emptyset$.
\vskip 8pt

\setlength{\tabcolsep}{00pt}
\begin{figure}
 \centering
\begin{tabular}{ c c }
 $\includegraphics[width=5.0in]{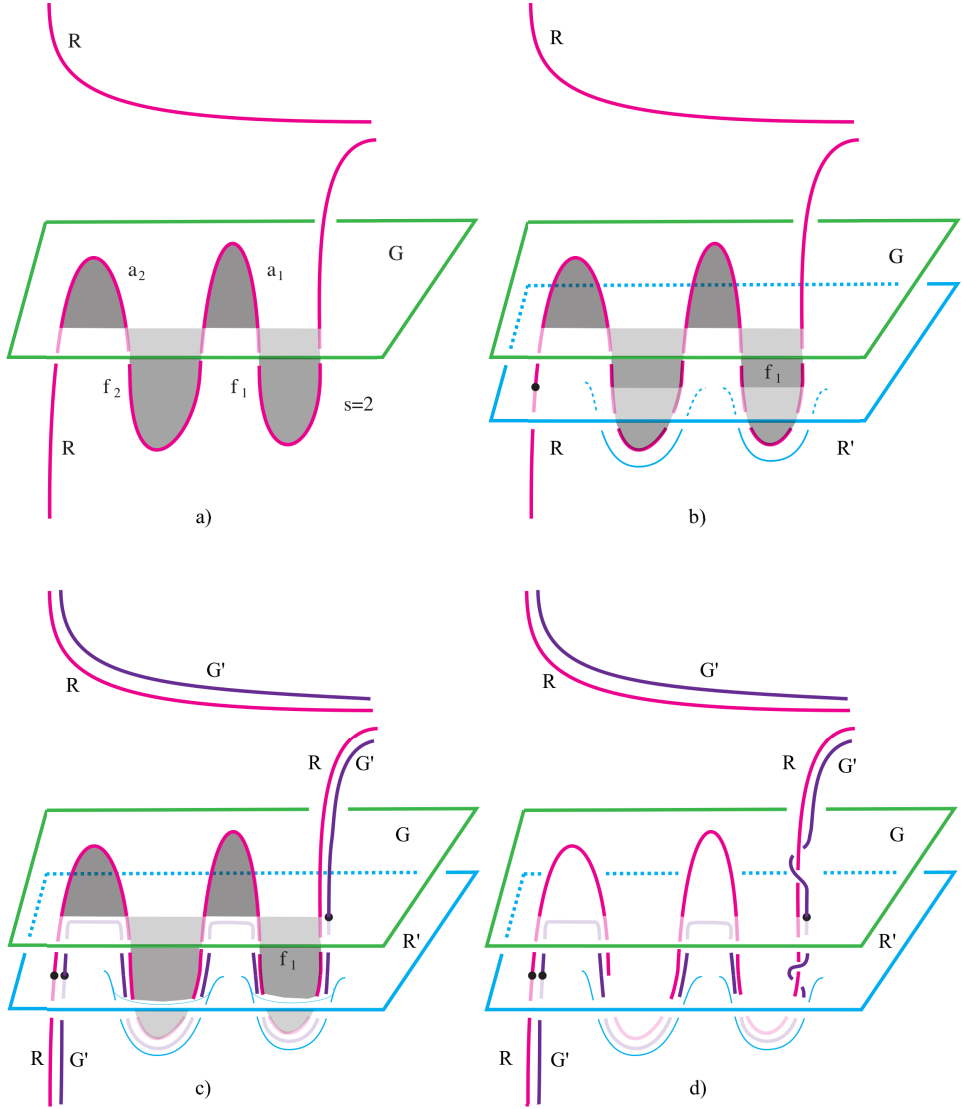}$  
\end{tabular}
 \caption[(a) X; (b) Y]{\label{Figure2,2}
 \begin{tabular}[t]{ @{} r @{\ } l @{}}
Constructing dual spheres
\end{tabular}}
\end{figure}

\noindent\emph{Proof of Step 4.}  As in the previous steps it suffices to show that there exists a Whitney disc $f_1'$ with $\partial f_1'=\partial f_1$ and $\tw(f_1',f_1)=(1,1)$ such that for $i>1$, $f_1'\cap f_i=\emptyset$ and in $\tilde V_k$, if $\tilde f_1'$ and $\tilde f_1$ are lifts of $f_1'$ and $f_1$, then $\tilde f_1'\cap \tilde f_1\neq\emptyset$ if and only if $\partial \tilde f_1'=\partial \tilde f_1$.   

To simplify notation denote $G_1$ and $R_1$ by $G$ and $R$.  By Lemma \ref{arm} we can assume that $\mR$ is in finger hand form with respect to $\mG$. See Figure \ref{hand}.  Since the winding of the hand of $R$  containing $f_1$ equals $0$, it follows that if this hand has $s\ge 1$ fingers, then there are $s$ alternative Whitney discs $a_1, \cdots, a_s$ as in Figure \ref{Figure2,2} a), where the $s=2$ case is shown.  Note that by choosing the $a_i$'s appropriately we can assume that the hand's fingers and alternate discs appear as in the figure, in particular $f_1$ is the first in the indicated sequence of $a_i$'s and $f_j$'s. 

Now construct the geometrically dual sphere $R'$ for $R$.  As before construct $R'$ by applying Whitney moves to $G$ using the applicable discs of $\mF$.  After a slight isotopy $R'\cap (\mG\cup \mF)=\emptyset$ and  $|R'\cap \mR|=|R'\cap R|=1$.  See Figure \ref{Figure2,2} b).  Next construct the geometrically dual sphere $G'$ for $G$ by first doing Whitney moves to $R$ along $a_1, \cdots a_s$  as well as Whitney moves to all the relevant discs in $\mF\setminus \{f_1,\cdots, f_s\}$.  Construct $G'$ to lie close to $R$ and these discs, in particular $G'\cap R'$ is very close to $R\cap R'$.  After a slight isotopy we can assume that $|G'\cap \mG|=|G'\cap G|=1$ and $G'\cap (\mR\cup \mF)=\emptyset$.  See Figure \ref{Figure2,2} c).

Next isotope $G'$ so that its intersection with a 3-ball containing $f_1$ has a full right hand twist about $R$ near $f_1$ which is compensated by a full left hand twist as in Figure \ref{Figure2,2} d). This twist gets undone when moving both in the past and in the future and creates an intersection with $f_1$.
 
\setlength{\tabcolsep}{60pt}
\begin{figure}
 \centering
\begin{tabular}{ c c }
 $\includegraphics[width=5.5in]{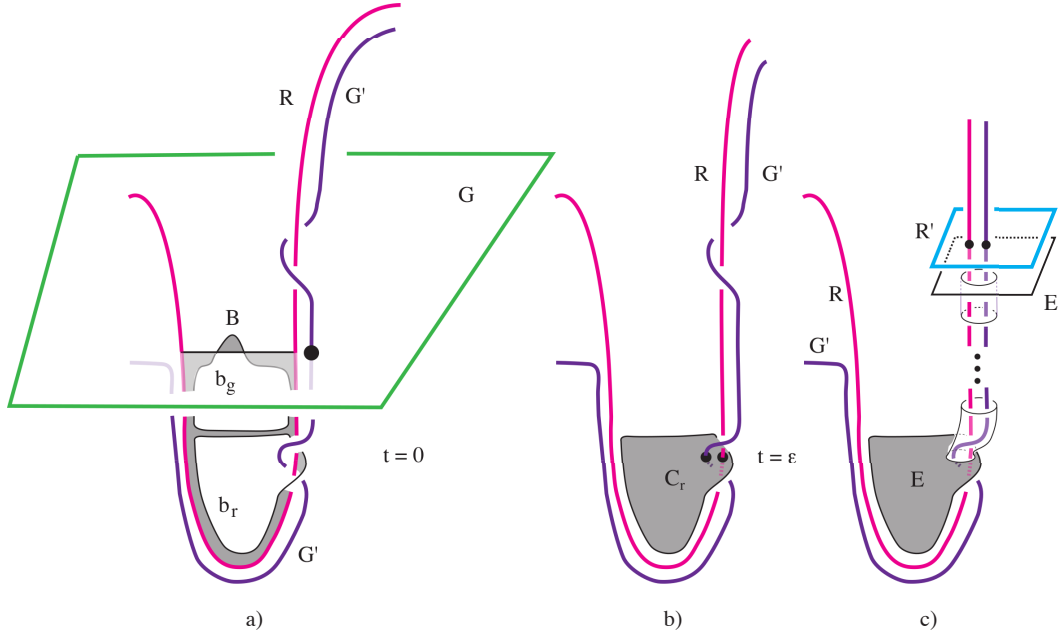}$  
\end{tabular}
 \caption[(a) X; (b) Y]{\label{Figure2,3}
 \begin{tabular}[t]{ @{} r @{\ } l @{}}
Constructing a (1,1)-twisted Whitney disc
\end{tabular}}
\end{figure}

We now construct $f_1'$.  Let $A$ be an annulus having $(1,1)$-twisting near $\partial f_1$ with one component of $\partial A$ equal to $\partial f_1$.  By construction, $A\cap G'=\emptyset$.   Next add an untwisted band to $A$ to create $B$ as in Figure \ref{Figure2,3} a).   Again $B\cap G'=\emptyset$.   Let $b_r$ and $b_g$ denote the components of $\partial B$ that link around $R$ and $G$.  Construct an embedded correctly framed pre-Whitney disc $C$ by capping off these components with discs $C_r$ and $C_g$.  $C_r$ consists of an annulus starting at $b_r $ that moves directly into the future  which is then capped off by a disc in a horizontal time slice.  So if Figure \ref{Figure2,3} a) shows a $t=0$ slice, then Figure \ref{Figure2,3} b) shows the $t=\epsilon$ slice and the spanning disc.  The disc $C_g$ bounded by $b_g $ will have an excess intersection point with $G$.  Since $1+1=2$, for $C$ to have the correct framing  the construction of $C_g$ might require first creating a self intersection and then pushing off to create two extra intersections with $G$.  Note that the unique intersection point of $C\cap G'$ is next to the unique transverse intersection point of $C$ with $R$ as indicated in Figure \ref{Figure2,3} b).  Next using a tube that follows parallel arcs in $R$ and $G'$, tube off these two points using a copy of $R'$.    See Figure \ref{Figure2,3} c).  Let $E$ denote the resulting framed embedded pre-Whitney disc.  Note that $E\cap G'=\emptyset$.  Finally tube off the one (or three) excess points of $E\cap G$ with copies of $G'$ to construct the desired Whitney disc $f_1'$.  \qed

\noindent\emph{Step 5}: The lemma holds when $(p_i,q_i)=(0,0)$ for $ i\neq 1$.
\vskip 8pt

\noindent\emph{Proof of Step 5}  First assume that $p_1$ is even, and hence by parity so is $q_1$.  Then there exists a sequence $(0,0)=(p^1,q^1), (2,0)=(p^2, q^2),\cdots, (p^r,q^r)=(p_1,q_1) $ so that for $1\le m\le r-1, (p^{m+1},q^{m+1})-(p^m,q^m)=(\pm 2, 0)$ or $(0,\pm 2)$.  Let $\mF=\mF^1, \mF^2, \cdots, \mF^r$ be complete systems of Whitney discs for $\mR$ and $\mG$ such that for $1\le m\le r, \mF^m=\{f^m, f_2, \cdots, f_n\} $ where $f^1= f_1, \tw(f^{2}, f^1)=(2,0)$, for $2\le m\le r-1, \tw(f^{m+1}, f^m)=((p^{m+1},q^{m+1})-(p^m,q^m))$, and for all $i, \phi(\mG, \mR, \mF^i, \mF^{i+1})$ is $S$-equivalent to the identity.  The existence of $\mF^2$ follows by Step 1 and the existence of $\mF^i, i>2$ follows from Step 2.  By Lemma \ref{twist sum} $\tw(f^r, f_1)=(p_1,q_1)$ and by factorization $\phi(\mG,\mR, \mF, \mF^r)=\phi(\mG,\mR, \mF^{r-1}, \mF^r)\circ \phi(\mG,\mR, \mF^{r-2}, \mF^{r-1})\circ\cdots\circ \phi(\mG,\mR, \mF^1, \mF^2)$.  Since each of the factors is $S$-equivalent to the $ \id$, the result follows.

Next assume that $p_1$ is odd.  Then there exists a sequence $(0,0)=(p^1,q^1), (1,1)=(p^2, q^2),\cdots, (p^r,q^r)=(p_1,q_1)$ so that for $2\le m\le r-1$, $(p^{m+1},q^{m+1})-(p^m,q^m)=(\pm 2, 0)$ or $(0,\pm 2)$.  Argue as above, except use one of Steps 3 or 4 in place of Step 1.\qed

\vskip 8pt

\noindent\emph{Step 6}:  The lemma holds when $(p_i,q_i)=0$ for $i\neq s$.\qed

\vskip 8pt
\noindent\emph{Proof of Step 6}:  This follows from Step 5 after reordering the elements of $\mF$.\qed
\vskip 8pt
\noindent\emph{Step 7}: General Case
\vskip 8pt \emph{Proof of Step 7:} Through repeated uses of Step 6, Lemma \ref{twist sum} and factorization, inductively construct $\mF_0=\mF, \mF_1, \cdots, \mF_n=\mF' $ so that $\phi(\mG,\mR,\mF,\mF_i)$ is $S$-equivalent to $\id$ and if $\mF_i=(f_1^i, \cdots, f_n^i)$, then $\tw(f_j^i,f_j)= (p_j,q_j)$ if $j\le i$ and $(0,0)$ otherwise.\end{proof}

\section {Germs of finger and Whitney discs}

The goal of this section is to prove the following:

\begin{proposition}  \label{germs} Let $\phi\in \Diff_0(\sonesthree)$, then $\phi$ is $S$-equivalent to $\psi\in \Diff_0(\sonesthree)$ arising from a boundary germ coinciding  $F|W$ system.  \end{proposition}

\begin{lemma}  \label{g boundary} Given $(\mG,\mR,\mF,\mW)$ there exists $(\mG,\mR,\mF,\mW')$ such that $\phi(\mG,\mR,\mF,\mW')$ is $S$-equivalent to $\id$ and $\mW\cap \mG=\mW'\cap \mG$.\end{lemma}

\begin{proof}  View $\mG$ and $\mR$ in AHF form. See Figure \ref{hand} which shows the restriction of $\mR\cup\mG$ to the i'th $\stwostwo$ factor, the finger discs which intersect $G_i$ as well as spheres $R_i' $ and $G_i'$ respectively geometrically dual to $R_i$ and $G_i$.  Letting $\mR'=\cup_{R_i\in \mR} R_i'$, note that $\mR'\cap \mG=\emptyset$ and is geometrically dual to $\mR$.   Let $w\in \mW$.   If $w\cap G_i\cap R_j\neq \emptyset$ and $\beta_w:=w\cap G_i$, then since $\partial w$ is homotopically trivial in $V_k$, $\partial \beta_w$ lies in a single hand $H_w\subset R_j$.  Let $\gamma_w$ be an embedded loop of the form $\beta_w\cup\alpha_w$ where  $\alpha_w \subset H_w$ and has interior disjoint from the \emph{finger arcs} $\mF\cap\mR$.  The loops $\cup_{w\in \mW} \gamma_w$ can be chosen to be pairwise disjoint.  Each $\gamma_w$ bounds an immersed disc $D_w$ contained in it's $\stwostwo$ factor whose interior is disjoint from $\mG\cup \mR'$.  These discs may have intersections and self intersections and possibly $\inte(D_w)\cap \mR\neq\emptyset$.  At the cost of creating additional intersections with $\mR$, these discs, which continue to be called $D_w$'s, can be made embedded and pairwise disjoint.  By boundary twisting near the $\alpha_w$'s, thereby creating further intersections with $\mR$, these discs can be  correctly framed.  They may fail to be Whitney discs only because the interior of the $D_w$'s intersect $\mR$ transversely.  Finally eliminate these intersections by tubing with parallel copies of components $\mR'$ using tubes that follow paths in $\mR$. Let $\mW'$ denote the resulting collection of Whitney discs.  Since $\mR'\cap(\mW'\cup \mF)=\emptyset$, it follows that $\phi(\mG,\mR,\mF, \mW')$ is $S$-equivalent to $\id$.  \end{proof}

\begin{lemma}  \label{r boundary} Given $(\mG,\mR,\mF,\mW)$ such that $\mF\cap\mG=\mW\cap\mG$ there exists $(\mG,\mR,\mF,\mW'')$ such that $\partial \mW''=\partial \mW$ and $\phi(\mG,\mR,\mF,\mW'')$ is $S$-equivalent to the identity.\end{lemma}

\begin{proof}  View $\mG$ in AHF form with respect to $\mF$  and let $\mG'$ denote the geometrically dual spheres to $\mG$ as above disjoint from $\mR$.    If $w\in \mW$, then since $w\cap \mG\subset\mF\cap \mG$, i.e. is a finger arc, $\partial w$ bounds an immersed disc whose interior is disjoint from $\mR \cap \mG'$.  As in the previous lemma these discs can be modified to construct a family of Whitney discs $\mW''$ with $\partial \mW''=\partial \mW$ and $\mW''\cap \mG'=\emptyset$ and hence $\phi(\mG,\mR, \mF, \mW'')$ is $S$-equivalent to the identity.  \end{proof}

\begin{lemma} \label{boundary}  Given $(\mG,\mR,\mF,\mW)$ there exists $(\mG,\mR,\mF_1,\mW_1)$ such that $\phi(\mG,\mR,\mF,\mW)$ is $S$-equivalent to $\phi(\mG,\mR, \mF_1, \mW_1)$ and $\partial \mF_1=\partial \mW_1$.\end{lemma}

\begin{proof} Apply Lemma \ref{g boundary} to find $(\mG,\mR,\mF,\mW')$ such that $\mW'\cap \mG=\mW\cap \mG$ and $\phi(\mG,\mR,\mF, \mW')$ is $S$-equivalent to the $\id$.  By  Lemma \ref{factorization}, $\phi(\mG, \mR, \mF, \mW)$ is isotopic to $\phi(\mG, \mR, \mW', \mW)\circ \phi(\mG, \mR, \mF, \mW')$ and hence $\phi(\mG, \mR, \mF, \mW)$ is $S$-equivalent to $\phi(\mG,\mR, \mW', \mW)$.

Now apply Lemma \ref{r boundary} to find $(\mG, \mR, \mW', \mW'')$ such that $\phi(\mG, \mR, \mW', \mW'')$ is $S$-equivalent to $\id$ and $\partial \mW=\partial \mW''$.  Again by factorization  $\phi(\mG,\mR, \mW', \mW)$ is $S$-equivalent to $\phi(\mG, \mR, \mW'', \mW)\circ \phi(\mG, \mR, \mW', \mW'')$.  It follows that $\phi(\mG,\mR,\mF,\mW)$ is $S$-equivalent to $\phi(\mG,\mR, \mW'', \mW)$  where $\partial \mW=\partial\mW''$.\end{proof}

\vskip 8pt

\noindent\emph{Proof of Proposition} \ref{germs}.  Let $\phi\in\Diff_0(\sonesthree)$.  Suppose that it is represented by $(\mG,\mR,\mF,\mW)$ which is supported on $V_k$.  Staying within the $S$-equivalence class of $\phi$, we can additionally assume by Lemma \ref{boundary} that $\partial \mF=\partial \mW$.  By passing to a finite sheeted cover of $\sonesthree$ and hence $V_k$, an operation preserving $S$-equivalence, we can assume that every hand from an $R_i$ to its $G_i$ has winding $0$.  Indeed, any finite cover of degree greater than the maximal winding over all hands from an $R_i$ to its $G_i$ suffices.  Order the elements $(f_1, \cdots, f_n)$ of $\mF$ and $(w_1, \cdots, w_n)$ of $\mW$ so that for each $i$, $\partial f_i=\partial w_i$.    Let $(p_i,q_i)=\tw(w_i,f_i)$. Now apply Lemma \ref{twist} to find $\mF'$ so that  $\phi(\mG, \mR, \mF, \mF')$ is $S$-equivalent to the identity and $\tw(\mF',\mF)=((p_1,q_1), \cdots, (p_n,q_n))$.  It follows that $\phi (\mG,\mR,\mF,\mW)$ is $S$-equivalent to $\phi (\mG,\mR,\mF',\mW)\circ \phi (\mG, \mR, \mF, \mF')$ and hence to $\phi(\mG,\mR,\mF',\mW)$.  By Lemma \ref{twist sum}, $\tw(\mF',\mW)=\tw(\mF', \mF) + \tw(\mF, \mW) = \tw(\mF', \mF)-\tw(\mW, \mF)=((0,0), \cdots, (0,0))$.\qed

It follows from Theorem \ref{main sequivalent} and Proposition \ref{germs}

\begin{theorem}  \label{germs theorem} The Schoenflies conjecture is true if and only if every boundary germ coinciding $F|W$ system interpolates to the trivial $F|W$ system.  \end{theorem}

\section{Homotopic Whitney and Finger Discs}

The following is the main result of this section.

\begin{proposition} \label{whitney homotopy} If $\phi \in \Diff_0(\sonesthree)$, then up to $S$-equivalence $\phi$ is represented by an $F|W$ system on some $V_k$ such that if $\mF=(f_1,\cdots, f_m)$ and $\mW=(w_1, \cdots, w_m)$, then for every $i$ the boundary germ of $w_i$ coincides with that of $f_i$ and $w_i$ is homotopic to $f_i$ via a homotopy fixing $N(\partial w_i)$ pointwise and supported in $V_k\setminus \mG\cup \mR$.\end{proposition}

\begin{remarks} \label{hands unique} Using Proposition \ref{germs} we will start with an $F|W$ system  $\mF=(f_1,\cdots, f_m)$ and $\mW=(w_1, \cdots, w_m)$ satisfying the boundary germ conclusion.   We will assume that $\mR$ is in AHF form with respect to $\mF$.  Also for $i\neq j, R_i\cap G_j$ is contained in a single hand, since this can be achieved by passing to a finite cover, an $S$-equivalence preserving operation.  In addition, we will assume that for all $i$, $R_i\cap G_{i+1}\neq\emptyset$.  If necessary, achieve this by adding extra hands with equal  finger and Whitney discs.  Here indices in $\BZ$ are modulo $k$.  \end{remarks} 

\begin{notation}  We continue to use notation  as in \ref{winding}.  Define $X:=N(\mG\cup\mR)$, $Y:=V_k-\inte X$, $U_0=\cup_{i=1}^k S^2\times S^2_i$ and $Z:=V_k\setminus \cup_{i=1}^k \inte(S^2\times S^2_i)$.  Let $A_{i,j}$ denote the arm that goes from $R_i$ to $G_j$, provided one exists and define $A_{i,j}^Z:=Z\cap A_{i,j}$.  \end{notation}

\begin{lemma} \label{high cover} By passing to another finite cover and isotoping the $A^Z_{i,j}$'s we can assume that $ \pi(A_{i,j}^Z)\subset$ exactly one of $(i,j)$ or $(j,i), \diam\pi(A_{i,j}^Z)<k/16$ and     for each $w_p\in \mW$, $\diam\pi(w_p\cap Z) <k/16$.  \qed\end{lemma}

\begin{notation}  Let $I_{i,j}$ denote the short subinterval of $S^1$ bounded by $i, j$.  In what follows this interval will usually have length $\le k/4$.  Define a directed graph $G$ whose vertices are the $\Sigma_i$'s and whose edges $E:=\{e_{i,j} \}$ are the $A_{i,j}$'s with $i\neq j$, where $e_{i,j}$ points from $\Sigma_i $ to $\Sigma_j$.  Let $C$ denote the cycle formed by $e_{0,1}, \cdots, e_{k-1,k}$.   
For $e_{i,j}\in E$, let $H_{i,j}\subset Z$ denote a 1-handle with attaching discs in $\Sigma_i$ and $\Sigma_j$ such that $\pi(H_{i,j})\subset I_{i,j}$.  Assume that the $H_{i,j}$'s are pairwise disjoint.  Let $U_1=U_0\cup_{i,j} H_{i,j}$.  \end{notation}  

\begin{lemma} The $A^Z_{i,j}$'s can be naturally isotoped into  the $H_{i,j}$'s after which $U_1$ deformation retracts to $X\cup N(\mF\cup\mA)$, where $\mA =\{a_1, \cdots, a_n\}$ are the auxiliary discs. \end{lemma}

\begin{proof}  Thicken the finger discs to 2-handles and expand them to first \emph{fill} the fingers, then the hands and then the arms.  Finally add $N(\mA)$.  The result is isotopic to $N(\mG^{\std}\cup\mR^{\std})$ union 1-handles, one for every arm from $R_i$ to $G_j, i\neq j$.  Note that an arm from an $\Rs_i$ to $\Gs_i$ together with its finger and auxiliary discs gets absorbed into $N(\mGs\cup \mRs)$. After the $A^Z_{i,j}$'s have been naturally isotoped to lie in the $H_{i,j}$'s we see that $U_1$ deformation retracts to $X\cup N(\mF\cup\mA)$.\end{proof}

\begin{lemma} There exists a system of discs $\mD:=\cup_{e\in E\setminus C}D_e:=\{D_1,\cdots, D_r\}$ called \emph{arm rest discs} such that

i) The elements of $\mD$ are pairwise disjoint and properly embedded in $Y$

ii) $\mD\cap(\mF\cup \mA)=\emptyset$

iii) If $U=X\cup N(\mF\cup\mA\cup \mD)$, then $S^1\times S^3\setminus \inte(U)$ is isotopic  to a vertical $S^1\times B^3$. \end{lemma}

\begin{proof} The arm rest discs will have the property that if $D_i$ corresponds to $e_{p,q}$, then $D_i$ runs over  $H_{p,q}$ exactly once and is disjoint from all the other $H_{p,q}$'s except those of the form $H_{i, i+1}$. Assuming conclusions i) and ii), it follows that when thickened to 2-handles the $D_i$'s cancel the $H_{i,j}$'s of $U_1$ with $ j\neq i+1$ and hence the resulting $U$ is isotopic to $\cup_{i=1}^k S^2\times S^2_i\cup H_{i,i+1}$ and so its closed complement is isotopically a vertical $S^1\times B^3$. We detail the construction of a special case, say $D_1:=D_{e_{3,1}}$, from which the general construction may be deduced.  Our $D_1$ is bounded by $ \alpha*\beta$ where $\alpha\subset \partial N(A_{3,1})$ with initial point in $N(\Rs_3)$ and final point in  $N(A_{3,1})\cap N(\Gs_1)$.  Then $\beta$  follows a path  $\subset\partial N(\Gs_1\cup \Rs_1\cup A_{1,2}\cup \Gs_2\cup \Rs_2\cup A_{2,3}\cup\Gs_3\cup \Rs_3)$.  \end{proof}

\begin{definition}  Let $\mE:=\mF\cup\mA\cup\mD$.  Fix orientations on the elements of $\mE$ and then induce orientations on the elements of $\mW$ from those of $\mF$ using the boundary germ condition.   \end{definition}

\begin{lemma} \label{zero homology} If $S$ closed oriented surface in $Y$, then $[S]=0\in H_2(Y)$ if and only if for each $E\subset \mE, \langle S,E\rangle=0$.\end{lemma}

\begin{proof} The forward direction follows from the fact that algebraic intersection number is a homological invariant.  Conversely, if all the intersection numbers equal $0$, then $S$ is homologous to a surface disjoint from $U$ and hence is homologically trivial, since $H_2(S^1\times B^3)=0$.\end{proof} 

\begin{definition}  For $i\in \{1, \cdots, k\}$, let $R_i'\subset S^2\times S^2_i$ be a dual sphere to $R_i$ as in Figure \ref{hand}, i.e. is  obtained by choosing a parallel copy of $\Gs_i$ that intersects $R_i$ exactly once and $R_i'\cap \mF=\emptyset$.  Construct oriented pairwise disjoint \emph{linking spheres} to the elements of $\mE$, i.e. to $E\in \mE$ we define an oriented embedded sphere $T_E$ with trivial normal bundle such that $T_E\cap (\mE\cup\mG\cup\mR)=T_E\cap E$ is a single point of  positive sign.  We do this first on the finger and auxiliary discs.  Let $E$ be such a disc with say $E\cap R_j\neq\emptyset$.  Let $T'_E$ be a 2-sphere consisting of an annulus disjoint from $R_j$ that intersects $E$ once in its interior together with two discs that each intersect $R_j$ once of opposite sign.  These $T'_E$'s can be chosen to be oriented, pairwise disjoint and so that $T'_E\cap \mE=T'_E\cap E = 1$ positive point.  To obtain $T_E$ tube off the two intersections with parallel copies of $R_j'$ where the tubes follow arcs in $R_j$ \cite{No}.  This can be done maintaining pairwise disjointness and so that each $T_E$  is disjoint from all the $R'_j$'s.     Given an arm rest disc $D$ that intersects $H_{i,j}, j\neq i+1$, construct the oriented linking sphere $T_D\subset Z$ using a sphere that links $H_{i,j}$.   \end{definition} 

\begin{remark} \label{diameter} Note that each $T_E$ can be constructed so that $\diam(\pi(T_E\cap Z))<k/16$.\end{remark}

\begin{lemma} \label{homology} i) $H_2(Y)$ is freely generated by $\{[T_f]|f \in \mF \}\cup \{[T_a]|a\in \mA\}\cup\{[T_D]|D\in \mD\}$.

ii) $H_2(Y,\partial \mE)$  is freely generated by $\{[E]|E\in \mE\}$ and $\{[T_f]|f\in \mF\}\cup \{T_a|a\in \mA\}\cup\{T_D|D\in \mD\}.$ \qed\end{lemma}

\begin{lemma} If $w_i\in\mW$, then the coefficient of the $[T_{f_i}]$ term of $[w_i]\in H_2(Y, \partial \mE)$ equals zero.  \end{lemma}

\begin{proof} If this coefficient was non zero, then $w_i$ would have a framing inconsistent with that of a Whitney disc with the same boundary germ as $f_i$.\end{proof}

\begin{lemma} \label{homology pairs} If $n_{p,q}$ denotes the coefficient of the $[T_{f_q}]$ term of $[w_p] \in H_2(Y,\partial \mE)$, then $n_{i,j}=-n_{j,i}$.  \end{lemma}

\begin{proof}  Since $w_i\cap w_j=\emptyset$ it follows that $0=\langle [w_i], [w_j] \rangle=n_{i,j} + n_{j,i}$.  The latter equality follows since $[w_i]=[f_i]+n_{i,j} T_{f_j}$ and $ [w_j]=[f_j]+ n_{j,i} T_{f_i}$ plus other terms that do not contribute to interection number.\end{proof}

\begin{lemma} \label{whitney homology} There exists a system $\mF'=(f_1', \cdots, f_m')$ of Whitney discs with the same boundary germs as $\mF$ such that for each $p$, $ [f'_p]=[w_p]\in H_2(Y, \partial \mE)$ and $\phi(\mG, \mR, \mF, \mF')$ is $S$-equivalent to $\id$.  Finally, for each $p$, $\diam(\pi(f'_p\cap Z) ) \le k/4$. \end{lemma}

\begin{proof}  Construct $f^1_p$ by  first tubing $f_p$ to $n_{p,q}$ copies of $T_{f_q}$, when $q<p$ and then taking the disjoint union with $n_{p,q}$ parallel copies of $T_{f_q} $ when $p<q$.   The $f^1_p$'s can be constructed to be embedded, disjoint from each $R'_j$.

Construct $f^2_1, f^2_2, \cdots, f^2_m$ as follows.  For $q>p$ tube the disc component of $f^1_p$ to its $T_{f_q}$ components using tubes that follow arcs in $f^1_q$ connecting oppositely oriented points of $f^1_p\cap f^1_q$.  This can be done so that the $f^2_p$'s are pairwise disjoint, embedded and disjoint from each $R'_j$.  

Finally obtain $f'_p$ by tubing $f^2_p$ to $\langle w_p, a_r\rangle$ parallel copies of $T_{a_r}$ and $\langle w_p, D_s \rangle$ parallel copies of $T_{D_s}$.  This whole construction  can be done so that the $f'_p$'s are pairwise disjoint, disjoint from each $R'_j$ and each $\diam(\pi(f'_p\cap Z))\le k/4$.   By construction $[f'_p]=[w_p]$ all $p$.   Since for all $j$,  $(\mF\cup \mF' )\cap R'_j=\emptyset$ it follows that $\phi(\mG,\mR, \mF, \mF')$ is  $S$-equivalent to $\id$.   \end{proof}

\noindent\emph{Proof of Proposition} \ref{whitney homotopy}.  By Proposition \ref{germs}  $\phi$ is $S$-equivalent to $\phi(\mG, \mR,\mF, \mW)$ where the boundary germs of $\mF$ and $\mW$ coincide.  After lifting to finite cover we can assume that the conclusion of Lemma \ref{high cover} holds.  Now let $\mF'$ be as in Lemma \ref{whitney homology}.  By factorization $\phi$ is $S$-equivalent to $\phi(\mG, \mR, \mF', W)\circ \phi(\mG, \mR, \mF, \mF')$ and hence $\phi$ is $S$-equivalent to $\phi(\mG, \mR, \mF', \mW)$.

Since the inclusion $Y\to V_k$ is a $\pi_1$-isomorphism, it follows that $\tilde Y=\tilde V_k\setminus N(\tilde G\cup \tilde R)$.  Since $\diam(\pi(f'_p\cup w_p))\cap Z\le k/4$ it follows that $f'_p$ is homologous to $w_p$ by a chain disjoint from  $\piinv(j+1/2)$ for some  $j\in\{1, 2 ,\cdots k\}$.  It follows if $\tilde f'_p, \tilde w_p$ are lifts with common boundary germs, then they are  homologous in $H_2(\tilde Y, \partial \tilde\mE)$ and hence are homotopic rel partial in $\tilde Y$.  Therefore, for all $p$, $f'_p$ is homotopic rel $\partial$ to $w_p$ in $Y$. \qed


\section {Homotopy implies Concordance}

In \cite{MM} Maggie Miller shows that in a 4-manifold whose fundamental group contains no 2-torsion,  homotopic 2-spheres are concordant provided one of them has a dual 2-sphere and under suitable hypothesis this holds more generally.  Here we use \S 4.1 \cite{MM} to show that homotopic $F|W$ systems are S-equivalent to ones with $\mW$ \emph{standardly concordant} to $\mF$.

\begin{figure}
 \centering
\begin{tabular}{ c c }
 $\includegraphics[width=5.5in]{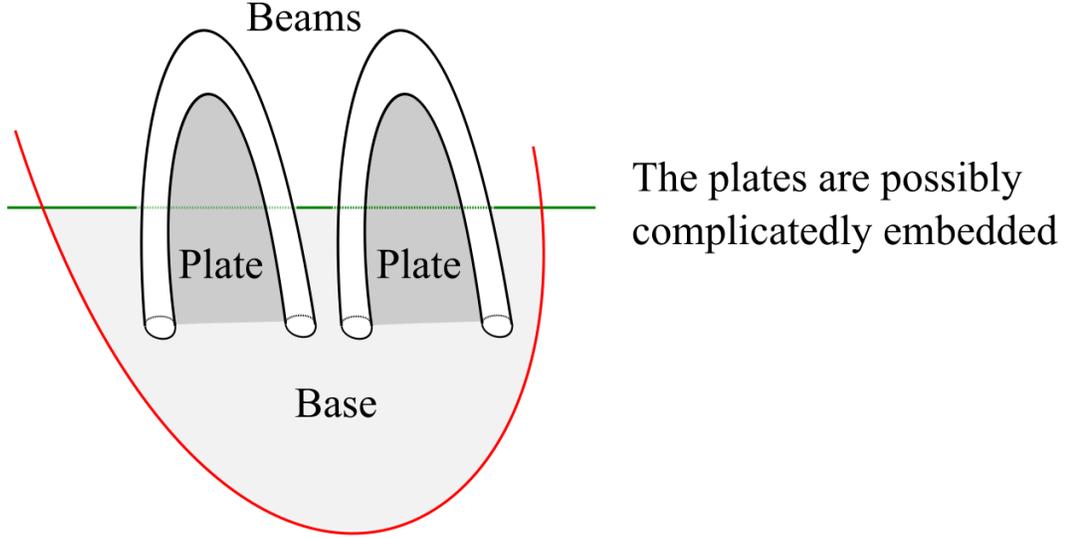}$  
\end{tabular}
 \caption[(a) X; (b) Y]{\label{base beam plate}
 \begin{tabular}[t]{ @{} r @{\ } l @{}}
A standardly concordant Whitney disc
\end{tabular}}
\end{figure}

\begin{definition}  The $\partial$-germ coinciding $F|W$ system $(\mG,\mR,\mF,\mW)$ has $\mW$  \emph{standardly concordant} to $\mF$ if it is obtained by starting with $\mF\times [0,1]$, attaching standardly embedded cancelling 3-dimensional 1- and 2-handles to $\mF\times 1$ and then reimbedding the 2-handles to obtain an immersed $\mF\times [0,1]$ whose $\mF\times 1:=\mW_1$ is embedded. The reimbedded 2-handles are required to coincide with the original ones near their attaching regions. Finally $\mW$ is obtained from $\mW_1$ by isotoping slightly to regain the $\partial$-germ coinciding condition. This essentially describes a standard concordance as a critical level embedding    

The data for a standard concordance is given by  \emph{bases, beams and plates}.  See Figure \ref{base beam plate}.  The bases are $\mF$. The beams are the 3-dimensional 1-handles slightly extended to attach to the bases and the plates are the cores of the reimbedded 2-handles, slightly extended to attach to the beams and bases.  While the plates are embedded, their interiors intersect the beams in parallel copies of the beams' cores, disjoint from beams' lateral surfaces and are otherwise disjoint from the bases.   \end{definition}

\setlength{\tabcolsep}{60pt}

\begin{proposition} \label{concordance} Let $(\mG,\mR,\mF,\mW)$ by a $F|W$ system on $V_k$ with $\mW$ homotopic to $\mF$.  Then there exists a $F|W$ system $(\mG,\mR,\mF,\mW')$ with $\mW'$ standardly concordant to $\mF$ and $\phi(\mG,\mR,\mF,\mW)$ S-equivalent to $\phi(\mG,\mR,\mF,\mW')$.\end{proposition}

\begin{proof}  Let $\mN_f$ (resp. $\mN_w$) be a system of dual spheres to $\mF$ (resp. $\mW$), which exists by Lemma \ref{whitney duals}.  If $|\mF|=n$, then let $f_t:\cup_{i=1}^n D^2\to V_k$ be a regular homotopy \cite{Sm} supported away from $N(\mG\cup\mF)$ such that $f_0=\mF$ and $f_1=\mW$, where we abuse notation by identifying a map with its image.  We will assume that all the finger moves occur at $t=1/4$ and all the Whitney moves when $t=3/4$.  Since the support of the finger moves can be chosen to lie in a small neighborhood of the union of $\mF$ and \emph{finger} arcs we can assume that $f_{.251}$ is \emph{dual} to $\mN_f$.  This means that $f_{.251}(D_i)$ intersects $\mN_f$ exactly once and at the component  of $\mN_f$ that intersects $f_0(D_i)$. Since $f_1$ is obtained by first applying Whitney moves to $f_{.749}$ and then isotopy, it follows that $f_{.749}$ is obtained from $f_1$ by finger moves and isotopy.  Thus we can assume that $f_{.749}$ is dual to $\mN_w$ and that there is an ambient isotopy from $f_{.251}$ to $f_{.749}$ which starts out dual to $\mN_f$ and ends dual to $\mN_w$.  

Let $X_k$ denote $V_k\setminus N(\mR\cup \mG\cup \mF)$.  Since the map  from  $\pi_1(X_k)$ to $\pi_1(V_k)$ induced by inclusion is an isomorphism we can assume that after passing to finite cover of $V_k$ that each finger arc represents the trivial element of $\pi_1(X_k)$ and hence corresponds to a standard arc from $\mF$ to itself.  Note that an arc may go from one component of $\mF$ to another.

Each finger move $f_k$ creates two points of self intersection of $f_{.251}$ which we eliminate by replacing two discs on one component of $f_{.251}$ by a tube $T_k$ as in Figure 6 \cite{MM} that lies very close to an arc $\alpha_k\subset f_{.251}$.  Let $R_1$ denote the embedded surface obtained by modifying $f_{.251}$ by these tubings.  The ambient isotopy induces a diffeomorphism $g_1:V_k\to V_k$ with $g_1(f_{.251})=f_{.749}$ and an isotopy of $R_1$ to $g_1(R_1):=R_2$, where we assume that $g_1(T_k)$ is a tube lying close to $g_1(\alpha_k)$, i.e. we can assume that $g_1$ takes a small regular neighborhood of $\alpha_k$ containing $T_k$ to a small regular neighborhood of $g_1(\alpha_k)$.  The core circle of $T_k$ bounds a small 2-disc $A_k$ that intersects $f_{.251}$ in a point $a_k\in \alpha_k$.  Let $g_1(A_k):=E_k$.  Since $f_{.251}$ (resp. $f_{.749}$) is dual to $\mN_f$ (resp. $\mN_w$) we obtain an embedded disc $B_k$ (resp. $F_k$) by tubing $A_k$ (resp. $E_k$) with a copy of a component of $\mN_f$ (resp. $\mN_w$).  The tube follows an embedded arc in $f_{.251}$ (resp. $f_{.749})$ from $a_k$ (resp. $g_1(a_k)$) to a point near $\mN_f$ (resp. $\mN_w$).  
Thicken $F_k$ to a 3-dimensional 2-handle $\omega_k$ and let $\mW'$ be the result of embedded surgery of $R_2$ along the $\omega_k$'s.  By construction, $\mW'$ $\partial$-coincides with $\mW$ and is dual to $\mN_w$.  By Lemma \ref{common duals} $\phi(\mG,\mR,\mF,\mW)$ is S-equivalent to $\phi(\mG,\mR,\mF,\mW')$.

To complete the proof we show that up to isotopy $\mW'$ is standardly concordant to $\mF$.  After an isotopy we using $\mN_w$, can assume that the 1-handles are short and standard.  As in \S 4.2 \cite{MM} up to isotopy $R_1$ is obtained from $f_{.249}$ by embedded surgery along  3-dimensional 1-handles, one 1-handle for each finger move.  See Figure 6 \cite{MM}.  Up to isotopy the thickened $B_k$'s denoted $\tau_k$'s are cancelling 2-handles.  Now reimbed each $\tau_k$ using $\ginv_1(\omega_k)$.  Since $\mW'$ is isotopic to the surface obtained by embedded surgery to $R_1$ along the $\ginv_1(\omega_k)$'s the result follows.  \end{proof}

\begin{remark} The dual spheres $\mN_f$  and the $B_k$'s were mentioned  to make the argument more symmetric and transparent, however as in \cite{MM} they are not needed in the proof.\end{remark}

\section{Finger$|$Whitney-carving/surgery presentations of Schoenflies spheres}

We start by defining the notion of a carving/surgery presentation of a Schoenflies sphere.  We then define a specialized form of this presentation called a $F|W$-carving/surgery presentation.  The main result of the next section is that every Schoenflies sphere has a $F|W$-carving/surgery presentation.  There we assume the conclusion of Proposition \ref{germs}.  The subsequent section shows that a $F|W$ system satisfying Proposition \ref{concordance} has an \emph{optimized} $F|W$-carving/surgery presentation.

\begin{definition} Let $\BS_0\subset S^4$ denote the standard 3-sphere and let $X^0_S, X^0_N$ denote its 4-ball closed complementary regions, respectively called the \emph{southern} and \emph{northern 4-balls}.    A \emph{carving/surgery presentation} $(\mL,\mD)$ of a Schoenflies sphere $\BS$ in the 4-sphere with closed complementary regions $\Delta_S$ and $\Delta_N$ consists of 

i) a framed link $\mL=\{k_1, \cdots, k_n\}\subset \BS_0$ that surgers $\BS_0$ to $S^3$ such that each component $k_i$ bounds an embedded disc $D_{k_i}$ such that $D_{k_i}\cap \BS_0\subset N^3(\mL)$ and induces the given framing on $k_i$.  Here $N^3(\mL)$ denotes a regular neighborhood of $\mL\subset \BS_0$.  If a neighborhood of $\partial D_{k_i} \subset X^0_S$ (resp. $X^0_N$), then label $k_i$ with a $\bullet$ (resp. 0).  To minimize notation we will usually denote $D_{k_i}$ by $D_i$.  Let $\mD=\{D_1,\cdots, D_n\}$.

ii) If $\inte(D_{i})\cap N^3(k_j)\neq\emptyset$, then the components of  $D_{i}\cap N^3(k_j)$ bound pairwise disjoint discs in $D_{i}$  that are parallel to $D_{j}$ and contained in a small neighborhood of $D_{j}$.  This induces a partial order $>$ on the $k_i$'s and the $D_{i}$'s, with the minimal elements corresponding to discs $D_k$ such that $D_k\cap \BS_0=k$.  We depict partial orders as forests with oriented edges and the simplicial metric.  Call a $k\in L$ or its corresponding disc $D_k$ a \emph{depth-m} element if the minimal distance to a leaf is $m-1$, so minimal elements are depth-1.  Let $M$ be the maximal depth over all the components of $L$.  

The CS-presentation $(\mL, \mD)$ gives rise to Schoenflies balls $\Delta_S, \Delta_N$ as follows.  

i) Let $X_S^1\subset S^4$ denote the compact manifold obtained from $X^0_S$ by carvings \cite{Ak} corresponding to the minimal discs contained in $X^0_S$ and by adding 2-handles corresponding to the minimal discs contained in $X^0_N$.  Let $X_N^1$ denote the closed complementary region of $X_S^1$.  Let $X_S^2\subset S^4$ denote the compact manifold obtained from $X_S^1$ by carving the depth-2 discs lying in $X_S^1$ and adding 2-handles corresponding to the depth-2 discs lying in $X_N^1$.  Let $X_N^2$ denote the closed complementary region to $X_S^2$.  In a similar manner construct $X_S^3, \cdots, X_S^M$ and $X_N^3, \cdots, X_N^M$.

ii)  Define $\Delta_S=X_S^M$, $\Delta_N=X_N^M$ and $\BS=\partial \Delta_S=\partial \Delta_N$.  By construction $\BS=S^3$ and $\Delta_S$ and $\Delta_N$ are its Schoenflies balls.\end{definition}

\begin{definition} \label{nesting2}  If $k$ is depth $m$ and labeled with a 0, then let $D^0_k$ denote the associated 2-handle attached to $X^{m-1}_S$, while if labeled with a $\bullet$, then let $D^\bullet_k$ denote the associated 2-handle  attached to $X^{m-1}_N$.  We will abuse notation by calling $D^0_k$ (resp. $D^\bullet_k$) a \emph{2-handle of $\Delta_S$} (resp. $\Delta_N$) and a \emph{carving of $\Delta_N$} (resp. $\Delta_S$).  Viewed in $\Delta_S$, $D^0_k$ will be a 2-handle with a $D^2\times U$ deleted, where $U$ is a possibly empty open subset of the cocore.  E.g. if $k$ is depth-1 and $M=2$ then $D^0_k$ may have finitely many parallel copies of its core carved out.  \end{definition}

\begin{notation} \label{reindex} Recall from Notation 3.7 that $V_k=\sonesthree\#_k\stwostwo$ with $V_\infty$ its universal cover where the lifts of the $\stwostwo$ summands are denoted $\widetilde\stwostwo_i, i\in\BZ$.  In what follows $\tilde \BS_0\subset V_\infty$ will denote a 3-sphere of the form $pt\times S^3$ that separates the $\widetilde\stwostwo_i$'s where $i>0$ from the $\widetilde\stwostwo_j$'s where $j\le 0$.  To make the notation more symmetric about 0, reindex $\widetilde\stwostwo_i$ by $\widetilde\stwostwo_{i-1}$ when $i\le 0$.  Similarly reindex $\tilde R_i, \tilde G_i$, etc when $i\le 0$.  A finger (resp. Whitney) disc from $\tilde R_i$ to $\tilde G^{\rm{std}}_j$ will be denoted $\tilde f_{ij}$ (resp. $\tilde w_{ij}$), though there may be more than one such disc for a given $i$ and $j$. \end{notation}

\begin{definition}\label{nesting lodging} We say that $D_i$ \emph{nests in $D_j$} if some component of $D_i\cap N(k_j)$ is outermost in $D_i\cap (\cup_{p\neq i} N(k_p))$.  We similarly say that $k_i$ or $D_i$'s 2-handle or carving \emph{nests in} $k_j$ or $D_j$'s carving or 2-handle.  In this case we also say that  $D_j $ or its 2-handle or carving or $k_j$ \emph{lodges} $D_i$ or its carving or 2-handle or $k_i$.\end{definition}

\begin{definition}  Given the pairwise disjoint knots $k_1, \cdots, k_r $ in $\BS_0$, then $k'_1, \cdots, k'_r$ are called \emph{linking circles} if they bound pairwise disjoint discs in $\BS_0$ that $\delta_{ij}$ intersect the $k_i$'s.    Once defined, usually implicitly, they and their corresponding annuli $\subset \BS_0\setminus \cup \inte(N_{k_i}\cup N_{k_i'})$ are fixed once and for all, with subsequent geometric operations and constructions avoiding these annuli.\end{definition}

\begin{definition}  \label{fwcs} A \emph{$FW$-carving/surgery presentation} (\emph{FWCS-presentation}) is a CS-presentation whose link $\mL=B_L\sqcup S_L\sqcup N_L$ is a disjoint union of \emph{knots} $L_k$ and their linking circles $L'$ such that all the linking circles are 0-framed.  Furthermore $L_k=B_L^k\sqcup S_L^k\sqcup N_L^k$ and $L'=B'\cup S'\cup N'$ with the partial order and 0, $\bullet$ labeling arising as follows.  
\vskip 8pt

1a) $B_L=B_L^k\cup B'$ where $B_L^k$ is the unlink $U=\{b_{i_1}, \cdots, b_{i_p}\}$ and $B'=\{b'_{i_1}, \cdots, b'_{i_p}\}$ are its linking circles.  Here $i_j\in \BZ_{\neq 0}$ and $i_j\neq i_k$ unless $j=k$.  If $i<0$ (resp. $>0$), then $b_i$ is labeled with a 0 (resp. $\bullet$) and it's linking circle $b'_i$ is labeled with a $\bullet$ (resp. 0).  The bulleted (resp. 0'd) discs bound the disc unlink in $X^0_S$ (resp. $X^0_N$).   Within $B_L$, there  is an arrow from each bulleted linking circle to each knot with label 0 and an arrow from each linking circle with label 0 to each bulleted knot.  \vskip 8pt

\emph{Motivation} (where this comes from):  A $b_i, i<0$ (resp. $i>0$) will arise if there is a finger disc of the form $\tilde f_{ji}$ with $j>0$ (resp. $j<0$).
\vskip 8pt

1b) $S_L=\sqcup S_{ij}$ where $i<0$ and $j\in \BZ_{\neq 0}$.  Associated to each $S_{ij}$ is a finite union $\mW_{ij}$  of 2-discs each of which contains finitely many pairwise disjoint simple closed curves whose union is denoted $S^k_{ij}$.  The elements of $S^k_{ij}$ are in 1-1 correspondence with  a set of knots, also denoted $S^k_{ij}\subset\BS_0$ and each knot comes with a linking circle.  $S_{ij}$ is the union of these knots and their linking circles  $S'_{ij}$.   \vskip 8pt

\emph{Motivation} (where this comes from):  The discs associated to $S_{ij}$ are the $\tilde w_{ij}$ discs with $i<0$ and the simple closed curves are intersections of the $\tilde w_{ij}$'s with $\tilde \BS_0$.
\vskip 8pt

1b) continued:  The partial ordering restricted to the knots is induced by the inclusion relation of the simple closed curves within the discs, with the innermost curves being the minimal ones.  The linking circles are given the opposite partial ordering.  If two knots are connected by an arrow, then one is labeled with a 0 and the other with a $\bullet$, subject to the condition that the maximal knots are labeled with 0's. If a knot is labeled with a 0 (resp. $\bullet$), then its linking circle is labeled with a $\bullet$ (resp. 0).    See Figure \ref{ij,Family}.
The boxes on the knots remind us that they may be knotted and linked with other knots of $L$.
\setlength{\tabcolsep}{60pt}
\begin{figure}
 \centering
\begin{tabular}{ c c }
 $\includegraphics[width=5.5in]{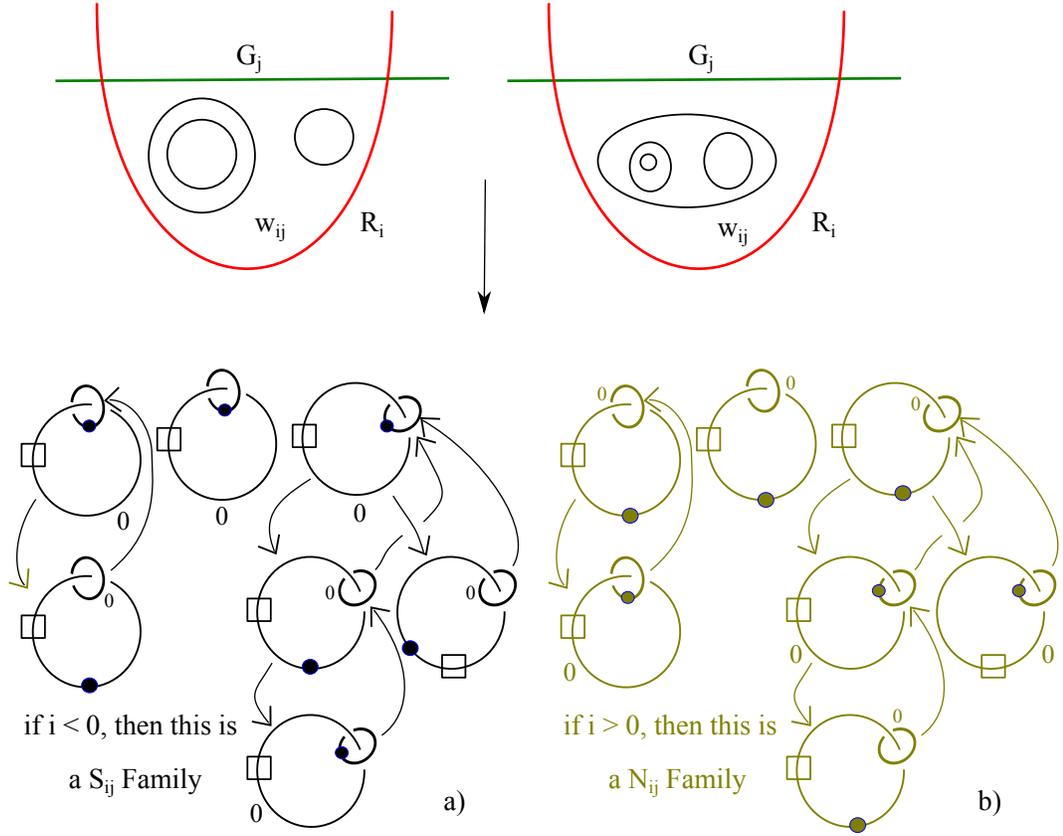}$  
\end{tabular}
 \caption[(a) X; (b) Y]{\label{ij,Family}
 \begin{tabular}[t]{ @{} r @{\ } l @{}}
Constructing the $S_{ij}$ and $N_{ij}$ Families from Intersection Data
\end{tabular}}
\end{figure}
\vskip 8pt
1c) $N_L=\sqcup N_{ij}$ where $i>0$ and $j\in \BZ_{\neq 0}$.  Associated to each $N_{ij}$ is a finite union $\mW_{ij}$ of discs each of which contains pairwise disjoint simple closed curves whose union is denoted $N^k_{ij}$.  The elements of  $N^k_{ij}$ are in 1-1 correspondence with a set of knots, also denoted $N^k_{ij}\subset \BS_0$ and each knot comes with a linking circle.  $N_{ij}$ is the union of these knots and their linking circles $N'_{ij}$.  \vskip 8pt
\emph{Motivation} (where this comes from):  The discs associated to $N_{ij}$ are the $\tilde w_{ij}$ discs with $i>0$ and the simple closed curves are intersections of the $\tilde w_{ij}$'s with $\tilde \BS_0$.
\vskip 8pt
1c) continued:  The partial ordering restricted to the knots is induced by the inclusion relation of the simple closed curves within the discs, with the innermost curves being the minimal ones.  The linking circles are given the opposite partial ordering.  If two knots are connected by an arrow, then one is labeled with a 0 and the other with a $\bullet$, subject to the condition that the maximal knots are labeled with bullets. If a knot is labeled with a 0 (resp. $\bullet$), then its linking circle is labeled with a $\bullet$ (resp. 0).   See Figure \ref{ij,Family}.
\vskip 8pt
1d) Let $\alpha$ be a knot of $S_{ij}$ or $N_{ij}$.  We say that  $\alpha$ is at  \emph{level n} if there exists an arc from $\alpha$ to $\partial \mW_{ij}$ which intersects $N^k_{ij}\cup S^k_{ij}$ $n$ times where $n$ is the minimal possible, so maximal knots are exactly those of level-1.
\vskip 8pt
 2a) \emph{Order relations involving $B_L$ and $N_L\cup S_L$}:  For every $i,j,k\in \BZ_{\neq 0}$ with $j<0$ and $k>0$ construct directed edges according to Figure \ref{BL,edges} a).  I.e. for each maximal knot of $N_{kj}$ there is an arrow from its linking circle to the linking circle $b'_j$ of $b_j$.  For every maximal knot of $S_{ji}$ there is an arrow from $b'_j $ to the knot.

 In addition, for every knot of $N_L\cup S_L$ labeled with a $\bullet$, maximal or not, construct an arrow to each $b_j$ where $j<0$.  

For every $i,j,k\in \BZ_{\neq 0}$ with $j>0$ and $i<0$ construct directed edges according to Figure \ref{BL,edges} b).  I.e. For each maximal knot of $S_{ij}$ there is an arrow from its linking circle to $b'_j$.   For each maximal knot of $N_{jk}$ there is an arrow from $b'_j$ to the knot.  

In addition, for every knot of $N_L\cup S_L$ labeled with a 0, maximal or not, construct an arrow to each $b_j$ where $j>0$. 

Finally, for every linking circle $k'$ of a maximal knot $k$ of $S_L$ (resp. $N_L$) construct an arrow to each $b_j$, where $j<0$ (resp. $j>0$).

\setlength{\tabcolsep}{60pt}
\begin{figure}
 \centering
\begin{tabular}{ c c }
 $\includegraphics[width=2.5in]{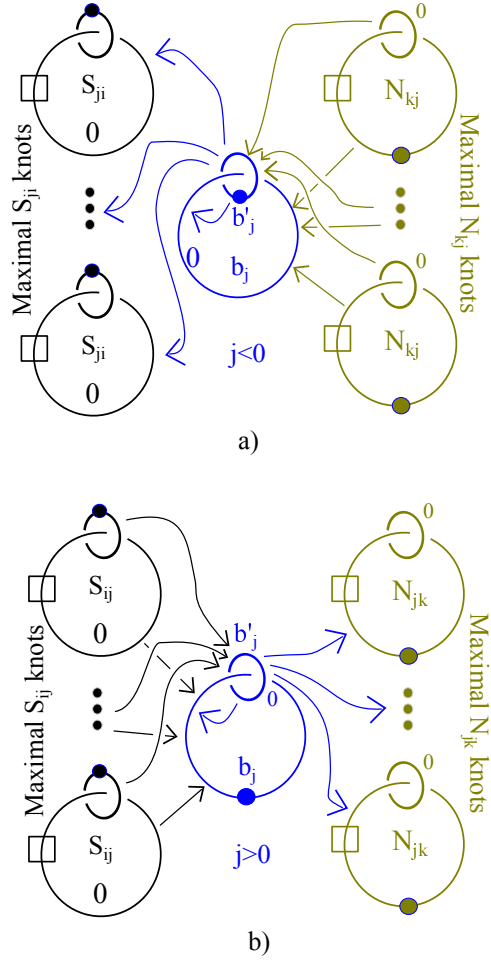}$  
\end{tabular}
 \caption[(a) X; (b) Y]{\label{BL,edges}
 \begin{tabular}[t]{ @{} r @{\ } l @{}}
Order relations involving $B_L$.  Non maximal knots not shown.
\end{tabular}}
\end{figure}
\vskip 8pt
2b) \emph{Additional order relations involving $S_L$}:  For every maximal knot $k$ of $S_{ji}$ and every maximal knot $k_1$ of $S_{ik}$ construct an arrow from the linking circle $k'$ of $k$ to the knot $k_1$.  In particular for every maximal knot $k$ of $S_{ii}$ with linking circle $k'$ there is an arrow from $k'$ to $ k$ and if $k_1$ is another maximal knot of $S_{ii}$ with linking circle $k_1'$, there is an arrow from $k_1'$ to $k$ and an arrow from $k'$ to $k_1$.   
\vskip 8pt
2c) \emph{Additional order relations involving $N_L$}:  For every maximal knot $k$ of $N_{ji}$ and every maximal knot $k_1$ of $N_{ik}$ construct an arrow from the linking circle $k'$ of $k$ to the knot $k_1$.  In particular for every maximal knot $k$ of $N_{ii}$ there is an arrow from $k'$ to $k$ and if $k_1$ is another maximal knot of $N_{ii}$ with linking circle $k_1'$, there is an arrow from $k_1'$ to $k$ and an arrow from $k'$ to $k_1$.   \end{definition}

This completes the description of the additional combinatorial structure needed for a CS-presentation to be a FWCS-presentation, subject to checking that the order relation is a partial order.

\begin{lemma} The order relation on the link $L$ of a FWCS-presentation is a partial order.\end{lemma}

\begin{proof}  The knots of $B_L$ are minimal elements.  After deleting these elements, all remaining directed edges from the knots of $S_L$ (resp. $N_L$) go to knots of $S_L$ (resp. $N_L$) and the ordering of these knots is induced from the partial ordering on embedded simple closed curves on discs.    Thus it suffices to prove the lemma for the ordering restricted to $L'$.  Here the elements of $B'$ are minimal elements and again the ordering of what remains is induced from the partial ordering on embedded simple closed curves on discs.\end{proof}

\begin{figure}
 \centering
\begin{tabular}{ c c }
 $\includegraphics[width=2.5in]{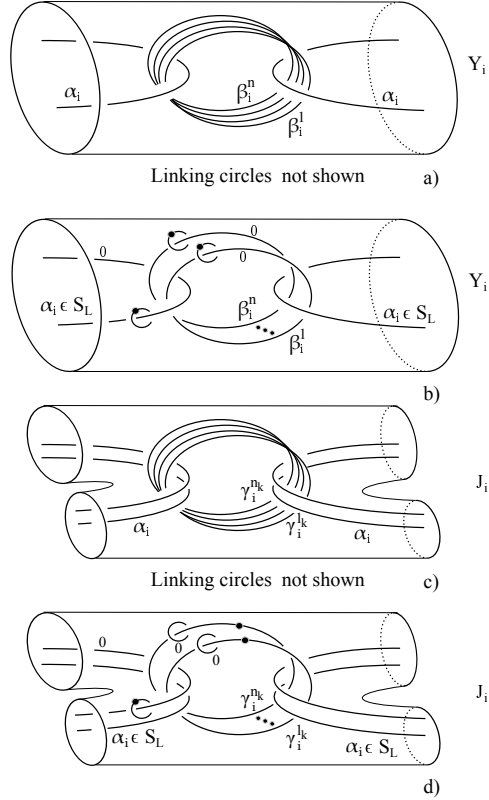}$  
\end{tabular}
 \caption[(a) X; (b) Y]{\label{optimized1}
 \begin{tabular}[t]{ @{} r @{\ } l @{}}
Local view of an optimized FWCS-Presentation near $S_L\cup N_L$.  
\end{tabular}}
\end{figure}

\begin{definition}  We continue to call a CS-presentation an FWCS-presentation if it satisfies the conditions of Definition \ref{fwcs} except that the arrows are a proper subset of those stated in that definition.\end{definition} 
\begin{definition}  \label{optimized} An \emph{optimized} FWCS-presentation is a FWCS-presentation with the following additional features.
\vskip 8pt
i) The knots of $S_{ij}$ and $N_{ij}$ are of level at most 2 (and so the level of a knot is determined by whether or not it's $\in S^k_L$ and whether or not it's labeled with a $0$). 

ii) If $k \in L_k$, then $D_k$ induces the 0-framing on $k$.  

iii) The knots of $S_L\cup N_L$ are the disjoint union of $\mA$, $\mB$ and $\mC$ where $\mB\cup \mC$ includes all the level-2 knots and if $k\in\mB\cup \mC$, then $k$ is a minimal element with respect to the partial order on $L$, i.e. no arrows point out of $k$.  

iv) If $\mA=\{\alpha_1, \cdots, \alpha_m\}$, then there are pairwise disjoint solid tori $Y_1, \cdots, Y_m\subset \BS_0$ such that $Y_i\cap (B^k_L\cup \mC)=\emptyset, Y_i\cap\mA=\alpha_i$ and $\mB\subset \cup Y_i$.   Let $Y_i\cap \mB=\{\beta_i^{i_1}, \cdots, \beta_i^{i_m}\}$. Then $\alpha_i, \beta_i^{i_1}, \cdots, \beta_i^{i_m}$ appear as in Figure \ref{optimized1} a), i.e. the $\beta_i^j$'s are parallel unknots and $\alpha_i$ is obtained from two parallel copies of the core of $Y_i$ banded together by the trivial band that goes through the $\beta_i^j$'s.    If $\alpha_i\in N^k_L$, then $\mB\cap Y_i\cap N^k_L $ (resp. $\mB\cap Y_i\cap S^k_L$) are level 1 (resp. level 2) curves.  If $\alpha_i\in S_L^k$, then $\mB\cap Y_i\cap S^k_L$ (resp. $\mB\cap Y_i\cap N^k_L$) are level 1 (resp. level 2) curves.  Figure \ref{optimized1} b) shows the case where $\alpha_i\subset S^k_L$.

v) There are pairwise disjoint handlebodies $J_1, \cdots, J_m\subset \BS_0$ such that $(J_1\cup\cdots\cup J_m)\cap B^k_L=\emptyset$ and $\mC\subset J_1\cup\cdots\cup J_m$.  Here  $J_i$ is obtained by adding 1-handles to a slightly thickened $Y_i$  and near each 1-handle, $\mC$ intersects $J_i$ in parallel unknots that clasp $\alpha_i$ as in Figure \ref{optimized1} c).    If $\alpha_i\in N^k_L$, then $\mC\cap J_i\cap N^k_L $ (resp. $\mC\cap J_i\cap S^k_L$) are level 2 (resp. level 1) curves.  If $\alpha_i\subset S^k_L$, then $\mC\cap J_i\cap S^k_L$ (resp. $\mC\cap J_i\cap N^k_L$) are level 2 (resp. level 1) curves.  Figure \ref{optimized1} d) shows the case where $\alpha_i\subset S^k_L$.  
 \end{definition}
 
 \setlength{\tabcolsep}{60pt}

\begin{remarks}  i) $\mA\cup\mC\cup B^k_L$ is the unlink as is $\mB\cup \mC\cup B^k_L$.

ii) For a given $\alpha_i$, there is exactly one set of parallel $\beta_i^j$'s that link it, however there are finitely many sets of curves in $\mC$ that link it, so that $\genus(J_i)\in \BZ_{\ge 0}$.

iii)  The handlebodies $J_i$ may be knotted and may link each other and $B^k_L$.  \end{remarks}

\section{Schoenflies spheres have Finger$|$Whitney-carving/surgery presentations}

The following is the main result of this section.

\begin{theorem}\label{carving theorem} Every smooth 3-sphere in the 4-sphere has a FWCS-presentation.\end{theorem}

\begin{proof}  By Proposition \ref{equivalence} every 3-sphere $\Sigma'\subset S^4$ corresponds to the S-equivalence class of some $\phi\in \Diff_0(\sonesthree)$. Viewing $S^4$ as $(\BR\times S^3)\cup\{S,N\}$, then $\Sigma'$ is isotopic to $\tilde\phi(\pt\times S^3)$ where $\tilde\phi$ is the lift to $\widetilde \sonesthree$.  By Lashoff - Shaneson \cite{LS} and Sato \cite{Sa} there is a pseudo-isotopy $f$ from $\id$ to $\phi$ which by \cite{HW} arises from a 1-parameter Hatcher-Wagoner family $(q_t, v_t)$.  We will assume that $(q_t, v_t)$ induces a $F|W$ system satisfying the conclusion of Proposition \ref{germs} and that $(q_t,v_t)$ has been normalized as in the second and third paragraphs of the proof of Theorem \ref{stable isotopy}.  In what follows $V$ will denote $S^1\times S^3$ and $V_k$ will denote $\sonesthree\#_k\stwostwo$ where $k$ is the number of components of the nested eye.  It suffices to prove the theorem for the 3-sphere $\Sigma$ obtained from $\Sigma'$ by reversing orientation, which by Proposition \ref{equivalence} is the class of $\tilde\phi^{-1}(\pt\times S^3)$.  

\vskip 10 pt
\noindent\textbf{Step 1:} Show how to construct $\phiinv(U)$ from the $F|W$-system where $U$ is a closed submanifold of $V$.

\vskip 8pt

We view the 1-parameter family $(q_t, v_t)$ as a smoothly varying family  of handle structures $h_t$ on $V\times I$.  See the second paragraph of P. 174 \cite{HW} and \S4 \cite{Qu}.  Here  $h_{1/4}$ corresponds to $k$ 2 and 3-handles in standard cancelling position and for $t\in [1/4,3/4]$ the handle structure changes according to the path of 2-sphere boundaries of the cores of the attaching 3-handles.  

We introduce some terminology  to keep track of the data, in particular both before and after the 2-handle attachments.  First $h_t$ denotes the handle structure on $V\times I\times t$.  For $t\in [1/4,3/4]$ the 2-handles are attached to $V\times[0,1/4]\times t$ along a set of $k$ 0-framed simple closed curves $\Omega\times 1/4\times t$, where $\Omega=\{\omega_1, \cdots, \omega_k\}$ and the $\omega_i$'s bound pairwise disjoint discs $\mD=\{D_1, \cdots, D_k\}$.  We abuse notation by also viewing $\Omega$ and the $D_i$'s as subsets of $V\times0\times t$.  To simplify notation we let $V_k\times t$ denote the result of attaching these 2-handles.  This $V_k\times t $ corresponds to the $V_k\times 1/4\times t$ in the proof of Theorem \ref{stable isotopy}. With terminology as in \S2, $\mGs=\{\Gs_1, \cdots, \Gs_k\}\subset V_k\times t$ denotes the boundaries of the 2-handle cocores, where $\Gs_i$ is the standard green sphere in the $i$'th $\stwostwo$ factor, i.e. it is of the form $S^2\times \pt$. When $t=1/4$, the 3-handles are attached along the set of  spheres $\mRs=\{\Rs_1, \cdots, \Rs_k\}$ where $\Rs_i$ is the $i$'th standard red sphere, i.e. is  of the form $\pt\times S^2$ in the $i$'th $\stwostwo$ factor and $\Rs_i$ flows to $D_i\subset V\times 0\times {1/4}$ under $-v_{1/4}$.  Fix regular neighborhoods $N(\Omega)$ and $N(\mGs)$ so that $V_k\times t\setminus \inte(N(\mGs))$ is diffeomorphic to $V\times 0\times t\setminus \inte(N(\Omega))$ with the diffeomorphism $\lambda$ induced by $-v_t$.  Note that both are diffeomorphic to $\sonesthree\#_kS^2\times D^2$.  Thanks to our initial normalization, much of the data corresponding to the projections to $V_k$ and $V$ are independent of $t\in [1/4,3/4]$, e.g. $\Omega, \mGs, \mRs$, the $D_i$'s and $\lambda$.  To help keep track of $\lambda|N(\mGs)$ we recall the 
\vskip8pt

\noindent\emph{Section-Fiber Rule}:  For the $i$'th 2-handle $D^2\times D^3$ we have $\partial\textrm{ core}=\omega_i=\partial D^2\times 0$ and $\partial\textrm{ cocore}=\Gs_i=0\times \partial D^3$.  Here, with $\omega_i$ viewed in $V\times 0\times t$ and $\Gs_i$ viewed in $V_k\times t$, then $\partial D^2\times \partial D^3$ can be simultaneously viewed as $\partial N(\Gs_i)=N^1(\Gs_i)\subset V_k\times t$ and as $\partial N(\omega_i)=N^1(\omega_i)\subset V\times 0\times t$ with $\lambda$ equating the two.  Here $N^1$ denotes the unit normal bundle of the space in question. Under $\lambda$ a section $\pt\times \Gs_i$ of $N^1(\Gs_i)$ is identified to a normal fiber $\pt\times S^2$ of $N^1(\omega_i)$ and a normal fiber $S^1\times \pt $ of $N^1(\Gs_i)$ is identified with a section  $\omega_i\times\pt$ of $N^1(\omega_i)$.

\vskip8pt

Figure \ref{section fiber} a) shows a 3-dimensional slice of $N^1(\Gs_j)\cup \Rs_j\subset V_k\times t$.  The two points in a) are the intersection with a circle fiber and the dark 2-sphere on top and the light one below are 2-sphere sections of $N^1(\Gs_j)$.  On the other hand, we can view $N(\omega_j)\subset V\times 0\times t$ using cylindrical coordinates $(x,y,r,\theta)$ and Figure \ref{section fiber} b) shows the $\theta=0, \pi$ slices with two 2-sphere fibers and the intersection with a section.  

\vskip8pt
\setlength{\tabcolsep}{60pt}
\begin{figure}
 \centering
\begin{tabular}{ c c }
 $\includegraphics[width=4.5in]{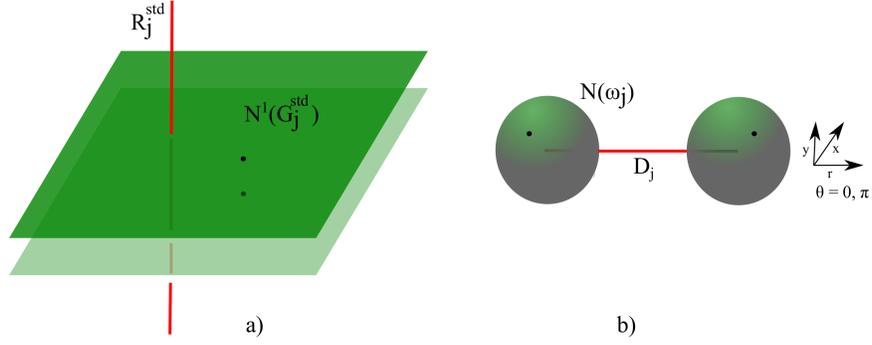}$  
\end{tabular}
 \caption[(a) X; (b) Y]{\label{section fiber}
 \begin{tabular}[t]{ @{} r @{\ } l @{}}
Sections and Fibers\end{tabular}}
\end{figure}
Let $\mR_t:=\{R_{1,t}, \cdots, R_{k,t}\}\subset V_k\times t$
 denote the boundary of the cores of the 3-handles, where $\mR_{1/4}=\mRs$.  When projected to $V_k$, $\mR_t$ first undergoes finger moves with $\mGs$ and then Whitney moves according to its $F|W$ system.  Let $t_0\in [1/4,3/4)$ be just after the Whitney moves are completed.  Let $\mR^0_t=\{R_{1,t}, \cdots, R_{k,t}\}$ denote $\lambda(\mR_t\setminus \inte N(\mGs))$.  We can assume that $\mR_t\cap N(\mGs)$ is a union of fibers of $N(\mGs)$ except during the finger and Whitney moves.   We can assume that for $t\in [1/4,t_0], \mR^0_t$ coincides with $\mR^0_{1/4}$ near $\partial \mR^0_{1/4}$, where $\mR^0_{1/4}=\mD\setminus \inte N(\Omega)$.  Having completed the Whitney moves, $\mR^0_{t_0}$ is a union of discs.  In what follows we abuse notation by viewing $\mR_t\subset V_k$ and $\mR^0_t\subset V$, rather than $V_k\times t$ and $V\times 0\times t$.  Under isotopy extension the path $\mR_t$ extends to an ambient isotopy $\eta_t:V_k\to V_k, t\in [1/4, t_0]$.

 We now show how to construct $\phiinv(U)$, up to isotopy,  from the $F|W$-system where $U$ is a closed submanifold of $V$.  Our $\phiinv(U):=U^0_1$ is obtained by flowing  $U\subset V\times 1\times 1$ via $-v_1$ to $V\times 1\times 0$.   To determine $U^0_1$ from our pseudo-isotopy we start with $U\subset V\times 1\times 0$ and keep track of the result $U^0_t$ of flowing, via $-v_t$, $U$ down to $V\times 0\times t$.   At times when $U$ transversely intersects the ascending spheres of the 2-handles, discs (if $\Dim(U)=2$) or solid tori (if $\Dim(U)=3$) do not make it down to  $V\times 0\times t\setminus \inte(N(\Omega))$.  After the deaths, discs or solid tori are restored producing $U^0_1$, up to isotopy.  Here are more details in the case $\Dim(U)=3$.

Since $\Dim(U)=3$ we can assume that $U$ is disjoint from the union of ascending spheres of the 2 and 3-handles of $V\times I\times 1/4$ since the latter is a disjoint union of $k$ 2-discs $\subset V\times 1\times 1/4$.  Let $U_t$ denote the result of flowing $U$ down to $V_k\times t$ under $-v_t$.  By construction $U_{1/4}$ is disjoint from $N(\mGs\cup \mRs)$ and for $t\in [1/4, t_0], U_t=\eta_t(U_{1/4})$.  We can assume that $U_{t_0}$ is transverse to $\mGs$ and hence  $U_{t_0}\cap N(\mGs)$ is a union of solid tori, each of which is a union of normal discs.  By construction it is disjoint from $\eta_t(N(\mRs))$.  Define $U_t^0:=\lambda(U_t\setminus \inte(N(\mGs)))\subset V\times 0\times t\setminus \inte(N(\Omega))$.  

We now cancel the 2 and 3-handles of $h_{t_0}$ or equivalently introduce deaths to directly extend  the 1-parameter family $(q_t, v_t), t\in [0, t_0]$ to a Hatcher - Wagoner family $(q_t',v'_t), t\in [0,1]$, where $(q'_t, v'_t)=(q_t, v_t)$ for $t\in [0, t_0] $ and $(q'_t, v_t')$ has no excess 3/2 intersections in $[t_0,1]$ and $(q_1,v_1)$ is nonsingular.  The induced pseudo-isotopy is from $\id$ to a $\phi':V\to V$ which by Lemma \ref{chenciner} is isotopic to $\phi$.  By Lemma \ref{loop to diff} the isotopy class of $\phi'$ is independent of the isotopy extension $\eta_t$.

Under  cancellation of the 2 and 3-handles of $h_{t_0}$ as in \cite{Mi}, see P. 174 \cite{HW}, $U^0_{t_0}$ is modified as follows to obtain $U^0_1=\phi^{-1}(U)$, up to isotopy.  If $C$ is a normal fiber of $U_{t_0}\cap N^1(\mGs)$, then $\lambda(C)$ is a section of $N^1(\Omega)$ which is a circle that  is capped off by the union of an annulus in $N(\Omega)$ and a parallel copy of a component of $\mR^0_{t_0}$.  Note that if $K$ is any compact set $K'$ of fibers of $N^1(\mGs)$ disjoint from $\mRs$, e.g. $U_{t_0}\cap N^1(\mGs)$, then $\lambda(K)$ is capped off by a $K'$ family of annuli and a $K'$ family of parallel copies of components of $\mR^0_{t_0}$.  See Figure \ref{glissade figure}.  This completes the construction of Step 1.

\setlength{\tabcolsep}{60pt}
\begin{figure}
 \centering
\begin{tabular}{ c c }
 $\includegraphics[width=4.5in]{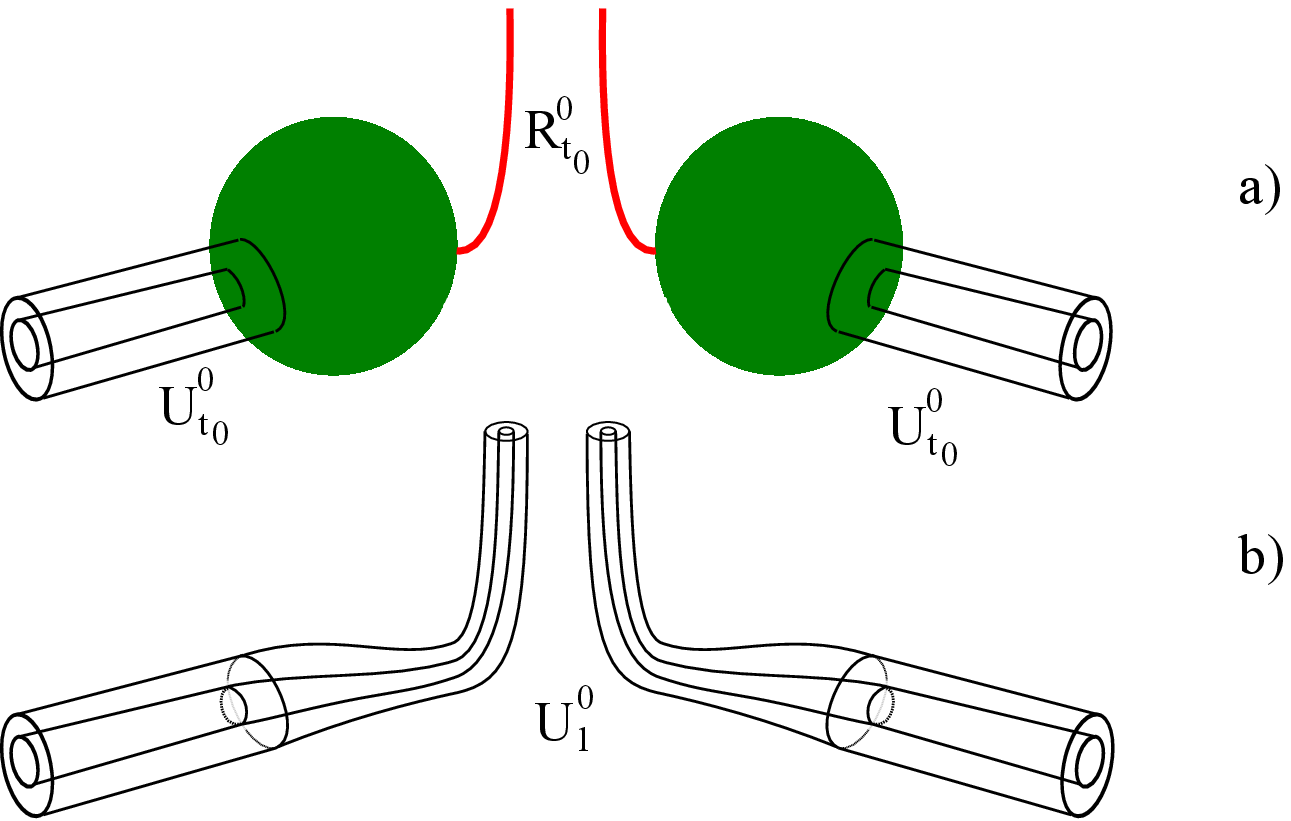}$  
\end{tabular}
 \caption[(a) X; (b) Y]{\label{glissade figure}
 \begin{tabular}[t]{ @{} r @{\ } l @{}}
A Glissade\end{tabular}}
\end{figure}

\begin{definition}  We call the operation of replacing $U^0_{t_0}$ by $U^0_1$ a \emph{glissade}.   \end{definition}

\begin{notation}  \label{rereindex} As in Notation \ref{reindex} we modify Notation \ref{winding} by reindexing $\widetilde\stwostwo_i, \widetilde\Rs_i, \widetilde \Gs_i$, etc. so that each $i\le 0$ is replaced by $i-1$.  Therefore $\tilde\pi^{-1}[1-2\epsilon, -1+2\epsilon]=[1-2\epsilon,-1+2\epsilon]\times S^3$ with $\rho$ modified accordingly.

Let $\tilde h_t$ denote the lift of the handle structure $h_t $ to $\tilde V\times I\times t$.  In a similar manner we reindex its data.  Let $\BS_0:=0\times S^3\subset \tilde V$ and $\BS_0':=0\times S^3\subset \tilde V_k$.  With $S^4$ identified with $\tilde V\cup\{S,N\}, \BS_0$ is our standard 3-sphere in $S^4$.\end{notation}

\noindent\textbf{Step 2:}  Apply the construction of Step 1 to $\BS_0$ to construct a sphere $\BS_1$ isotopic to $\phiinv(\BS_0)$. 
\vskip 8pt
Let $\BS'_{1/4}=\BS'_0\subset \tilde V_k$ and  $\BS_{1/4}=\BS_0 \subset \tilde V$.  In what follows we will assume that the finger (resp. Whitney) moves occur during  $t\in (3/8-\epsilon, 3/8)$ (resp. $t\in (5/8-\epsilon, 5/8)$) and $t_0=5/8$ and that the finger moves associated to a given arm are done simultaneously.  We will also assume that the handle cancellations occur during $t\in (11/16-\epsilon, 11/16)$ and so $\BS_{11/16}=\BS_1$.  We now construct an explicit isotopy $\BS'_t\subset \tilde V_k$ with $\widetilde \mR_t\cap \BS_t'=\emptyset$ for $t\in [1/4, t_0]$ and use it to construct $\BS_1$.  
\vskip 8pt
\noindent i) \emph{ Prepare for the finger moves}:  For each arm-hand of $\tilde\mF$ from $\widetilde\Rs_i$ to $\widetilde\Gs_j $ with $i<0$ and $ j>0$ isotope $\BS'_{1/4}$ to $\BS'_{5/16}$ by doing a finger move from $\BS'_{1/4}$ to $\widetilde\Gs_j$ so  that $\BS'_{5/16}$ is disjoint from the track of that arm-hand.  This isotopy creates a torus component of $\BS'_{5/16}\cap\partial N(\widetilde\mGs_j)$.  The corresponding passage from $\BS_{1/4}$ to $\BS_{5/16}:=\lambda(\BS'_{5/16}\setminus \inte (N(\widetilde\mGs)))$ involves the removal of an open solid torus.  In a similar manner isotope  $\BS'_{1/4}$ to $\BS'_{5/16}$ when $i>0$ and $j<0$.  The supports of all these isotopies need to be disjoint.   Figures \ref{finger move} a), b) show the before and after local pictures in $\tilde V_k$, i.e. at times $t=5/16-\epsilon$ and $t=5/16$ where $5/16-\epsilon$ is a time just before $\BS'_t$ locally intersects $\widetilde\mGs$.
\setlength{\tabcolsep}{60pt}
\begin{figure}
 \centering
\begin{tabular}{ c c }
 $\includegraphics[width=4.5in]{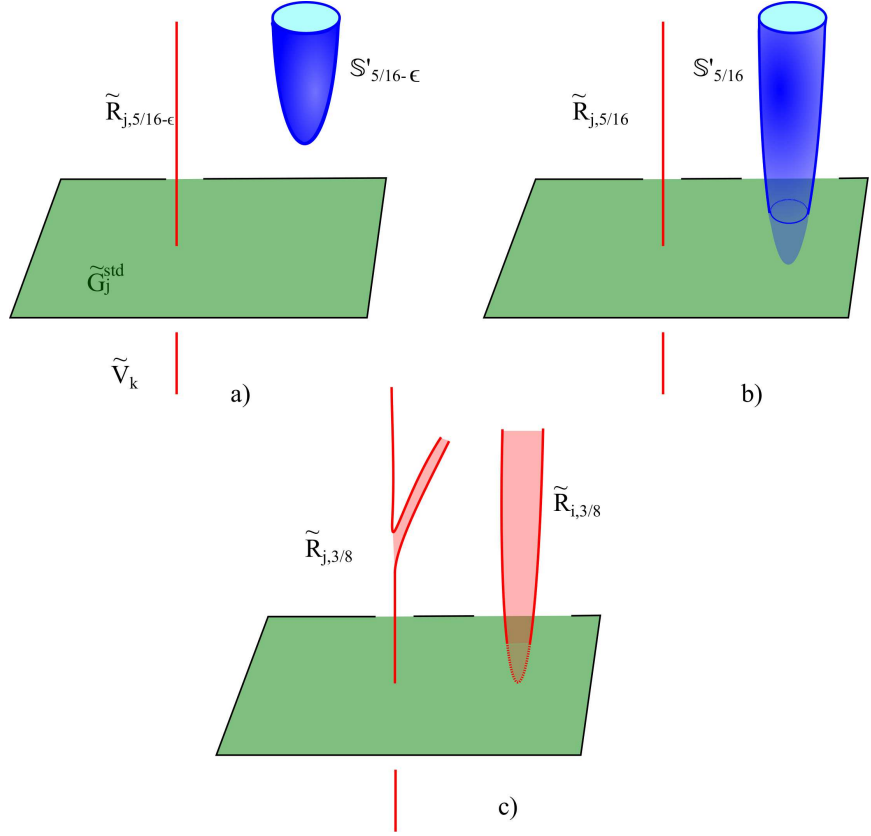}$  
\end{tabular}
 \caption[(a) X; (b) Y]{\label{finger move}
 \begin{tabular}[t]{ @{} r @{\ } l @{}}
Preparing for and doing the Finger move, as seen in $\tilde V_k$\end{tabular}}
\end{figure}
\vskip 8pt
\noindent ii) \emph{Do the finger moves}: For $t\in [3/8-\epsilon, 3/8)$ isotope $\tilde \mR_t$ according to its arms, hands and fingers.  See Figure \ref{finger move} c).  Here $\BS'_{3/8}$ is not shown.  
\setlength{\tabcolsep}{60pt}
\begin{figure}
 \centering
\begin{tabular}{ c c }
 $\includegraphics[width=4.5in]{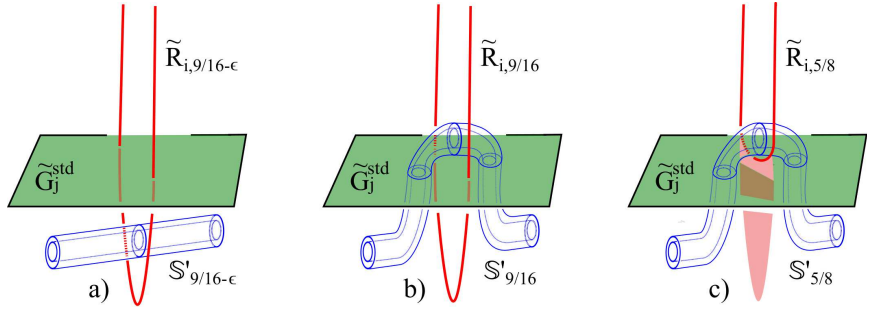}$  
\end{tabular}
 \caption[(a) X; (b) Y]{\label{whitney move}
 \begin{tabular}[t]{ @{} r @{\ } l @{}}
Preparing for and doing the Whitney move, as seen in $\tilde V_k$\end{tabular}}
\end{figure}

\vskip 8pt
\noindent iii)\emph{ Prepare for the Whitney moves}:  If $\tilde w_{ij}\subset\tilde V_k$ is a Whitney disc between $\tilde R_{i,1/2}$ and $\widetilde \Gs_j$, then $\tilde w_{ij}\cap \BS'_{1/2}$ is eliminated by an isotopy $\BS_t', t\in [(9/16-\epsilon, 9/16)$ as in Figures \ref{whitney move} a), b) which shows a 3D-slice.  Note that each component of $\tilde w_{ij}\cap \BS'_{1/2}$ gives rise to two torus components of $N^1(\widetilde\mGs_j)\cap \BS'_{9/16}$.  In that figure both $\BS'_t$ and $\tilde R_{i,t}$ extend into the past and the future.

\vskip 8pt
\noindent iv) \emph{Do the Whitney moves}: For $t\in (5/8-\epsilon, 5/8)$ isotope $\tilde R_t$ according to its Whitney discs to obtain $\tilde R_{5/8}$.  See Figure \ref{whitney move} c).  The shaded region indicates the local projection of $\tilde R_{i, 5/8}$ to the present.  
\vskip 8pt
\noindent v) \emph{Do the glissade}:  We can assume that $\BS'_{5/8} \cap N(\widetilde\mGs)$ is a union of normal discs.  If $C\subset N^1(\widetilde \Gs_j)$ is the boundary of one such normal disc, then under the glissade $\lambda(C)$ is capped off by a disc which is the union of an annulus $\subset N(\tilde w_j)$ and a copy of $\tilde R^0_{j, 5/8}$. This completes the construction of Step 2.

\vskip 10pt
Since each $\tilde \mR_{5/8}$ is isotopic to $\tilde \mR_{1/4}$, Theorem 10.1 \cite{Ga1} applies to any finite set of components of $\tilde \mR_{5/8}$.  Applying $\tilde \lambda$ to this isotopy we obtain the following result.

\begin{lemma} \label{pie cocore isotopy} For any finite set $J\subset \BZ$, there exists an ambient isotopy $\kappa_t$ of $\tilde V$ fixing $\cup_{j\in J} N(\tilde \omega_j)$ pointwise such $\kappa_0=\id$ and $\kappa_1(\cup_{j\in J} \tilde R^0_{5/8})=\cup_{j\in J} \tilde R^0_{1/4}$.\qed\end{lemma}

\begin{remark} Actually the construction of Step 2 can be done at the $V$ and $V_k$ level so that there is a $\BZ$-equivariant ambient isotopy of $\tilde \mR^0_{5/8}$ to $\tilde \mR^0_{1/4}$ that fixes $N(\tilde\Omega)$  pointwise. \end{remark} 

\vskip 2pt

The next two steps introduce technology to visualize the passage of $\BS_0$ to $\BS_1$ more concretely.  In Step 5 we redo Step 2 from this change of perspective, using a slightly different isotopy extension $\eta_t$, after which we see that $\BS_1$ has a CS-presentation.  Step 6 gathers the data recorded in Step 5 to show that we actually have a FWCS-presentation.
\vskip 10pt

\noindent\textbf{Step 3:} Learn to work with anellini discs.
\begin{definition} An \emph{anellini disc} is a $\solidtorus$ viewed as an $\epsilon N^1(D)$ where $D$ is an embedded 2-disc in a 4-manifold.  Here $\epsilon$ is very small, so that the fibers are the boundaries of tiny normal discs.    In what follows the $\epsilon$ will be suppressed.\end{definition}

\begin{remark}An anellini disc is the 3-dimensional analogue of an annulus viewed as a thin tube which makes the surface look 1-dimensional.  Anellini discs make part of the 3-manifold look 2-dimensional.  As it is often beneficial to work with surfaces having subsurfaces a union of tubes, anellini discs are useful for working with 3-dimensional submanifolds of 4-manifolds.  In particular a reimbedding of a solid torus corresponding to an anellini disc can be determined by a reimbedding of the anellini disc.\end{remark}

Our next definition expresses the well-known operation of Figure \ref{acompression} in terms of carvings, surgeries and anellini discs.  See for example the knotted 3-spheres of  \S8 \cite{BG}.
\setlength{\tabcolsep}{60pt}
\begin{figure}
 \centering
\begin{tabular}{ c c }
 $\includegraphics[width=4.5in]{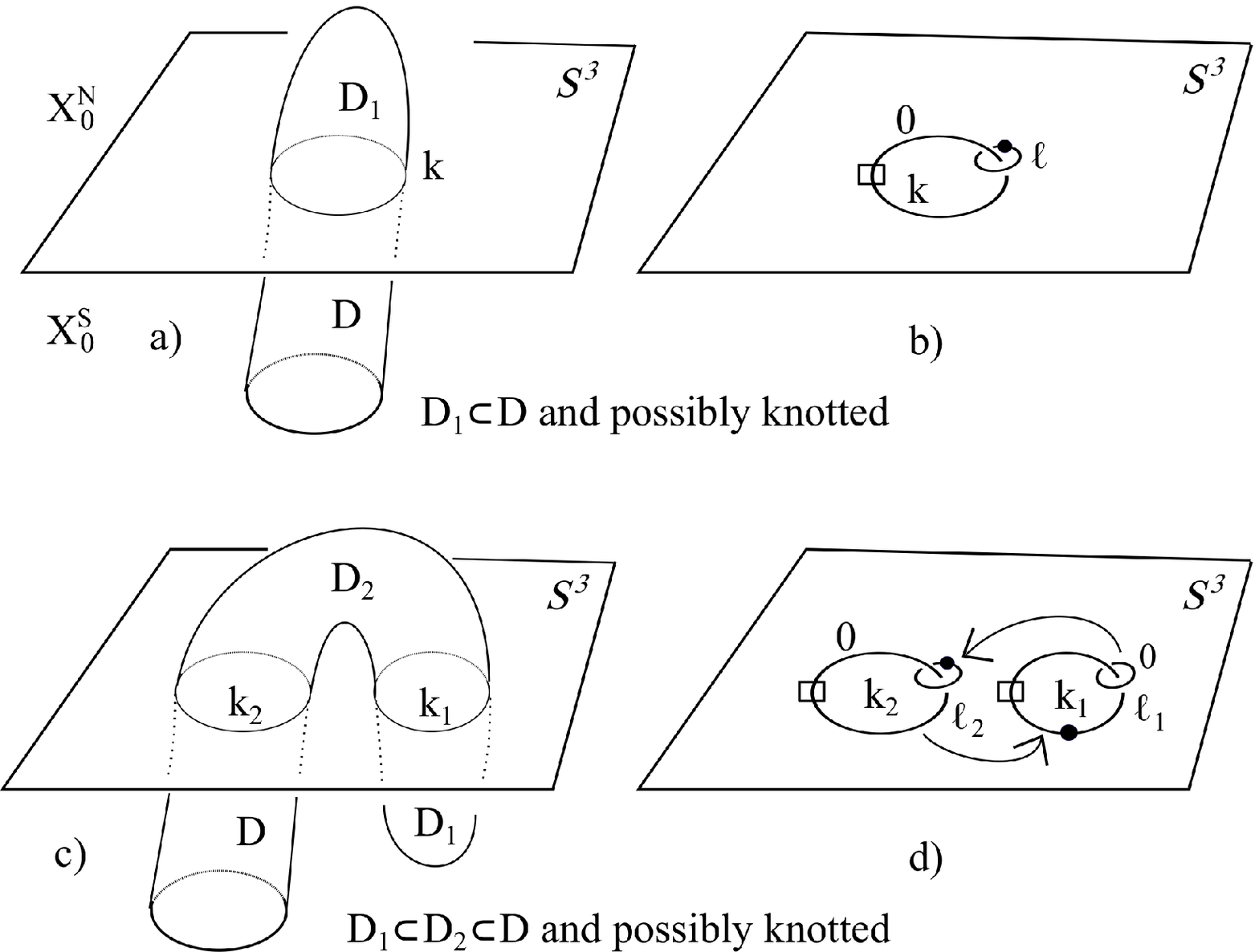}$  
\end{tabular}
 \caption[(a) X; (b) Y]{\label{acompression}
 \begin{tabular}[t]{ @{} r @{\ } l @{}}
Anellini compression\end{tabular}}
\end{figure}
\begin{definition} \label{anellini} Figure \ref{acompression} a) schematically shows a disc $D$ intersecting $\BS_0$ in the knot $k$, which bounds the subdisc $D_1\subset X^N_0$, the northern 4-ball.  $\BS_0$ is isotopic to $S^3_1$, the result of embedded surgery along a 2-handle with core $D_1$ and carving the standard disc $D_\ell\subset X_0^S$ bounded by the linking circle $\ell$. See Figure \ref{acompression} b). We call this operation \emph{anellini compression}.  While it creates two anellini discs, we single out  $D_\ell$ as the \emph{anellini disc} resulting from this operation.  See Remark \ref{ztheta invariance}.  Note that $|D_\ell\cap D|=1$ where $D_\ell$ is viewed as a disc. If instead we had $D_1\subset X^S_0$, then $S^3$ is isotopic to a presentation with $D_1$ carved and with embedded surgery along $D_\ell$ with again $D_\ell$ the anellini disc.  While we described anellini compression from the $X^S_0$ perspective, we just as well could have used $X^N_0$, where carving and embedded surgery switch roles.   

Figure \ref{acompression} c) shows the case when $D\cap \BS_0 =k_1\cup k_2$.  Here $S^3$ is isotopic to the presentation of Figure \ref{acompression} d), where $\ell_1, \ell_2$ are the linking circles.  Note that a $\bullet$ (resp. 0) corresponds to a carving along a disc whose boundary germ points into $X^S_0$ (resp. $X^N_0$).  While the induced framings on linking circles are always 0-framed, the induced framings on the knots are not necessarily 0-framed even though they might be labelled with a 0.  In our case $\BS_0$ is first carved along the disc $D_1$ and then embeddedly surgered along the disc $D_2$.  It is also carved along the standard disc $D_{\ell_2} $ followed by an embedded surgery along the disc $D_{\ell_1}$.  Here $D_{\ell_1}$ consists of a tube that starts at $k_1$ and follows an arc in $D_2$ that is capped off with a copy of $D_{\ell_2}$. We say that $D_2$ \emph{nests} in $D_1$ and $D_{\ell_1}$ \emph{nests} in $D_{\ell_2}$, i.e. it enters the hollow created by the carving or embedded surgery.  Nesting is combinatorially recoded as in Figure \ref{acompression} d).  If $\partial D\cup D_1\subset X^N_0$, then we have the similar situation with 0's replaced by bullets and vice versa.  We call $D_{\ell_1}, D_{\ell_2}$ the resulting anellini discs.

In general a disc $D$ may intersect $\BS_0$ in a finite set $k_1, \cdots, k_n$ of simple closed curves.  Here $\BS_0$ is isotopic to $S^3_1$ which is represented by a set of carvings and embedded surgeries along the $k_i$'s and their linking circles and $S^3_1$.  Figure \ref{Nesting} b) shows a \emph{nesting diagram} that schematically indicates the carvings and surgeries for the disc $D$ of Figure \ref{Nesting} a) with $\partial D\subset \inte(X^S_0)$.  The boxes indicate that the knots can be knotted and linked with each other.\end{definition}

\setlength{\tabcolsep}{60pt}
\begin{figure}
 \centering
\begin{tabular}{ c c }
 $\includegraphics[width=4.5in]{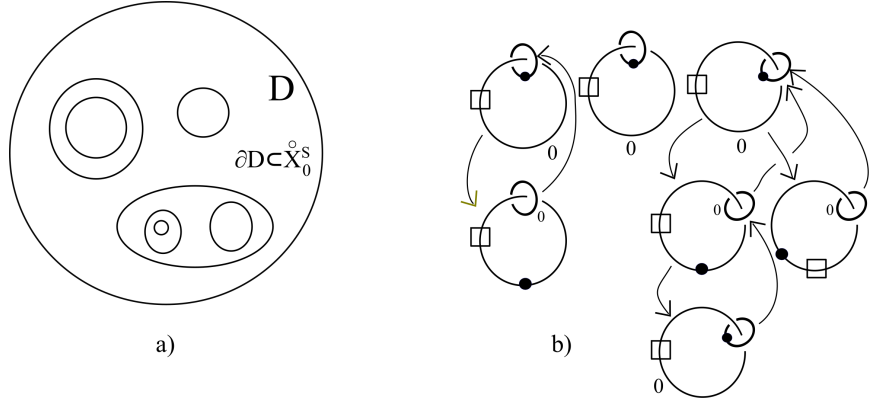}$  
\end{tabular}
 \caption[(a) X; (b) Y]{\label{Nesting}
 \begin{tabular}[t]{ @{} r @{\ } l @{}}
The Nesting Partial Order\end{tabular}}
\end{figure}

\begin{remarks} \label{ztheta invariance}  We single out the $D_{\ell_i}$'s as the resulting anellini discs because in what follows, e.g. the next result,  it is these discs that will under go transformations as surfaces.

If $D_{\ell_1}$ nests in $D_{\ell_2}$, then they are disjoint as anellini discs, however when viewed as discs, $D_{\ell_2}$ is a subdisc of $D_{\ell_1}$.  In what follows, the distinction \emph{as an anellini disc} vs. \emph{as a disc} should be clear from context.  \end{remarks}

To complete Step 3, we now show how to glissade from the point of view of anellini discs.  Notation is as in v) of Step 2).  

\begin{lemma}  \label{anellini glissade}  If $\BS'_{5/8}\cap N(\widetilde\mGs)$ is contained in anellini discs $A'_{5/8}$ and the intersection with these anellini discs, viewed as discs, is a union of fibers of $N(\widetilde\mGs)$, then $\tilde \lambda(A'_{5/8})\setminus \inte(N(\Omega)) \subset\tilde V $ is a union of discs denoted $A_{5/8}$ with open subdiscs removed, one for each component of $A'_{5/8}\cap N(\widetilde\mGs)$.  The glissade creates the anellini discs $A_1$ obtained by capping off the boundary circles of $A_{5/8}$ with annuli in $N(\tilde\Omega)$ together with copies of components of $\tilde \mR^0_{5/8}$.

Equivalently, first obtain $A'_1$ from $A'_{5/8}$ by tubing off each intersection with $\widetilde\mGs$ with a parallel copy of a component of $\widetilde \mR_{5/8}$ and then let $A_1=\lambda(A_1')$. \qed\end{lemma}

\begin{remark}\label{multi glissade}  The glissade may involve nested anellini discs which coincide, as discs, along $N(\widetilde\mGs)$, in which case boundary circles of $A_{5/8}$ are capped off by coinciding copies of components of $\tilde \mR^0_{5/8}$.  In particular $A_1$ will be correspondingly nested.\end{remark}

\vskip 10 pt
\noindent\textbf{Step 4:}  Change our perspective of $N(\tilde \mD)$ to pies with fillings and crust.
\vskip 8pt

The advantage of this point of view is that the finger move and glissade are particularly simple and it gives a way of seeing how the transformation $\BS_0$ to $\BS_1$ can be realized via regular homotopy.

\vskip 8pt
\begin{definition} \label{pie and filling} We now view the component of  $N(\tilde\Omega\cup \tilde R^0_{1/4})$ associated to $\tilde w_j$ as the thin \emph{pie} $D^2_j\times D^2_\epsilon:= P_j$.  See Figure \ref{pie} a).  This is just Figure \ref{section fiber} b) flattened out, with a corner added,  $\tilde R^0_{1/4}$ thickened and the cylindrical coordinates $x,y,r,\theta$ rotated.

Here $N(\tilde R^0_{j,1/4})$ corresponds to $D^2_j\times \frac{1}{2}D^2_\epsilon$ and is called the \emph{filling} and we call $D^2_j\times (D^2_\epsilon\setminus\inte(\frac{1}{2} D^2_\epsilon))$ the \emph{crust}.    Circles of the form $u\times \partial D^2_\epsilon$, $ u\in D^2_j $ are $\lambda$ images of fibers of $N^1(\widetilde \Gs_j)\setminus \inte(N(\widetilde\Rs_j))$ and are called \emph{crust fibers}.  Also discs of the form $u\times \frac{1}{2}D^2_\epsilon$ are parallel copies of $\tilde R^0_{1/4}$.  We call discs of the form $u\times D^2_\epsilon$ \emph{pie cocores} and discs of the form $D^2_j\times v, v\in \frac{1}{2} D^2_\epsilon$ are called \emph{sections of the filling}.  See Figure \ref{pie} b).  \end{definition}

\setlength{\tabcolsep}{60pt}
\begin{figure}
 \centering
\begin{tabular}{ c c }
 $\includegraphics[width=4.5in]{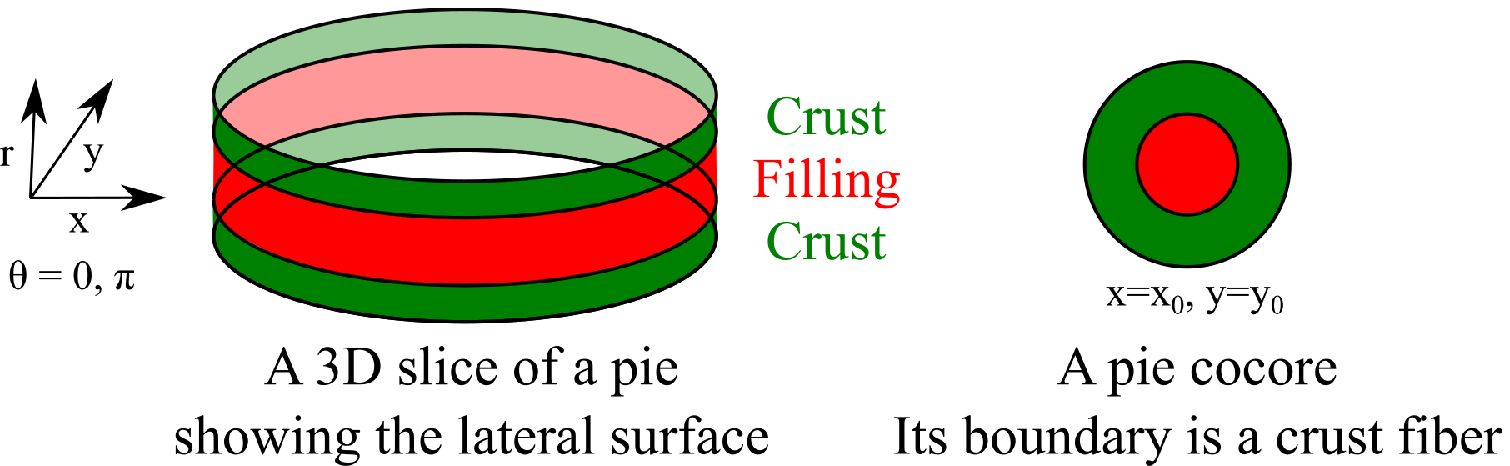}$  
\end{tabular}
 \caption[(a) X; (b) Y]{\label{pie}
 \begin{tabular}[t]{ @{} r @{\ } l @{}}
A Pie\end{tabular}}
\end{figure}

\begin{remark}\label{pie finger} \emph{Finger moves from the pie view}:  Figure \ref{finger pie} shows various 3D slices, before and after, of the $\lambda$ images $\subset \tilde V$ of a finger move.  Figures \ref{finger pie} a), b) show the $\theta = 0, \pi$ slices, while Figures c), f) and d), g) and e), h)  respectively show $x=-a$, $x\in (-a, a)$ and $x=a$ slices.  Note that figures b), f), h) have the interiors of the pie deleted from the fingers, i.e. two open discs.  Figures c)-g) correspond to a finger move with corners.  In Figures a), b)  solid lines denote the part of $\tilde\mR^0_t$ in the present, while the shading denotes the parts in the past or the future.  From this point of view we see  that a finger move in $\tilde V$ is the result of a regular homotopy followed by the deletion of two open discs.  \end{remark}

\setlength{\tabcolsep}{60pt}
\begin{figure}
 \centering
\begin{tabular}{ c c }
 $\includegraphics[width=4.0in]{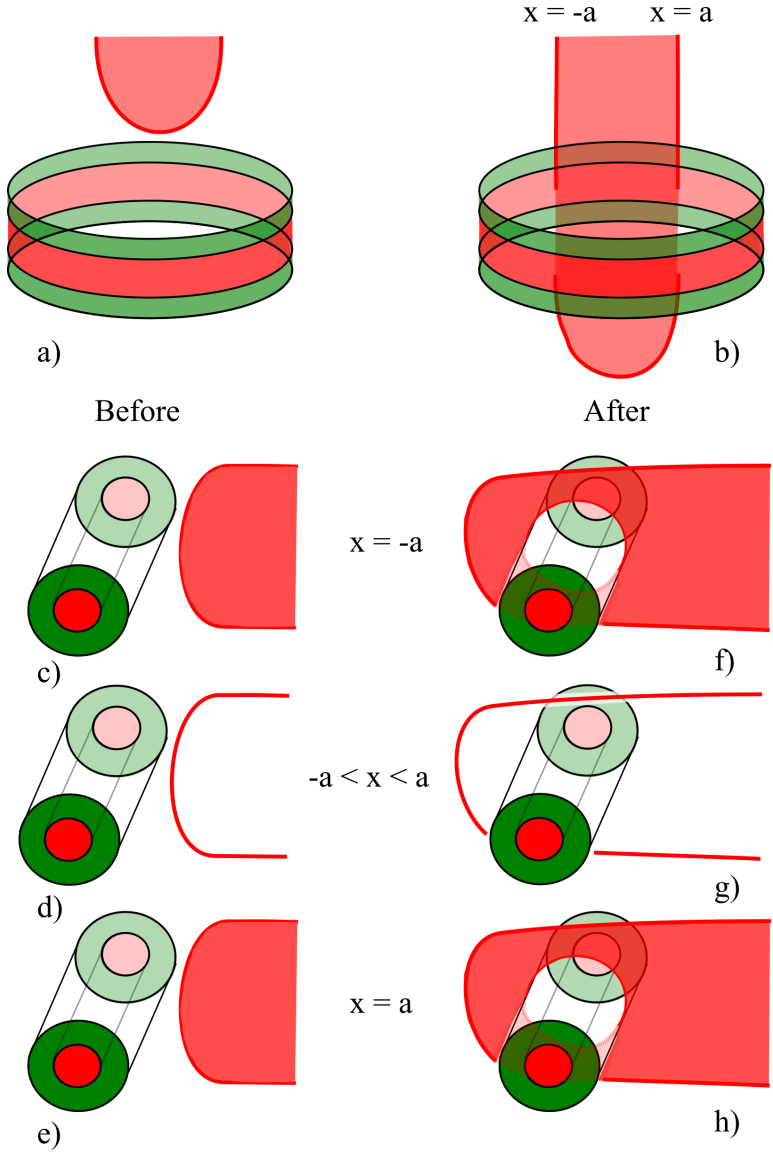}$  
\end{tabular}
 \caption[(a) X; (b) Y]{\label{finger pie}
 \begin{tabular}[t]{ @{} r @{\ } l @{}}
Finger Move: Pie View\end{tabular}}
\end{figure}

\begin{remark} \label{pie glissade}  \emph{The glissade from the pie view}:  Since $\tilde R^0_{j,5/8}$ coincides with $\tilde R^0_{j,1/4}$ near $\partial \tilde R^0_{j,1/4}$ we can assume that $N(\tilde R^0_{j,5/8}):=D^2_j\times \tilde R^0_{j,5/8}$ and fiberwise coincides with $N(\tilde R^0_{j,1/4}):=D^2_j\times \frac{1}{2}D^2_\epsilon$ near $D^2_j\times \partial \frac{1}{2}D^2_\epsilon$.  If $A_{5/8}$ are $\lambda$ images of anellini discs that intersect $\partial(D^2_j\times D^2_\epsilon)$ in $K_j\times \partial D^2_\epsilon$, then the glissade gives $$A_1=A_{5/8}\cup_j(K_j\times (D_\epsilon^2\setminus \inte \frac{1}{2} D^2_\epsilon))\cup_j (K_j\times \tilde R_{j,5/8}^0).$$  See Figure \ref{glissade pie}.  Here $A_{5/8}$ intersects $P_j$ in a single crust fiber.  Note that $\tilde R_{j,5/8}^0$ and more generally $\tilde \mR^0_{5/8}$, may intersect $P_j$ in sections of the filling.  These are not shown in the figure. This completes Step 4. \end{remark}

\setlength{\tabcolsep}{60pt}
\begin{figure}
 \centering
\begin{tabular}{ c c }
 $\includegraphics[width=4.5in]{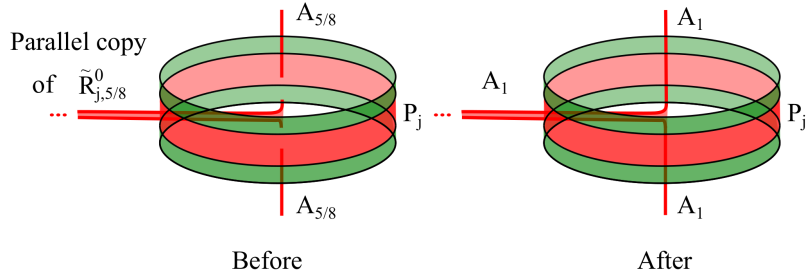}$  
\end{tabular}
 \caption[(a) X; (b) Y]{\label{glissade pie}
 \begin{tabular}[t]{ @{} r @{\ } l @{}}
The Glissade: Pie View\end{tabular}}
\end{figure}
\noindent\textbf{Step 5:}  Construct $\BS_1$ from the pie perspective using anellini compressions. 
\vskip 8pt
In this section we go through the construction of Step 2 using the technology developed in Steps 3 and 4.  We also record as \emph{data memos} the  data needed for our $F|W$-carving/surgery presentation.   \vskip 8pt
\setlength{\tabcolsep}{60pt}
\begin{figure}
 \centering
\begin{tabular}{ c c }
 $\includegraphics[width=4.5in]{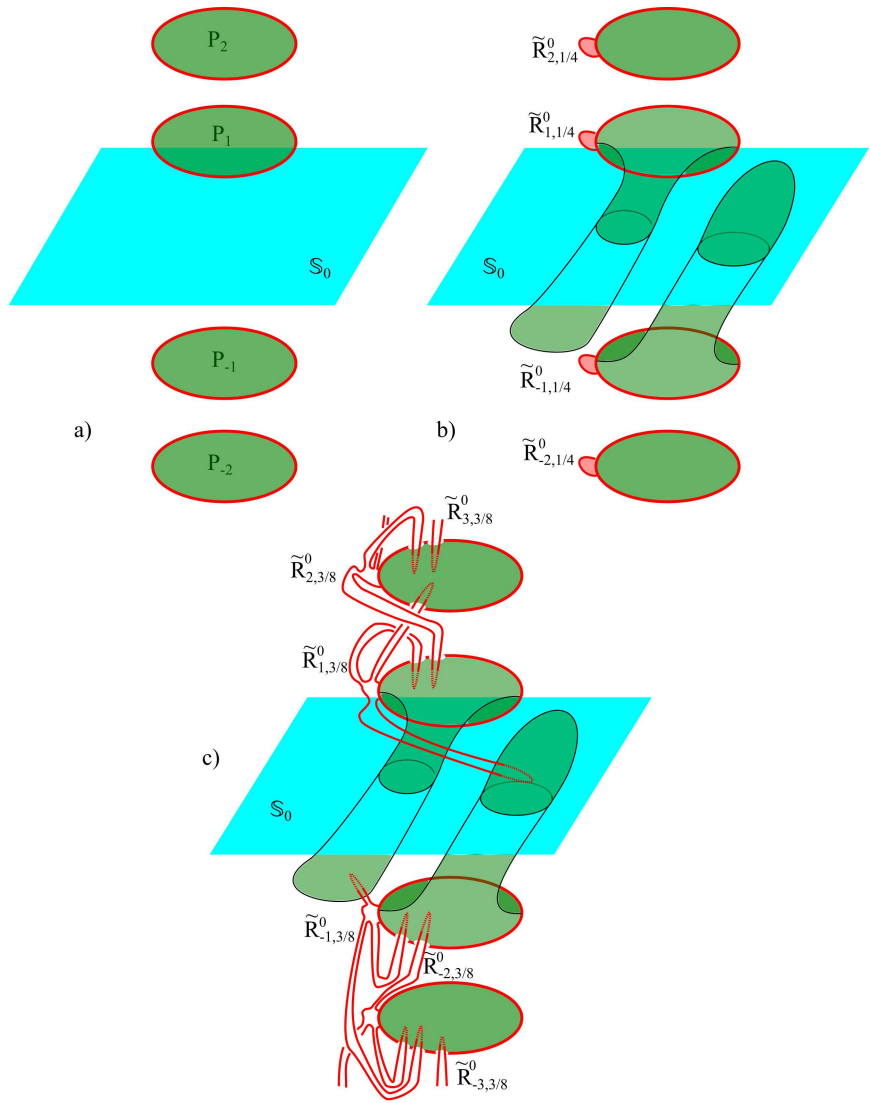}$  
\end{tabular}
 \caption[(a) X; (b) Y]{\label{finger pie move}
 \begin{tabular}[t]{ @{} r @{\ } l @{}}
Finger move in $\tilde V$: Pie view\end{tabular}}
\end{figure}

\setlength{\tabcolsep}{60pt}
\begin{figure}
 \centering
\begin{tabular}{ c c }
 $\includegraphics[width=5in]{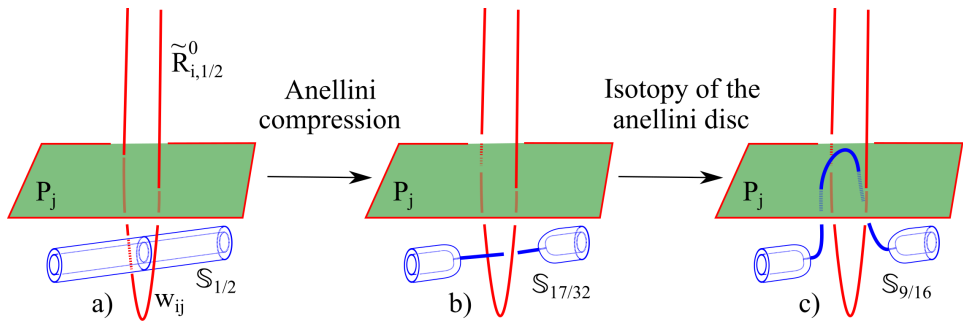}$  
\end{tabular}
 \caption[(a) X; (b) Y]{\label{whitney prepare}
 \begin{tabular}[t]{ @{} r @{\ } l @{}}
Preparing for the Whitney move\end{tabular}}
\end{figure}

\noindent i) \emph{Prepare for the finger moves}:  Here we take the point of view of finger moves of $\widetilde \Gs_j$ into $\BS_0'$ rather than the other way around with the corresponding modifications in $\tilde V$.  We abuse notation by continuing to call the modified $\widetilde \Gs_j$ by the same name.  For each $j>0$ (resp. $j<0$) for which there is an arm from $\widetilde \Rs_i$ to $\widetilde \Gs_j$ with $i<0$ (resp. $i>0$), do a finger move from $P_j$ into $\BS_{1/4}=\BS_0$.   See Figures \ref{finger pie move} a), b).  In Figure b) $\tilde\mR^0_{1/4}$ has been isotoped to pop slightly outside of $P_j$.  The dark lines indicate the present and the shading the projection to the present.  To keep indices consistent we let $\BS_{1/4}$ denote $\BS_0$ even though it has not been modified.

\begin{data memo} \label{memo1} a) Viewed in $\BS_0, P_j\cap \BS_{1/4}$ is an unknotted solid torus $N^3(g_j)$.  The knots $g_j$ arising from this preparatory move are unknotted and unlinked in $\BS_0$.

b) Call a pie \emph{plunged} if it intersects $\BS_0$.  Let $P_{n_1}, \cdots, P_{n_p}$ denote the plunged pies with $n_i>0$ and $P_{s_1}, \cdots, P_{s_q}$ those  with $s_j<0$ and let  $g_{n_1}, \cdots, g_{n_p}, g_{s_1}, \cdots, g_{s_q}$ denote the corresponding knots.

c) If $j<0$ (resp. $j>0$), then $P_j\cap X^N_0=\sigma_j$ (resp. $P_j\cap X^S_0$) is the product  of a standard 2-disc $D_{g_j}$ with pie cocores, where  $\partial D_{g_j}=g_j$.  In particular, those 2-discs respectively in $X_0^N$ and $X_0^S$ are unknotted and unlinked.\end{data memo} 
\vskip 8pt
\noindent ii) \emph{Do the finger moves}:  Figure \ref{finger pie move} c) shows $\tilde\mR^0_{3/8}$ after finger moves of $\tilde R^0_{i,1/4}$ into $P_{i-1}, P_i, P_{i+1}$, $i\in \BZ$,  where indices need to be adjusted as in Notation \ref{rereindex}.  The figure only shows the intersection of $\tilde\mR_{3/8}^0$ with the present.  Here we modify the time indices to 1/2 so in particular we rename $\BS_{1/4}$ as $\BS_{1/2}$.

\vskip 8pt
\noindent iii) \emph{Prepare for the Whitney moves}:  To do this we first do anellini compressions to obtain $\BS_{17/32}$ and then move the anellini discs off of their Whitney discs to obtain $\BS_{9/16}$.

Let $\tilde w_{ij}$ denote a Whitney disc from $\widetilde R_{i,1/2}$ to $\widetilde \Gs_j$ and $w_{ij}:=\lambda(\tilde w_{ij})\subset \tilde V$.  The number of such discs is the number of fingers in the corresponding hand and $\cup w_{ij}$ will denote the union of such discs.  For each $w_{ij}$ do anellini compressions to modify $\BS_{1/2}$ to $\BS_{17/32}$.   While anellini compression is realizable by regular homotopy, the resulting $\BS_{17/32}$ need not be embedded.  See Data Memo \ref{memo3.1} d).   

\begin{data memo} \label{memo3.1} a)  Each $w_{ij}$ gives rise to a nesting diagram as in Figure \ref{ij,Family}, which shows the case of a hand with two fingers, where in a) (resp. b))  $i<0$ (resp. $i>0$).  

b)  $\BS_{17/32}\cap w_{ij}$ is contained in the anellini discs which span the linking circles of the knots of the diagram.  These discs induce 0-framings on the linking circles.

c)  If $\gamma$ is a component of $\BS_{1/2}\cap w_{ij}$, then let $D_\gamma\subset w_{ij}$ denote the disc bounded by $\gamma$ and $\tilde D_\gamma=\lambda^{-1}(D_\gamma)$.  We let $\gamma$ also both denote this knot viewed in $\BS_0$ and the corresponding curve in $w_{ij}$. We denote its linking circle by $\gamma'$.  Let $D_{ij}:=\cup_\gamma D_\gamma \subset(\cup w_{ij})$ and $\tilde D_{ij}=\lambda^{-1} (D_{ij})$.  Give the $\gamma$'s the induced partial order coming from the nesting diagram knots, i.e. ordered by inclusion with innermost ones minimal and give the $\gamma'$'s the opposite partial order induced from the nesting of the $\gamma$'s. We partially order the anellini discs according to the partial order of their linking circles.  Let $D^\gamma_{ij}=D_\gamma\setminus\cup_{\beta<\gamma} \inte(D_\beta)$ and $\tilde D^\gamma_{ij}=\lambda^{-1} (D^\gamma_{ij})$.

d) Each point of $\tilde D^\gamma_{ij}\cap \widetilde\Rs_\ell$ gives rise to a section $\subset D^\gamma_{ij}\cap P_\ell$ of $P_\ell$'s filling.  See Figure \ref{section}.  If $P_\ell$ is a plunged pie, then $\BS_{17/32}$ will not be embedded.

e) Each point of $(\cup \tilde w_{ij}\setminus \tilde D_{ij})\cap \widetilde\Rs_\ell$ gives rise to a section $\subset (\cup w_{ij}\setminus D_{ij})\cap P_\ell$ of $P_\ell$'s filling.

f)  $\cup w_{ij}$ nests once in each $D_\gamma$, for maximal $\gamma\subset \cup w_{ij}$.

g) We call the knots arising  from these anellini compressions \emph{southern knots} (resp. \emph{northern knots}) when $i<0$ (resp. $i>0$). A given southern (resp. northern) knot is  denoted $S_{ij}$ (resp. $N_{ij})$ if it arises from a $w_{ij}$ Whitney disc.  When we say a knot (resp. linking circle) is \emph{northern maximal} or \emph{minimal} we mean that it is maximal or minimal among the northern knots (resp. linking circles).  Similarly we may refer to \emph{southern minimal} etc. knots or linking circles.  \end{data memo} 
\setlength{\tabcolsep}{60pt}
\begin{figure}
 \centering
\begin{tabular}{ c c }
 $\includegraphics[width=4.5in]{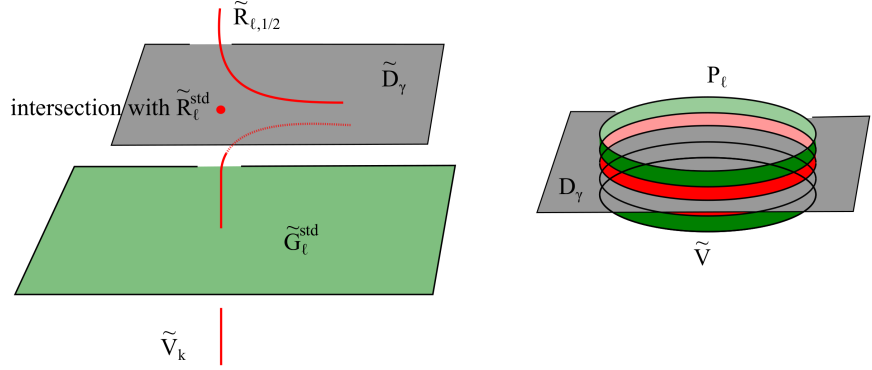}$  
\end{tabular}
 \caption[(a) X; (b) Y]{\label{section}
 \begin{tabular}[t]{ @{} r @{\ } l @{}}
A Section of a Pie Filling arising from an intersection with $\widetilde \Rs_\ell$\end{tabular}}
\end{figure}

Next isotope the northern and southern minimal anellini discs off of the Whitney discs to obtain $\BS_{9/16}$ at the cost of creating two oppositely signed intersections with $P_j$, viewed as a disc, for each such anellini disc.  See Figure \ref{whitney prepare}  which shows one anellini disc nested in another and compare with Figures \ref{whitney move} a) and b).

\begin{data memo}\label{memo3.2} a) Here each northern or southern minimal anellini disc intersects exactly one $P_j$ in two pie cocores.  To follow Step 2 exactly we need to  remove the interior of these cocores.  The advantage of keeping them is two fold.  We   can continue to work with the anellini discs as discs and we can extend the regular homotopy of $\BS_0$ to $\BS_{17/32}$ to a regular homotopy to $\BS_{9/16}$.  Note that this operation in conjunction with Data Memo \ref{memo3.1} d) may create additional self intersections of $\BS_{9/16}$.

b) Since an isotopy of an anellini disc extends to an isotopy of those discs nested inside, the  isotopy of the minimal anellini discs extends to all the anellini discs intersecting $w_{ij}$ and hence moves all of them off of $w_{ij}$.

c) An intersection of a minimal northern or southern anellini disc with $P_j$ is $\subset \sigma_j$  if and only if it arises from a $w_{ij}$ with $i<0$ and $j>0$ or $i<0$ and $j>0$.\end{data memo}

\noindent iv) \emph{Do the Whitney moves}:  

\begin{data memo} \label{memo4}  a)  $\tilde R^0_{i, 5/8}$ nests twice with opposite sign through each northern and southern maximal knot of the form $S_{ij}$ or $N_{ij}$.

b) Each point of $(\cup\tilde w_{ij}\setminus \tilde D_{ij})\cap \widetilde \Rs_\ell$ induces two oppositely oriented filling sections of $P_\ell\cap \tilde R^0_{i, 5/8}$.  \end{data memo}

\vskip 10pt

\noindent v) \emph{Do the glissade}:  First prepare for the glissade by doing anellini compressions to the $g_j$'s along the discs $D_{g_j} $ in the plunged pies to obtain $\BS_{21/32}$.   
 See Figures  \ref{glissade final} a), b).   It depicts the anellini compression along $D_{g_j}$ from the $P_j$ point of view, i.e. here $P_j$ looks flat rather than plunged.  In this figure only the intersections with the present are shown.  Figure \ref{glissade final} a) depicts two anellini discs while b) also includes a third, notably $D_{g'_j}$. This preparatory operation requires isotoping those anellini discs that intersect $D_{g_j}$ as in Figure \ref{glissade final} b).
 
 \vskip 8pt
Following the glissade preparatory move, $\BS_{21/32}$ intersects the pies in subdiscs of the anellini discs and these subdiscs are pie cocores.  The glissade involves replacing the fillings of these cocores with parallel copies of components of $\tilde \mR^0_{5/8}$ as in Remark \ref{pie glissade}.  The result transforms $\BS_{21/32}$ to $\BS_1$.  

All the self intersections of $\BS_{21/32}$ arise from intersections of these pie cocore fillings  with filling sections.  Let $\BS_{21/32}'$ denote $\BS_{21/32}$ with these pie cocore fillings removed.  Under the glissade these discs, are replaced by parallel copies of components of $\tilde\mR^0_{5/8}$,  which are embedded, pairwise disjoint and intersect $\BS_{21/32}'$ only along their boundaries.   It follows that the resulting $\BS_1$ is embedded. 

A pie cocore filling shares its boundary with a parallel copy of some $\tilde R^0_{j,5/8}$.  It follows by Lemma \ref{pie cocore isotopy} that  it is isotopic to $\tilde R^0_{j,5/8}$    fixing a neighborhood of its boundary pointwise, It follows that the glissade does not change the induced framings of the boundary of the annellini  discs and that the glissade can be realized by a regular homotopy.

 \setlength{\tabcolsep}{60pt}
\begin{figure}
 \centering
\begin{tabular}{ c c }
 $\includegraphics[width=4.5in]{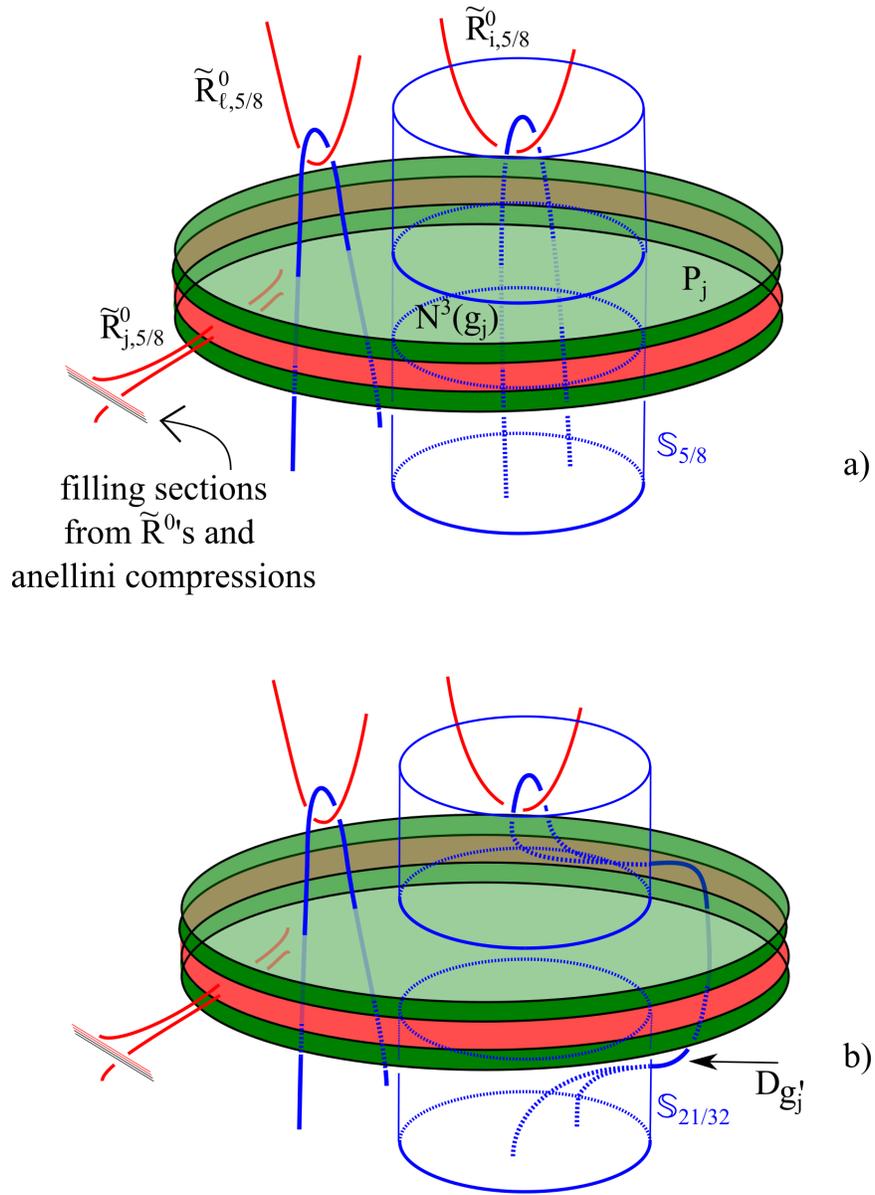}$  
\end{tabular}
 \caption[(a) X; (b) Y]{\label{glissade final}
 \begin{tabular}[t]{ @{} r @{\ } l @{}}
An Anellini compression after the Whitney moves\end{tabular}}
\end{figure}

\begin{data memo} \label{memo5}  a)    $\tilde R^0_{i,5/8}$ will nest $g_\ell$ if and only if $(\cup_j\tilde w_{ij}\setminus (\cup_j \tilde D_{ij})) \cap \widetilde \Rs_\ell\neq \emptyset$ and for $\gamma \subset w_{ij}$, $\gamma$   will nest in $g_\ell$ only if $\tilde D^\gamma_{ij}\cap \widetilde \Rs_\ell\neq \emptyset$.

b) Following the glissade, $g'_i$ will nest exactly those knots nested by $\tilde R^0_{i, 5/8}$, i.e. all the maximal knots of the form $S_{ij}$ or $N_{ij}$ plus those $g_\ell$ from a) above.  

c) If $j>0$ (resp. $j<0$), then the linking circle $g'_j$ lodges exactly the southern (resp. northern) minimial linking circles $S'_{ij}$ (resp. $N'_{ij}$). 

d) If $i<0$ and $j>0$, then a southern minimal  linking circle $S'_{ij}$ will only nest in $ g'_j$.  If instead $j<0$, then $S'_{ij}$ will nest in the maximal  $S_{ji}$ knots and those  $g_\ell$'s nested by $\tilde R^0_{i,5/8}$.

e) If $i>0$ and $j<0$, then a northern minimal  linking circle $N'_{ij}$ will only nest in $g'_j$.  If instead $j>0$, then $N'_{ij}$ will nest exactly  the maximal $N_{ji}$ knots and  those $g_\ell$'s nested by $\tilde R^0_{i,5/8}$.

\end{data memo}

This completes Step 5.  
\vskip 10pt

\noindent\textbf{Step 6:} Organize the data.
By construction $\BS_1$ is obtained from $\BS_0$ by a carving/surgery presentation and the link $\mL$ is a union of knots $L_k$ and linking circles $L'$.    The various data memos record the partial order details which we summarize now and show that we have a $F|W$-carving/surgery presentation.  

The passage of $\BS_0 $ to $\BS_1$ involved modifications of $\BS_0$ while: A) preparing for the Whitney moves, B) preparing for the glissade and C) doing the glissade.  The Whitney preparatory moves create the northern and southern knots and their linking circles.  Each northern or southern knot corresponds to a $\gamma\subset w_{ij}$ and the disc $D_\gamma\subset w_{ij}$ it bounds.  Thus the nesting relations among the northern or southern knots arise during A) and give the relations described in Definition \ref{fwcs} 1b) and 1c).   These $D_\gamma$'s may contain filling sections which in turn create nestings into $D_{g_j}$'s during B), which are recorded in the second and forth paragraphs of Definition \ref{fwcs} 2a).

The $g_j$'s were defined in Data Memo \ref{memo1} as were their spanning discs $D_{g_j}$, though they and their linking circles weren't created until B). These knots and their linking circles comprise $B_L$.  By construction the $g_j$'s are unknotted and unlinked in $\BS_0$ and the intersections of the $D_{g_j}$'s with either $X^N_0$ or $X^S_0$ are also unknotted and unlinked.   A $D_\gamma$ will nest a $D_{g_j}$ only if it intersects $P_j$ in filling sections which occurs only if $\tilde D_\gamma\cap \widetilde\Rs_j\neq\emptyset$.  See Data Memo's \ref{memo3.1} d).  These nestings give a possibly proper subset of the nestings  described in the second and fourth paragraphs of Definition \ref{fwcs} 2a).  

The spanning disc $D_{\gamma'}$ of the linking circle $\gamma'$ of a southern or northern knot $\gamma$ is constructed in 3 steps.  First,  we have $D^1_{\gamma'}$ the standard spanning disc constructed in the anellini compression operation during A). Second, the northern and southern minimal $D^1_{\gamma'}$'s are isotoped during A) which induces the isotopy of the non minimal ones. Let $D^2_{\gamma'}$ denote the isotoped $D^1_{\gamma'}$.  Note that it intersects some $P_j$ in two pie cocores of opposite sign.  Third, if these cocores $\subset \sigma_j$, then $D^2_{\gamma'}$ is isotoped during B) to nest $D^1_{g_j'}$, the standard disc spanning $g_j'$.   Its final position is then determined during C) by the glissade applied to $D^1_{g'}$.  Otherwise, these pie cocores fillings are replaced during C) by copies of $\tilde R^0_{j, 5/8}$ to obtain $D_{\gamma'}$.  By construction $D^1_{\gamma'}$ induces the 0-framing on $\gamma'$.  Isotopy does not change that framing, so $D_\gamma$ also induces the 0-framing.  Similarly the spanning disc $D_{g_j'}$ is first obtained by the anellini compression operation to obtain $D^1_{g'_j}$ which intersects $P_j $ in a single pie cocore.  $D_{g_j'}$ is then obtained by replacing the cocore filling with a copy of $\tilde R^0_{j, 5/8}$.  The previous argument shows that it induces the 0-framing on $g_j'$.   Thus $D_{g_j'}$ nests exactly what is nested by  $\tilde R^0_{j, 5/8}$ as in Data Memo \ref{memo5} b) and is accounted for in Definition \ref{fwcs} 1a) and the first and third paragraphs 
 of 2a). 

By construction no spanning disc of a knot nests in the spanning disc of a linking circle.  The minimal $N_{ij}'$ and $S_{ij}'$ nestings  that arise during B) are as in the first and third paragraphs of Definition \ref{fwcs} 2a).  The nestings created during C) are a possibly proper subset of those described in Definition \ref{fwcs} 2a) first, third and fifth paragraphs, 2b) and 2c), see Data Memo \ref{memo5} d) and e).  This completes of Step 6 and hence the proof of Theorem \ref{carving theorem}. \end{proof}

As noted all steps in the transformation of $\BS_0$ to $\BS_1$ can be achieved through regular homotopy and hence a new proof of the following result.

\begin{theorem}  Any smooth 3-sphere in $S^4$ is regularly homotopic to the standard 3-sphere.  \qed\end{theorem}

\section{Upgrading to an Optimized presentation}

\vskip 10pt
The goal of this section is to prove the following:

\begin{theorem}\label{optimized theorem}  If $\BS$ is a smooth 3-sphere in $S^4$, then $\BS$ has an optimized FWCS-presentation.  \end{theorem}

\begin{construction}\label{adaptive sphere}  With notation as in \ref{winding}, given the $\partial$-germ coinciding  $F|W$ system $(\mG,\mR,\mF, \mW)$ with $\mR$ in AHF form, then let $\BS'$ denote the 3-sphere $\piinv(1/2)\subset V_k.$  By passing to a finite cover we can assume that $\Max_{i,j}\{\Diam(\pi(R_i))$, $\Diam(\pi(w_j))\}< k/100$ for all $R_i\in \mR$ and $w_j\in \mW$.  If $R_i\cap G_j\neq\emptyset$, where $i\in [1,k/100] $ and $j\in [0, -k/100]$ (resp. $i\in [0,-k/100]$ and $j\in [1,k/100])$, then do a single finger move of $\BS'$ into $G_j$.  Let $\hat\BS$ be the result of these finger moves.  We can assume that they were done so that $\hat\BS\cap \mR=\emptyset$.  Let $\hat B_L$ denote $\hat\BS\cap\mG$.\end{construction}

\begin{remarks}  \label{adaptive remarks}  i)  There is an ambient isotopy $F_t$ of $V_k$ such that $F_1(\BS')=\hat \BS$ and $F_1(\mRs)=\mR$.

ii)  After passing to $\tilde V_k, \hat \BS$ will become our $\BS'_{3/8}$ in the proof of Theorem \ref{carving theorem}.\end{remarks}

\begin{figure}
 \centering
\begin{tabular}{ c c }
 $\includegraphics[width=4.0in]{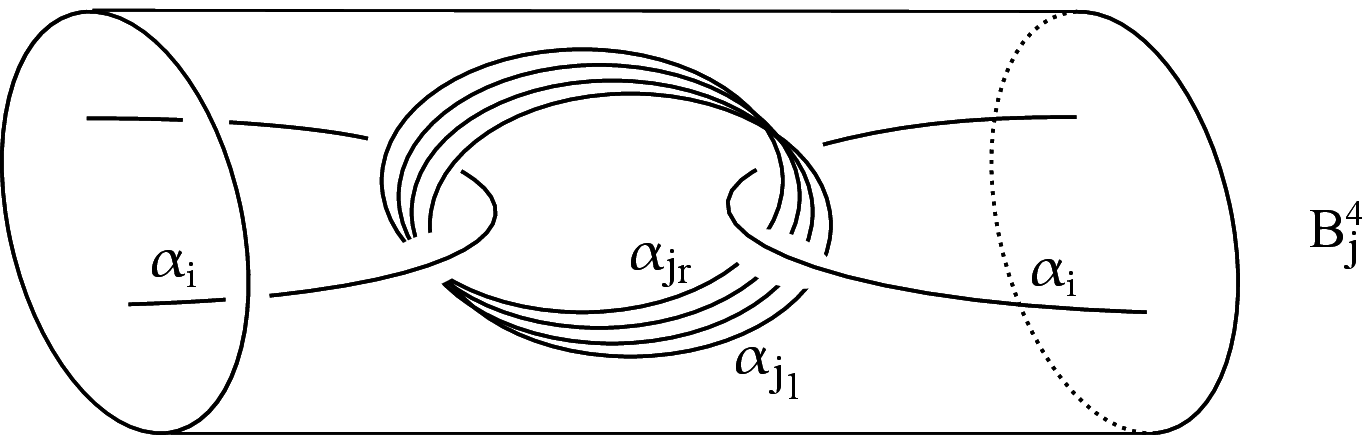}$  
\end{tabular}
 \caption[(a) X; (b) Y]{\label{adaptive1}
 \begin{tabular}[t]{ @{} r @{\ } l @{}}
Adapted Linking\end{tabular}}
\end{figure}

\begin{definition} \label{adaptive definition}  The $F|W$ system $(\mG,\mR, \mF, \mW)$ is $\hat\BS$ \emph{adapted} if $\hat \BS$ arises from Construction \ref{adaptive sphere} and for every $w_i\in \mW, w_i\cap \hat \BS$ is a union of simple closed curves $E_i=\{\alpha_{i_1}, \cdots, \alpha_{i_p}\}$ with $\mE_0=\cup E_i$ such that 

i) For all $j$, level $(\alpha_{i_j})\le 2$, where level is calculated within $w_i$.

ii) $\mE_0=\mA_0\sqcup\mC_0$, where $\mC_0$ contains all the level-2 curves.  If $\alpha\in \mE_0$, then let $D_\alpha\subset w_i$ denote the disc bounded by $\alpha$.  

iii) Here $\mC_0=\{C_0^1, \cdots, C_0^m\}$, where $C_0^j$ is a union of parallel curves $\alpha_{j_{1}}, \cdots, \alpha_{j_{r}}$ that link some $\alpha_{i}\in \mA_0$ as in Figure \ref{adaptive1}.  Furthermore, these linkings occur in disjoint 4-balls $B^4_1, \cdots B^4_m$ that are disjoint from $\mRs$ as well as $B^k_L$ and their spanning discs,  where all $B^4_j\cap D_{j_{1}}, \cdots, B^4_j\cap D_{j_{r}}$ are parallel to each other and to $\hat\BS$ and lie to the same side of $\hat\BS$ opposite that of $B^4_j\cap D_{i}$ and their projections to $\hat\BS$  have parallel clasp intersections.   \end{definition}

\setlength{\tabcolsep}{60pt}

\begin{definition} \label{adaptive plate definition}  Let $(\mG,\mR, \mF, \mW)$ be an $F|W$ with $\mW$ standardly concordant to $\mF$ and with $\mP=\{P_1, \cdots, P_s\}$ the plates.   We say that it is  $\hat\BS$ \emph{plate-adapted} if $\hat \BS$ arises from 
Construction \ref{adaptive sphere} and $P_i\cap \hat \BS$ is a union of simple closed curves $E_i=\{\alpha_{i_1}, \cdots, \alpha_{i_p}\}$ with $\mE_1=\cup E_i$ such that 

i) For all $j$, level $(\alpha_{i_j})\le 2$, where level is calculated within $P_i$.

ii) $\mE_1=\mA_1\sqcup\mC_1$, where $\mC_1$ contains all the level-2 curves.   If $\alpha\subset \mE_1$, then let $P_\alpha\subset P_i$ denote the disc bounded by $\alpha$.  

iii) Here $\mC_1=\{C_1^1, \cdots, C_1^m\}$, where $C_1^j$ is a union of parallel curves $\alpha_{j_{1}}, \cdots, \alpha_{j_{r}}$ that link some $\alpha_{i}\in \mA_1$ as in Figure \ref{adaptive1}.  Furthermore, these linkings occur in disjoint 4-balls $B^4_1, \cdots B^4_m$ that are disjoint from $\mRs$ as well as $B^k_L$ and their spanning discs,  where all $B^4_j\cap P_{j_{1}}, \cdots, B^4_j\cap P_{j_{r}}$ are parallel to each other and to $\hat\BS$ and lie to the same side of $\hat\BS$ opposite that of $B^4_j\cap P_{i}$ and their projections to $\hat\BS$  have parallel clasp intersections.   \end{definition}

  \begin{lemma}\label{adaptive lemma}  1) If $(\mG,\mR,\mF,\mW)$ is $\partial$-germ coinciding and  $\hat\BS$ is as in Construction \ref{adaptive sphere}, then $\mW$ can be isotoped to to be $\hat\mS$ adapted via an isotopy supported away from $\mR\cup\mG$.
  
  2)  If in addition $\mW$ is standardly concordant to $\mF$, then the plates can be isotoped to be $\hat\mS$ plate-adapted via an ambient isotopy that is supported away from the beams, bases and $\mR\cup\mG$.\end{lemma}

\begin{proof}  Here is the idea of 1), with that of  2) being similar.  We attempt to do ``half disc compressions" to $\mW$ along pairwise disjoint embedded half discs whose boundaries consist of one arc in $\mW$ and one in $\hat\BS$ to reduce $\mW\cap\hat \BS$  to a collection of level-1 curves called $\mA_0$.  In practice, among other things the interiors of these half discs may intersect $\mW$.  Isotoping these intersections away will create embedded half discs at the cost of new intersections denoted $\mC_0$ of $\mW$ with $\hat\BS$. See Figure \ref{adaptive2}. More precisely, to prove 1) let $O(\mW)\subset \mW$ consist of $\mW\cap \hat\BS$ together with those points separated from $\partial \mW$ by an odd number of components of $\mW\cap\hat\BS$.  Let $\beta_1, \cdots, \beta_n\subset O(\mW)$ denote pairwise disjoint embedded arcs that cut $O(\mW)$ into discs.  Let $E_1, \cdots, E_n$ be immersed half discs such that $\partial E_i=\beta_r\cup\gamma_r$ where $\gamma_r\subset \hat\BS, \inte(E_i)\cap \hat \BS=\emptyset$ and $\pi(E_i)\subset [-k/20,k/20]$.  These discs might not be suitable because they may a) intersect $\mG$, b) intersect $\mR$, c) intersect each other, d) have the wrong framing, and e) intersect $\mW$ in their interiors.  

\setlength{\tabcolsep}{00pt}
\begin{figure}
 \centering
\begin{tabular}{ c c }
 $\includegraphics[width=6in]{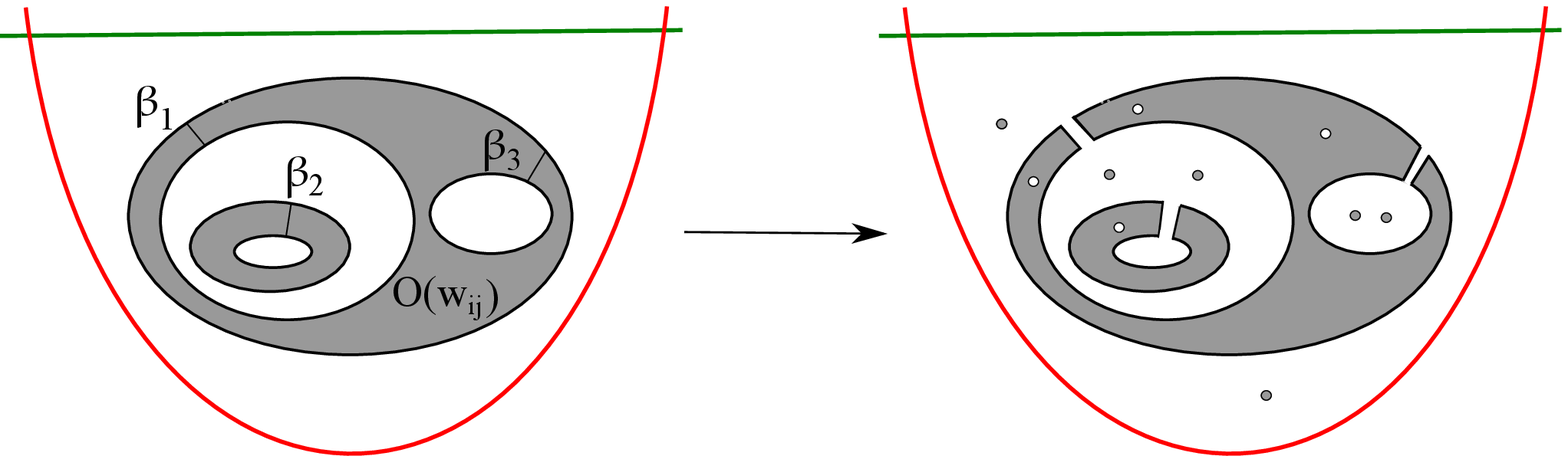}$  
\end{tabular}
 \caption[(a) X; (b) Y]{\label{adaptive2}
 \begin{tabular}[t]{ @{} r @{\ } l @{}}
Boundary Compressing\end{tabular}}
\end{figure}

Since $\pi_1(V_k)\setminus (\hat\BS\cup\mG)$ is freely generated by meridians of the discs in $\mG$ cut off during the finger moves we can
rechoose $\gamma_r$  so that a pushoff of the loop $\gamma_r\cup \beta_r$ is null homotopic in $V_k\setminus (\hat\BS\cup\mG)$.  This  addresses a).  If $x\in E_i\cap R_j$, then let $R'_j$ be a parallel copy of $G_j$ that intersects $\mR$ once and is disjoint from $\hat\BS$.  Tube off a neighborhood of $x\in E_i$ with a copy of $R'_j$ using an arc from $x$ to $R_j\cap R'_j$.  This addresses b) though will likely create new intersections with $\mW$ and the $E_j$'s.   To address c) do finger moves to the $E_i$'s to make them embedded at the cost of creating new intersections with $\mW$.  
To address d) do $\partial$-twisting as in \cite{FQ} (clarified in \cite{E}).  Again the cost is new intersections with $\mW$.  Finally for each $x\in \inte(E_i)\cap \mW$ choose pairwise disjoint arcs $\delta_x\subset E_i\setminus(\beta_i \cup \mRs)$  from $x$ to $\hat \BS$.  Do finger moves to $\mW$ along the $\delta_x$'s and then a bit further to make $(\cup\inte(E_i))\cap \mW=\emptyset$.  Denote the new components of $\mW\cap \hat\BS$ corresponding to $x$ by $\epsilon_x$, let $\mC_0$ be the union of these $\epsilon_x$'s and let $D_x$ denote the disc in $\mW$ bounded by $\epsilon_x$.  Let $\mA_0=(\mW\cap \hat\BS)\setminus \mC_0$.  Note that the $D_x$'s corresponding to a given $E_i$ can be isotoped to be parallel to each other and contained in a 4-ball as in Definition \ref{adaptive definition} iii).

Note that each component of $\mW\cap \hat\BS$ is now level-$\le 2$ and the level-2 ones are among the $\epsilon_x$'s.  Further, if $x$ corresponds to a $E_i$, $w_{rs}$ intersection where $\beta_i\subset w_{pq}$, then $\epsilon_x$ is level-2 if and only if $pr>0$.  For example, if $p<0$, then the half disc lies on the northern side of $\hat\BS$ causing $D_x$ to lie on the southern side and so $\epsilon_x$ will be level-2 if $r<0$ and level-1 otherwise.    This proves 1).
  
\vskip 8 pt

We now assume that $(\mG,\mR,\mF,\mW)$ satisfies the conclusion of Proposition \ref{concordance} and modify  the proof of 1) taking care to avoid the beams and bases.    To start with let $\mP$ denote the union of the plates and then then define $O(\mP)$ analogously as above to $O(\mW)$.  If $P_{ij}$ is a plate from $w_{ij}$, then $O(P_{ij})$ may intersect the beams in a finite set of arcs in $\inte(O(P_{ij}))$.  Note that if  $i<0$ (resp. $i>0$), then these arcs are cores of  beams of  $w_{pq}$'s, with $p>0$ (resp. $p<0$).  

Next choose $\beta_i$'s as above to be disjoint from the beams  and then choose $E_i$'s as above to be disjoint from $\mG$.  These $E_i$'s may intersect the bases and the beams in addition to the issues corresponding to b), c), d), e) above.   In this setting e) corresponds to $E_i$'s intersecting $\mP$ in their interiors.  The intersections with the bases and beams can be eliminated at the cost of new  b), d) and e) intersections.  Next eliminate b) intersections as before, again at the cost of more d) and e) intersections.  Again, the framing can be corrected at the cost of more e) intersections and the type d) intersections can be eliminated  at the cost of more e) intersections.

Finally for each $x\in \inte(E_i)\cap \mP$ choose pairwise disjoint arcs $\delta_x\subset E_i\setminus(\beta_i \cup \mRs\cup\textrm{beams})$  from $x$ to $\hat \BS$.  Do finger moves to $\mP$ along the $\delta_x$'s and then a bit further to make $(\cup\inte(E_i))\cap \mP=\emptyset$.  Denote the new components of $\mP\cap \hat\BS$ corresponding to $x$ by $\epsilon_x$, let $\mC_1$ be the union of these $\epsilon_x$'s and let $D_x$ denote the disc in $\mP$ bounded by $\epsilon_x$.   Let $\mA_1=(\mP\cap \hat\BS)\setminus \mC_1$.  Note that the $D_x$'s corresponding to a given $E_i$ can be isotoped to be parallel to each other and contained in a 4-ball as in Definition \ref{adaptive definition} iii).  Note that each component of $\mP\cap \hat\BS$ is now level-$\le 2$ and the level-2 ones are among the $\epsilon_x$'s.  Again, if $x$ corresponds to an $E_i$, $w_{rs}$ intersection with $\beta_i\subset w_{pq}$, then $\epsilon_x$ is level-2 if and only if $pr>0$.    \end{proof}

\noindent\emph{Proof of Theorem \ref{optimized theorem}}  Assume that $\mW$ satisfies the conclusion of Lemma \ref{adaptive lemma} 2) with $O(\mP), \mA_1, \mC_1$ as in the proof.  For each beam $b$, let $D_b$ denote a cocore disc.  Let $Q$ be a component of $O(\mP)$ and $\partial_e Q$ denote the outermost component of $\partial Q$ and $D_Q\subset \mP$ the disc bounded by $\partial_e Q$.  Note that $Q$ is a planar surface corresponding to two parallel copies $Q', Q'' \subset \mW$. Similarly define $\partial_e Q', \partial_e Q'', D_{Q'} $ and $D_{Q''}$. 

\begin{figure}
 \centering
\begin{tabular}{ c c }
 $\includegraphics[width=5.5in]{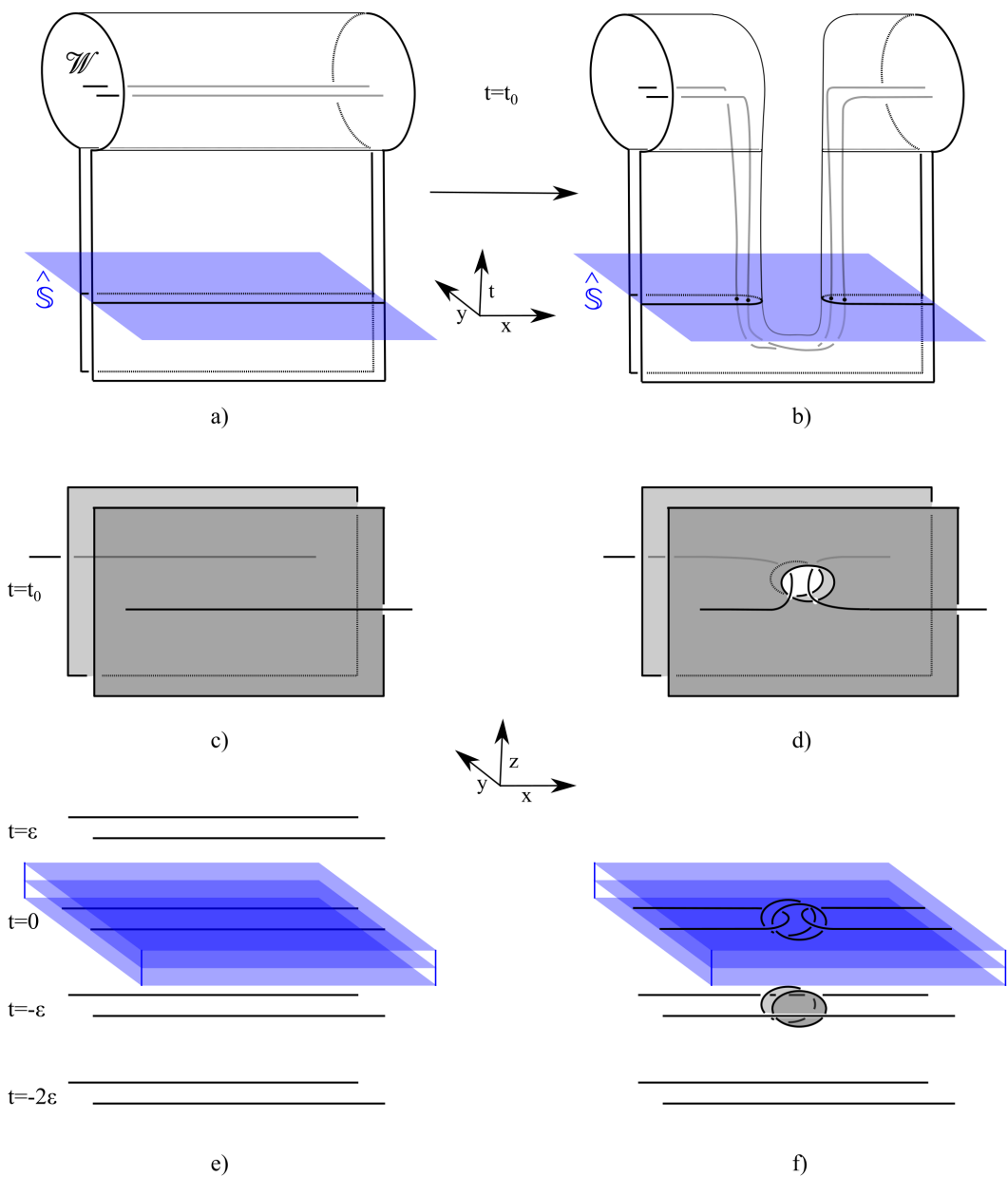}$  
\end{tabular}
 \caption[(a) X; (b) Y]{\label{beam compression}
 \begin{tabular}[t]{ @{} r @{\ } l @{}}
A beam compression\end{tabular}}
\end{figure}

For each $Q$ with $\partial_e Q\in \mA_1$ choose pairwise disjoint paths $\gamma_Q\subset P\setminus\inte(O(P))$ disjoint from $\mRs$  from $Q$ to a parallel copy $D_{b_Q}$ of some $D_b$ that intersect the beams only at their terminal endpoints.      Now for each such $Q$, do boundary compressions of $\mW$ into $\hat\BS$ following $D_{b_Q}\cup \gamma_Q$, simultaneously moving neighborhoods of $\inte D_{b_Q}\cap \mW$ out of the way, thereby creating new intersections of $\mW$ with $\hat\BS$ that are  denoted $\beta_x$'s.  The boundary compression bands together $Q'$ and $Q''$ to obtain $Q^*$ where the band connects $\partial_e Q'$ to $\partial_e Q''$.  Denote by $\partial_e Q^*$ its new outermost boundary and $D_{Q^*}\subset \mW$ its new disc.  See Figure \ref{beam compression}.  Figures a) and b), c) and d), e) and f) show various 3-dimensional slices before and after the boundary compression.  Note that if $\beta_x\subset w_{rs}$ and $Q^*\subset w_{pq}$, then $\beta_x$ is level-2 if and only if $pr<0$.

We now show that the resulting $\mW$ satisfies the conclusion of the lemma.  While strictly speaking the statements that follow are at the $V_k$ level, since $k$ was chosen to be very large they all hold when lifted to $\tilde V_k$ and under the construction in the proof of Theorem \ref{carving theorem} they continue to hold at the $S^4$ level.   In particular, since the description that follows occurs away from $\mG$ we can see this directly using the map $\lambda$ defined in Step 1 of the proof.  Let $\mA$ be the $\partial_e Q^*$'s which by construction are level-1.   Each $\alpha\subset \mC_1$ corresponds to two parallel curves in $\mW$.  The union of these curves is denoted $\mC$ and the union of the $\beta_x$'s is called $\mB$.  By construction the disc in $\mW$ bounded by any element of $ \mB \cup\mC$ has interior disjoint from $\mA\cup\mB\cup\mC\cup \mRs$ and hence conditions i) and iii) of Definition \ref{optimized} hold.  For each component of $\partial_e Q\in \mA_1$, let $Y_Q\subset \hat\BS$ be a small neighborhood of $\partial_e Q$ that contains $\partial_e Q'\cup\partial_e Q''$.    Following the boundary compression the union of these $Y_Q$'s satisfy condition iv).  The $B^4$'s that contain the $\epsilon_x$'s when restricted to $\hat\mS$ correspond to the 1-handles added to $Y_Q$'s to create the $J_i$'s which satisfy condition v).

If $\alpha\in \mB\cup\mC$, then the disc $D_\alpha\subset \mW$ bounded by $\alpha$ lies to one side of $\hat\BS$, indeed is boundary parallel, hence induces the 0-framing on $\alpha$.  We now show that if $\alpha\in \mA$, then $D_\alpha$ induces the 0-framing.  Such an $\alpha$ is of the form $\partial_e Q^*$.  Here $\partial_e Q'\cup\partial_e Q''$ bound a thin annulus $\subset \hat \BS$ with interior disjoint from $\mG\cup \mW\cup\mR\cup \mRs$. Also $\partial_e Q^*=\partial E_{Q^*}\subset \hat\BS$ where the disc $E_{Q^*} $ is the result of  cutting this annulus along an arc.  Since there is a regular homotopy of $D_{Q^*}$ to $E_{Q^*}$ fixing the boundary pointwise, $D_{Q^*}$ induces the 0-framing on $\partial_e Q^*$.   \qed 

\setlength{\tabcolsep}{60pt}

\vskip 10pt

We now summarize, for an optimized FWCS-presentation, the nesting and lodging relations for a given knot and linking circle.  
\begin{notation} We will often denote a knot (resp. linking circle) in $S_{ij}$ labeled by a 0 or $\bullet$ by $\overset\circ S_{ij}$ or $\overset\bullet S_{ij}$ (resp. $\overset\circ {S_{ij}'}$ or $\overset\bullet{S_{ij}'}$).\end{notation}

\begin{figure}
 \centering
\begin{tabular}{ c c }
 $\includegraphics[width=4.5in]{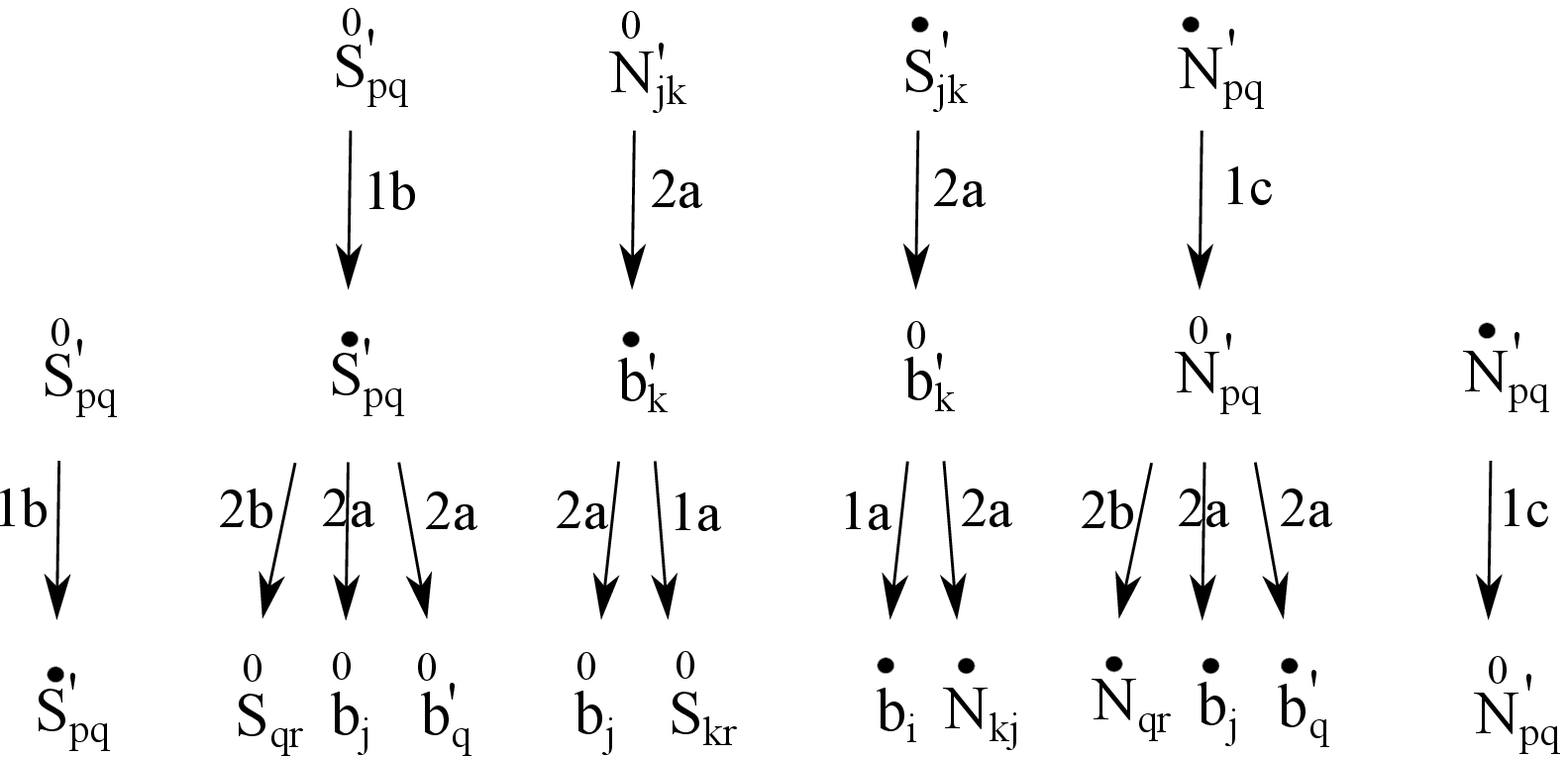}$  
\end{tabular}
 \caption[(a) X; (b) Y]{\label{diagram1}
 \begin{tabular}[t]{ @{} r @{\ } l @{}}
Nesting and lodging of linking circles in an optimized FWCS-presentation\end{tabular}}
\end{figure}

\begin{figure}
 \centering
\begin{tabular}{ c c }
 $\includegraphics[width=5.5in]{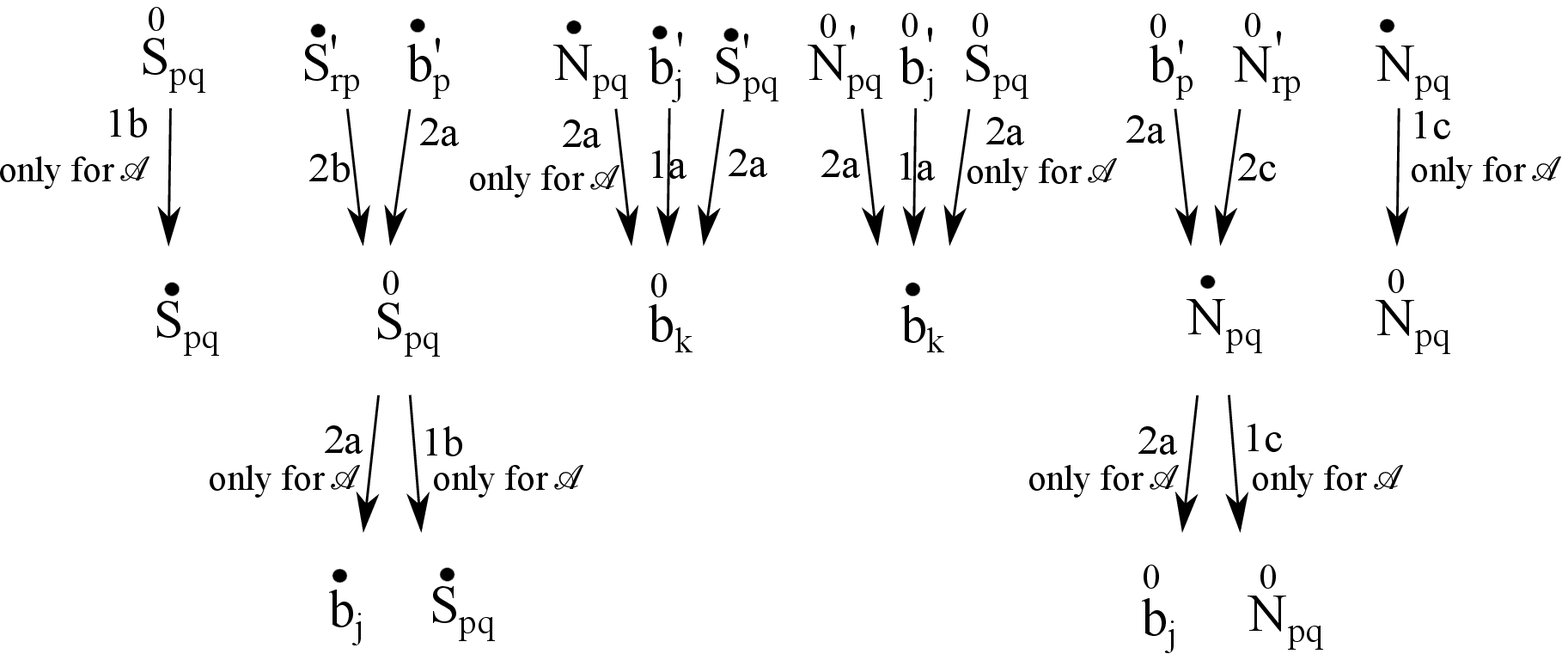}$  
\end{tabular}
 \caption[(a) X; (b) Y]{\label{diagram2}
 \begin{tabular}[t]{ @{} r @{\ } l @{}}
Nesting and lodging of knots in an optimized FWCS-presentation\end{tabular}}
\end{figure}

\begin{lemma} \label{dependency} The possible arrows in and out of a knot or linking circle of L for an optimized FWCS-presentation are stated in Figures \ref{diagram1} and \ref{diagram2}.\end{lemma}

\begin{proof}  These diagrams record the  partial order relations stated in Definition \ref{fwcs} subject to the conditions of Definition \ref{optimized}.\end{proof}

\setlength{\tabcolsep}{20pt}

\setlength{\tabcolsep}{20pt}
 
 \begin{remarks} i)  Being optimized, the knots and linking circles have level-$\le 2$.  It follows that a southern (resp. northern) knot is level-1 if and only if it is labeled with a 0 (resp. $\bullet$) and a southern (resp. northern) linking circle is level-1 if and only if it is labeled with a $\bullet$ (resp. 0).  Thus, there are six types of knots and linking circles.  The possible arrows in and out of each are shown in the diagram.  
 
 ii) Note that if a $b_j$ is labeled with a 0 (resp. $\bullet$),  then $j<0$ (resp. $j>0$) and if a $b_j'$ is labeled with a 0 (resp. $\bullet$), then $j>0$ (resp. $j<0$).  
 
iii)  The arrows indicate possible relations only where they can be defined.  For example a bulleted $S_{pq}'$ can have an arrow to a 0-labeled $S_{qr}$ only when $q<0$ and an arrow to a 0-labeled $b'_q$ only when $q>0$.  

 iv)  The data next to an arrow indicates where in Definition \ref{fwcs} the arrow appears.    \end{remarks}
 
\begin{corollary}  With respect to the partial order of an optimized FWCS-presentation for $A_{pq}'\in S'\cup N'$ and $A_{rs}\in S_L^k\cup  N_L^k$, then $A_{pq}'>A_{rs}$ if and only if $q=r$.\qed\end{corollary}


\section{Operations on CS-Presentations}

We now define a series of geometric operations on a possibly \emph{abstract} (defined below) CS-presentation $(\mL_0, \mD_0)$.    In what follows all the links $\mL$ of the presentation consist of knots $L^k$ and their linking circles $L' $ and all the handles are 0-framed.  We give the definitions for 0-abstract presentations, with the $\bullet$-abstract versions defined analogously.

\begin{operation}  \textbf{Concrete to abstract:}  Given a  CS-presentation of $(\mL_0, \mD_0)$ obtain a new description $(\mL_0, \mD_1)$ of $\Delta_S$ by declaring some set of 0-labeled handles to be \emph{abstract}.  By this we mean they correspond to  0-framed 2-handles attached to $\BS_0$ that are not concretely embedded in $S^4$.    The disc $D$ corresponding to a formerly 0-labeled $\alpha\in \mL$ is replaced by the core $D'$ of the abstract 0-handle.  For other discs $D_i\in \mD_0$, each parallel copy of $D_i$ is replaced by a parallel copy of $D_i'$.  Thus, abstract 0-handles no longer nest other carvings, though they  continue to have carved handles removed as before.  A CS-presentation with abstract handles will be called an \emph{abstract carving/surgery presentation}.  In a similar manner  starting with a CS-presentation we obtain a new description of $\Delta_N$ by declaring some bulleted 2-handles to be abstract.  A 2-handle labeled with a 0 (resp. $\bullet$) will be labeled with an $\oslash$ (resp. $\star$) when made abstract.  To clarify its type, an abstract CS-presentation for $\Delta_S$ (resp. $\Delta_N$) may be called a \emph{0-abstract CS-presentation} (resp. \emph{$\bullet$-abstract CS-presentation}). A non abstract CS-presentation will be called \emph{concrete}.  We naturally generalize the concrete to abstract  operation to allow for $(\mL_0, \mD_0)$ having abstract handles, however if $(\mL_0, \mD_0)$ is 0-abstract (resp. $\bullet$-abstract) only 0-labeled (resp. $\bullet$-labeled) elements of $\mL$ can be  made abstract.\end{operation} 

\begin{remark} \label{abstract remark}  A CS-presentation $(\mL_0, \mD_0)$ gives descriptions of both $\Delta_S$ and $\Delta_N$, however when made 0-abstract (resp. $\bullet$-abstract) we retain a description of $\Delta_S$ (resp. $\Delta_N$), however we lose contact with $\Delta_N$ (resp. $\Delta_S$).\end{remark}

\begin{operation} \textbf{Abstract to concrete:}  Suppose $A'=\{\alpha_1', \cdots, \alpha_m'\} \subset L'$ with all $\alpha_i'$'s labeled with $\oslash$.  Let $A=\{\alpha_1, \cdots, \alpha_m\}$ denote the corresponding knots.  Suppose further

i) There exist knots $B=\{\beta_1, \cdots, \beta_n\}$ such that every arrow into $A$ comes from an element of $B$.

ii) For all $1\le p\le m$ and $1\le q\le n$ we do not have $\beta_q'\ge\alpha_p'$, where $\beta_q'$ is the linking circle to $\beta_q$.  

Then we obtain a 0-abstract CS-presentation $(\mL_0, \mD_1)$ such that each $\alpha_i'\in A'$ is labeled with a 0.  The partial order on $(\mL_0, \mD_1)$ is induced from that of $(\mL_0, \mD_0)$ plus the following additional relations.  There is an arrow from $\alpha_p'\in A'$ to $\beta_q'$ if there is an arrow from $\beta_q\in B$ to $\alpha_p$. The abstract handles of $(\mL_0, \mD_1)$ are exactly those of $(\mL_0, \mD_0)$ with $A'$ removed.

Construct $\mD_1$ as follows.
Let $D_{0p}'\in \mD_0$ be associated to $\alpha_p'$.   
  Let $E_q\subset \mD_0$ (resp. $E_q'\subset \mD_0)$ denote the 2-disc corresponding to $\beta_q$ (resp, $\beta_q'$).  Construct the disc $D_{1p}'\in \mD_1$ for $\alpha_p'$ by starting with the standard 2-disc $A_p'\subset B_N^4$ spanning $\alpha_p'$ and then tubing off each intersection $x\in A_p'\cap E_q$ with a parallel copy of $E_q'$ together with a parallel copy of the annulus from $\partial N(\beta_q)$ to $\partial N(\beta_q')$.  The tube follows a path from $x$ to $\partial E_q$, with the various paths being pairwise disjoint.  The remaining discs of $\mD_1$ are obtained naturally from $\mD_0$.  In particular, those that nested in $D_{0p}$ now nest in $D_{1p}$.  \end{operation}

\begin{operation} \textbf{Knot - linking circle cancellation:}  Here we have $(\mL_0, \mD_0)$ with no arrows pointing into some $k \in L^k$ whose linking circle $k'$ is labeled with an $\oslash$.  Obtain $(\mL_1, \mD_1)$ by letting $\mL_1=\mL_0\setminus \{k, k'\}$.  Let $\hat E\subset B^4_N$  be a cocore of the carved 2-handle to $k$.  We can assume that the abstract 2-handle to $k'$ is attached along $\partial \hat E$ with $E'\in \mD_0$ denoting the disc corresponding to $k'$.  Obtain $\mD_1$ as follows.  If $D_0\in \mD_0$ with $\partial D\in \mL_1$, then obtain the corresponding disc $D_1\in \mD_1$ by replacing each copy of $E'\subset D_1$ with a copy of $\hat E$.  \end{operation}

\vskip 8pt
\noindent
\textbf{Remark.}  Equivalently, this operation is the result of first applying the abstract to concrete operation and then canceling via isotopy, the 2-handle corresponding to $k'$ with the carved 2-handle corresponding to $k$.
\vskip8pt

\begin{operation} \textbf{Hopf link cancellation:} Here $\mL_0$ has the split Hopf sublink  $k$, $k'$ with $k'$ labeled with a $\oslash$, i.e. there is a 3-ball $B^3 \subset \BS_0$ such that $B^3\cap \mL_0=k\cup k'$.  We obtain $\mL_1$ by deleting $k, k'$ from $\mL_0$.  The partial order on $(\mL_1, \mD_1)$ is the restriction of the partial order on $(\mL_0, \mD_0)$, i.e. the arrows in and out of $k$ and $k'$ are deleted.  Here $\mD_1$ is obtained in two steps.  First obtain $\mD_1'$ by replacing neighborhoods of subdiscs of $\mD_0$ that nest the carved $D_k$ by ones $\subset \inte(B^4_N)$ parallel to a standard $\BS_0$ spanning disc of $k$.  Then obtain $\mD_1$ from $\mD_1'$   as in knot - linking circle cancellation.   \end{operation}

\begin{operation} \textbf{a) 0-0 Handle sliding}  This operation involves sliding a 0 or $\oslash$ labeled knot $\alpha$ over a 0 or $\oslash$ labeled linking circle $\beta'$ where $\alpha\nless\beta'$.  The sliding path should not cross the annuli  between knots and linking circles.  Here $\mL_1$ and $\mD_1$ are obtained by the usual 4-manifold sliding operations; however, if $\beta'$ is $\oslash$ labeled, then  $\alpha$ becomes  $\oslash$ labeled and hence abstract.  We record the partial order changes.

i) $(\mL_1,\mD_1)$ has all the arrows of $(\mL_0,\mD_0)$, except that if $\alpha$ becomes abstract its outward arrows are deleted.  

ii) If there is an arrow from $\gamma$ (which may be a knot or linking circle) to $\alpha$, then in $(\mL_1,\mD_1)$ there is an arrow from $\gamma$ to $\beta'$.

iii) If there is an arrow from $\beta'$ to $\delta$ (which may be a knot or linking circle) and $\alpha$ is labeled with a 0, then in $(\mL_1,\mD_1)$ there is an arrow from $\alpha$ to $\delta$. 

 \vskip 8pt

\noindent\textbf{b) $\bullet-\bullet$ Handle sliding}  This operation involves sliding a $\bullet$ or $\star$ labeled knot $\alpha$ over a $\bullet$ or $\star$ labeled linking circle $\beta'$ where $\alpha\nless\beta'$.  The sliding path should not cross the annuli between knots and linking circles.  Here $\mL_1$ and $\mD_1$ are obtained by the usual 4-manifold sliding operations; however, if $\beta'$ is $\star$ labeled, then  $\alpha$ becomes  $\star$ labeled and hence abstract.  We record the partial order changes.

i) $(\mL_1,\mD_1)$ has all the arrows of $(\mL_0,\mD_0)$, except that if $\alpha$ becomes abstract its outward arrows are deleted.  

ii) If there is an arrow from $\gamma$ (which may be a knot or linking circle) to $\alpha$, then in $(\mL_1,\mD_1)$ there is an arrow from $\gamma$ to $\beta'$.

iii) If there is an arrow from $\beta'$ to $\delta$ (which may be a knot or linking circle) and $\alpha$ is labeled with a $\bullet$, then in $(\mL_1,\mD_1)$ there is an arrow from $\alpha$ to $\delta$.

\vskip 8pt

\noindent\textbf{c) $0-\bullet$ Handle sliding}  This operation involves sliding a $0$ or $\oslash$ labeled knot $\alpha$ over a $\bullet$ or $\star$ labeled linking circle $\beta'$ where $\alpha\nless\beta'$ and either $\alpha$ is 0 labeled or $\beta'$ is $\bullet$ labeled.  The sliding path should not cross the annuli between knots and linking circles.  Here $\mL_1$ and $\mD_1$ are obtained by the usual 4-manifold sliding operations.   If $\beta'$ is $\star$-labeled and there is an arrow from $\gamma$ (which is either a knot or linking circle) into $\alpha$, then in $(\mL_1, \mD_1)$, $\gamma$ becomes $\star$ labeled and hence abstract.  We record the partial order changes.

i) $(\mL_1,\mD_1)$ has all the arrows of $(\mL_0,\mD_0)$, except that if $\gamma$ as above becomes abstract its outward arrows are deleted.  

ii) There is an arrow from $\alpha$ to $\beta'$ in  $(\mL_1,\mD_1)$ unless $\alpha$ is $\oslash$ labeled in which case there are no arrows from $\alpha$.

iii) If there is an arrow from $\beta'$ to $\delta$ (which may be a knot or linking circle) and an edge from $\gamma$ to $\alpha$ as above, then in $(\mL_1,\mD_1)$ there is an arrow from $\gamma$ to $\delta$.

\vskip 8pt

\noindent\textbf{d) $\bullet-0$ Handle sliding}  This operation involves sliding a $\bullet$ or $\star$ labeled knot $\alpha$ over a $0$ or $\oslash$ labeled linking circle $\beta'$ where $\alpha\nless\beta'$ and either $\alpha$ is $\bullet$ labeled or $\beta'$ is $0$ labeled.  The sliding path should not cross the annuli between knots and linking circles.  Here $\mL_1$ and $\mD_1$ are obtained by the usual 4-manifold sliding operations.   If $\beta'$ is $\oslash$-labeled and there is an arrow from $\gamma$ (which is either a knot or linking circle) into $\alpha$, then in $(\mL_1, \mD_1)$, $\gamma$ becomes $\oslash$ labeled and hence abstract.  We record the partial order changes.

i) $(\mL_1,\mD_1)$ has all the arrows of $(\mL_0,\mD_0)$, except that if $\gamma$ as above becomes abstract its outward arrows are deleted.  

ii) There is an arrow from $\alpha$ to $\beta'$ in  $(\mL_1,\mD_1)$ unless $\alpha$ is $\star$ labeled in which case there are no arrows from $\alpha$.

iii) If there is an arrow from $\beta'$ to $\delta$ (which may be a knot or linking circle) and an edge from $\gamma$ to $\alpha$ as above, then in $(\mL_1,\mD_1)$ there is an arrow from $\gamma$ to $\delta$.
\end{operation}


\section {Embedding Poincare Balls in $S^4$}

By \emph{Poincare 4-ball} we mean a contractible 4-manifold whose boundary is $S^3$.  A \emph{Schoenflies 4-ball} is a Poincare 4-ball that embeds in $S^4$.
It is well known that the smooth 4-dimensional Poincare conjecture (SPC4) follows from the Schoenflies conjecture and the \emph{Poincare ball embedding conjecture}, i.e. ``every Poincare 4-ball embeds in $S^4$".

The goal of this section is to highlight three approaches towards the embedding conjecture.  One classical and the other two based on the methods of this paper.  

\begin{conjecture} \label{5ball} If $\Delta^4$ is a Poincare ball, then $\Delta^4\times I$ is diffeomorphic to $B^5$ and hence is a Schoenflies ball.  \end{conjecture}

\begin{remarks} i) Since a contractible 5-manifold with boundary $S^4$ is the 5-ball, p. 395 \cite{Sm2}, $\Delta^4\times I=B^5$ is equivalent to the double $D(\Delta^4)$ of $\Delta^4$ being diffeomorphic to $S^4$. Note that $\partial (\Delta^4\times I)=\Delta^4\cup_{\partial}\bar\Delta^4$.

ii) Conjecture \ref{negative equals reverse} is exactly that Conjecture \ref{5ball} holds for Schoenflies balls.  
 \end{remarks}
 
 The following is a special case of the \emph{Gluck conjecture}, that a Gluck twisted $S^4$ is diffeomorphic to $S^4$.

\begin{conjecture} \label{gluck}  If $\Delta^4$ is a Gluck ball, then $\Delta^4\times I=B^5$.\end{conjecture}

A classical approach to Conjecture \ref{5ball}, which we first learned  from Valentin Poenaru, is to prove the following two notorious conjectures.  See 4.89 \cite{Ki2}.

\begin{conjecture} \label{ac} i) If $\Delta^4$ is a Poincare ball, then $\Delta^4\times I$ has a handle decomposition with only 0, 1 and 2-handles.

ii) If $\Delta^5$ is contractible and built from 0, 1 and 2-handles, then $\Delta^5=B^5$.\end{conjecture}

\begin{remarks} A stronger form of i) is that $\Delta^4$ itself has a handle decomposition with only 0,1 and 2-handles.

ii) The Andrews - Curtis conjecture implies Conjecture \ref{ac} ii),  however, the Akbulut - Kirby presentations \cite{AK} are potential Andrews - Curtis counterexamples and all but the simplest ones are conjectured by Gompf \cite{Go} to be pairwise AC inequivalent. On the other hand, Gompf \cite{Go} has shown that $\Delta^5$'s arising from  these presentations are 5-balls.  It would be interesting to have more results in this direction.  A presentation of the trivial group gives rise to a unique 5-manifold with 0,1 and 2-handles inducing that presentation, so we have the following.\end{remarks}

\begin{problem} Show that the Miller-Schupp \cite{MS} presentation $\langle x,y|x^{-1}y^2 x=y^3, x=w\rangle$, where the exponent sum of $x$ in $w$ equals 0, gives $B^5$. \end{problem}

\begin{remark} For other presentations of the trivial group see \cite{Br} and \cite{Li}. \end{remark}

We  now describe a second approach to the Poincare ball embedding conjecture.

\begin{notation} \label{hcobordism notation}  Let $\mH_1=(W_1, B^4, \Delta^4,q_1,v_1)$ be a relative $h$-cobordism between a 4-ball and a Poincare ball with Morse function $q_1$ having only $k$ critical points of index-2 and $k$ of index-3 with gradient like vector field $v_1$ and all the index-2 critical points occur before those of index-3.  Express the middle level  as $B^4\#_k S^2\times S^2$ where $\mG^1=\{G_1, \cdots G_k \}$ are the ascending spheres of the 2-handles with $G_i=S^2\times y_0\subset S^2\times S^2_i$, the i'th $S^2\times S^2$ factor.  Let $\mR^1=\{R_1, \cdots, R_k\}$ denote the descending spheres of the 3-handles.  We can assume that with appropriate orientations $\langle R_i, G_j\rangle=\delta_{ij}$.  Let $\Rs_i$ denote $x_0\times S^2\subset S^2\times S^2_i$.  Let $\alpha$ a nonsingular flow line from $\inte(B^4)$ to  $\inte(\Delta^4)$ and let $B^4_0$, $\Delta^4_0$ respectively denote $B^4\setminus \alpha$, $\Delta^4\setminus \alpha$.\end{notation} 

\begin{definition} \label{stably trivial} The relative $h$-cobordism  $\mH_1$ is \emph{stably trivial} if it can be modified as in i), ii) below to satisfy condition iii).

i)  Let $W_2=W_1\setminus \alpha$.  Construct the the proper h-cobordism $\mH_2:=(W_2, B^4_0, \Delta^4_0, q_2, v_2)$ with $q_2$ and $v_2$ obtained by adding in the standard way an  infinite locally finite  sequence of canceling critical points of index-2 and 3.  The  middle level of $\mH_2$ is $B^4_0\#_\infty S^2\times S^2$,  $\mG^1$ extends to $\mG^2 =\{G_1, G_2, \cdots\}$ where $G_i=S^2\times y_0\subset S^2\times S^2_i$ and $\mR^1$ extends to $\mR^{2}:=\{R_1, R_2,\cdots\}$ where for $j\le k, R_j$ is as before and  for $j>k, R_j=x_0\times S^2\subset S^2\times S^2_j$.  

ii) Modify $v_2$ to $v_3$ so that the new set of descending spheres in the middle level becomes $\mR^{3} :=\{R^3_1, R^3_2, \cdots\}$ where $R^3_j$ is obtained from $R_j$ by applying finitely many finger moves into the $G_i$'s and the totality of finger moves is locally finite. Let $\mG^3$ denote $\mG^{2}$ and $q_3$ denote $q_2$.

iii) There exists a locally finite set $\mW^{3}$ of Whitney discs between $\mR^3$ and $\mG^3$ such that applying Whitney moves to $\mR^3$ using these discs yields $\mR^4$ where $R^{4}_i$ intersects $G^3_j$ geometrically $\delta_{ij}$.  \end{definition}

\begin{theorem} \label{poincare schoenflies} The Poincare ball $\Delta^4$ is a Schoenflies ball if and only if there is a relative $h$-cobordism $\mH=(W,B^4, \Delta^4, q,v)$ which is stably trivial.  \end{theorem}

\begin{proof}  First assume that $\mH$ is stably trivial.  Let $\alpha$ be a nonsingular  flow line of the glvf $v$ which goes from $\inte(B^4)$ to $\inte(\Delta^4)$.  Obtain a proper relative h-cobordism $\mH_1=(W_1, B^4_0, \Delta^4_0, q_1, v_1)$ by restricting to the complement of $\alpha$.  Next  change $(q_1,v_1)$ to $(q_2,v_2)$  to $(q_3, v_3)$  according to i) and ii) and then to $(q_4, v_4)$, where $q_2=q_3=q_4$, to realize the Whitney moves using the discs of iii). Finally, modify to $(q_5, v_5)$ by cancelling the critical points of index-2 and 3.  The resulting $q_5$ and $v_5$ are nonsingular and each flow line of $v_5$ is compact.  Indeed, all the modifications ($q_i, v_i)\to (q_{i+1}, v_{i+1})$ can be chosen to be  locally finite and  supported away from neighborhoods of $B^4_0$ and $\Delta^4_0$.  When $i=1$, the support is in  small neighborhoods of  nonsingular flow lines $\alpha_1, \alpha_2, \cdots$ that  approach $\alpha$.  When $i=2$, the support is in small neighborhoods of the arcs defining the finger moves.  When $i=3$, the support is in a small neighborhood of the Whitney discs and when $i=4$, the support is in a small neighborhood of the flow lines from the j'th index-2 critical point to the j'th index-3 critical point.  The nonsingular vector field $v_5$  induces a diffeomorphism between $B^4_0$ and $\Delta^4_0$.  It follows by Theorem \ref{tame end} that $\Delta^4$ is a Schoenflies ball.

We now prove the converse. By Proposition \ref{equivalence} there exists $\phi\in \pi_0(\Diff_0(\sonesthree))$ that gives rise to the Schoenflies balls $\pm\Delta^4$.  By Corollary \ref{phi to fw}, $\phi=\phi(\mG, \mR, \mF, \mW)$.  By Proposition \ref{germs} after possibly replacing $\phi$ by an S-equivalent one, we can assume that $\mF$ and $\mW$ coincide near their boundaries.   Let $f$ be the pseudo-isotopy from $\id$ to $\phi$ arising from this $F|W$ structure and $v$ the vector field on $S^1\times S^3\times I$ inducing $f$.   Recall that the green and red spheres of this $F|W$ structure arise as the ascending and descending spheres seen in the middle level of a handle structure on $S^1\times S^3\times I$ arising from a Morse function $q_1$ with glvf $v_1$.  This handle structure lifts to one on $\BR\times S^3\times I$.  Do Whitney moves to all the $\tilde\mW$ discs near one end and all the $\tilde\mF$ discs near the other to obtain $(\tilde q_2, \tilde v_2)$.  There should be far separation from the subsets of $\tilde\mW$ and $\tilde \mF$ used.  All but finitely many of the ascending spheres of the 2-handles meet the descending spheres of the 3-handles $\delta_{ij}$.  Cancel all of the $\delta_{ij}$ 2-handles with their corresponding 3-handles to obtain $(\tilde q_3, \tilde v_3)$.  Done appropriately, $\tilde v_3$ coincides with the vector field $\tilde v$ on the end where $\tilde \mW$ was used and with the vertical vector field on the other end.  Also $\tilde q_3$ is the standard projection near that end.  Therefore, $\BR\times S^3\times I$ compactifies to a relative an h-cobordism between $\hat B^4$ and one of $\pm\hat\Delta^4$, where ``hat" denotes ``remove an open 4-ball".  Further, $(\tilde q_3, \tilde v_3)$ extends to $(q_4, v_4)$ on the relative h-cobordism.  Fill in a 4-ball$\times I$ on the end where $\tilde\mF$ was used to obtain a relative h-cobordism between $B^4$ and $\pm\Delta^4$ with $(q_5, v_5)$.  Here $q_5$ is the standard projection on the filled in $B^4\times I$.  Note that had we switched which end to use the $\tilde\mF$ or $\tilde\mW$ discs we would have obtained $\mp\Delta^4$.  Uncompactifying the end corresponding to the $\mF$ Whitney moves i.e.  removing the $\inte(B^4)\times I)$, undoing the handle cancellations on that end and redoing the finger moves achieves i), ii).  The unused discs from $\tilde \mW$ provide the discs needed for iii).\end{proof} 

Here is a third approach using carvings.
\begin{definition}  A 4-manifold $M$ has a \emph{carving/2-handle presentation} if it is obtained from  $B^4$  by first attaching 2-handles to $\partial B^4$ to obtain $B'$ and then carving 2-handles from $B'$.  \end{definition}

\begin{remark}  We can assume that the boundary of the carved 2-handles $\subset \partial B^4$.  We allow the the carved 2-handles to pass through the attaching 2-handles.\end{remark}

\begin{theorem} \label{c2 presentation}  Every Poincare ball has a carving/2-handle presentation.\end{theorem}
\begin{proof} Let $P$ be a Poincare ball.  Since for $k$ sufficiently large $P\#_k \stwostwo$ is diffeomorphic to $B^4\#_k\stwostwo$, it follows that the manifold $P_k$ obtained by attaching 0-framed 2-handles to $k$ split Hopf links on $\partial P$ is diffeomorphic to the manifold $B_k$ obtained by attaching 0-framed 2-handles to $k$ split Hopf links on $\partial B^4$.  Let $\phi:P_k\to B_k$ denote such a diffeomorphism and let $C$ (resp. $C'$) denote the $2k$ cocores (resp. cores) of the first (resp. second) set of 2-handles.    It follows that we can obtain $P$ from $P_k$ by carving the 2k cocores $C$ of the $2k$ 0-framed 2-handles.   Therefore $P$ is obtained from $B^4$ by attaching the 2-handles $C'$ and then carving the 2-handles $\phi(C)$.\end{proof}

\begin{remark} To prove the Poincare ball embedding conjecture it suffices to show that a carving/2-handle presentation of the Poincare ball $P$ can be upgraded to a carving/surgery presentation.  However, even doing this in the simplest nontrivial case would be a great accomplishment.  I.e. where the union of the boundary of the core of the 2-handle and the boundary of the core of the carving form the Hopf link and the 2-handle is +1 framed.  A positive solution implies that Gluck balls embed in the 4-sphere.  Compare with Problem  4.23 of \cite{Ki1} and Question 10.16 \cite{Ga1}. \end{remark}

\enddocument